\RequirePackage{amsthm}

\documentclass[pdflatex,sn-mathphys-num]{sn-jnl}% Math and Physical Sciences Numbered Reference 

\usepackage{graphicx}%
\usepackage{multirow}%
\usepackage{amsmath,amssymb,amsfonts}%
\usepackage{amsthm}%
\usepackage{mathrsfs}%
\usepackage[title]{appendix}%
\usepackage{xcolor}%
\usepackage{textcomp}%
\usepackage{manyfoot}%
\usepackage{booktabs}%
\usepackage{algorithm}%
\usepackage{algorithmicx}%
\usepackage{algpseudocode}%
\usepackage{listings}%
\usepackage{longtable}%
\usepackage{subfigure}%
\usepackage{lmodern}%
\usepackage{anyfontsize}%
\usepackage{float}%

\theoremstyle{thmstyleone}%

\theoremstyle{thmstyletwo}%
\newtheorem{remark}{Remark}%

\theoremstyle{thmstylethree}%

\begin{document}

\title[Article Title]{Numerical Approaches for non-local Transport-Dominated PDE Models with Applications to Biology}

\author[1]{\fnm{Johan} \sur{Marguet}}\email{johan.marguet@univ-fcomte.fr}

\author*[1]{\fnm{Raluca} \sur{Eftimie}}\email{raluca.eftimie@univ-fcomte.fr}
\equalcont{These authors contributed equally to this work.}

\author*[1]{\fnm{Alexei} \sur{Lozinski}}\email{alexei.lozinski@univ-fcomte.fr}
\equalcont{These authors contributed equally to this work.}

\affil[1]{\orgdiv{Laboratoire de Mathématiques de Besançon}, \orgaddress{\street{16 Route de Gray}, \city{Besançon}, \postcode{25000}, \country{France}}}


\abstract{Transport-dominated partial differential equation models have been used extensively over the past two decades to describe various collective migration phenomena in cell biology and ecology. To understand the behaviour of these models (and the biological systems they describe) different analytical and numerical approaches have been used. While the analytical approaches have been discussed by different recent review studies, the numerical approaches are still facing different open problems, and thus are being employed on a rather ad-hoc basis for each developed non-local model. The goal of this review is to summarise the basic ideas behind these transport-dominated non-local models, to discuss the current numerical approaches used to simulate these models, and finally to discuss some open problems related to the applications of these numerical methods, in particular the finite element method. This allows us to emphasize the opportunities offered by this numerical method to advance the research in this field. In addition, we present in detail some numerical schemes that we used to discretize these non-local equations; in particular a new semi-implicit scheme we introduced to stabilize the oscillations obtained with classical schemes.}


\maketitle

\section{Introduction}\label{Introduction}

\begin{figure}[!ht]
    \centering
    \includegraphics[width=\textwidth]{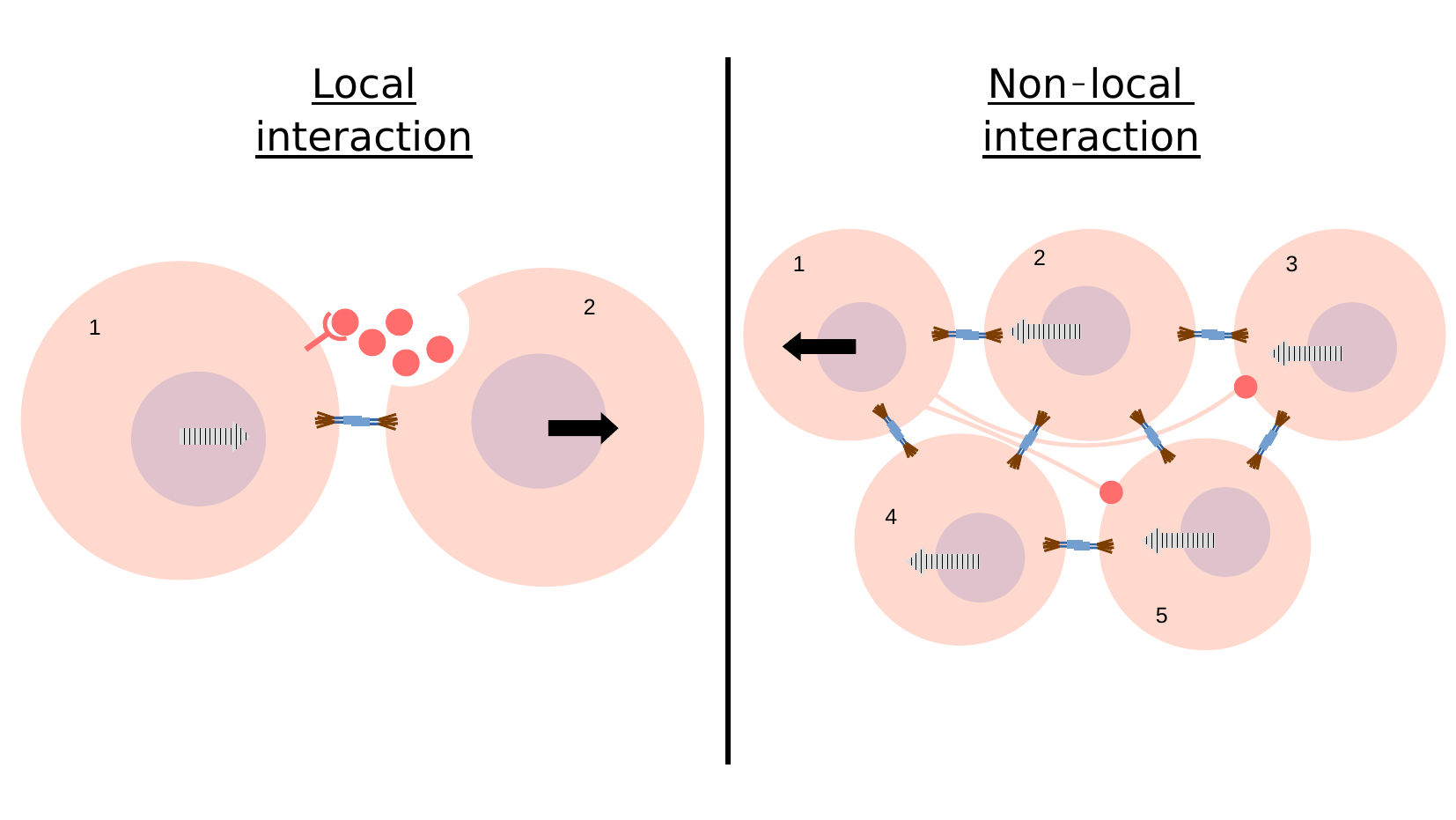}
    \caption{Difference between local and non-local interactions at the cellular level. In both cases, there are two types of interaction: bio-mechanical (represented by the solid and striped arrows) and bio-chemical (represented by bright pink molecules). On the left, we say that these interactions are local since they appear between directly neighboring cells. Whereas on the right, these interactions are non-local since the mechanical traction of cell 1 towards the left leads to the movement towards the left of cells 2 and 4 which are neighbors but also of cells 3 and 5 which are not directly neighboring. In addition, cell 1 can chemically interact with distant cells 3 and 5 via protrusions carrying a chemical messenger.}
    \label{Fig:SchemeInteraction}
\end{figure}

Collective behaviours of cells, bacteria, insects, animals and even humans~\cite{Vicsek_2012} have attracted human attention for at least two millennia~\cite{PlinyTheElder77AD}, due to the spatial and spatio-temporal patterns exhibited by these organisms. Many of these collective behaviours are leading to the directed, coordinated movement of organisms: from the directed migration of cells during health (e.g., tissue morphogenesis~\cite{FriedlGilmour2009,ScarpaMayor2016}) and disease (e.g., cancer~\cite{FriedlGilmour2009}, abnormal wound healing~\cite{adebayo2023mathematical}), to the directed movement of flocks of birds~\cite{hemelrijk2012schools} or herds of ungulates~\cite{yan2023collective} during their annual migrations or swarm of bees~\cite{bhagavan2016structural}. The interactions between the organisms inside these biological/ecological aggregations can be either local (when organisms interact only with immediate neighbours), 
or non-local (when organisms can interact with others further away), see Figure~\ref{Fig:SchemeInteraction}. 
For more than half a century, researchers have acknowledged and further investigated non-local social interactions that are the result of long-distance communication between the different members of these biological/ecological aggregations: from non-local animal-animal interactions due to long-distance auditory, chemical or visual signals~\cite{McComb2003,Payne1971_LongRangeAcousticWhales,Wilson1965_ChemicalCommunicInsects}, to non-local cell-cell interactions via bio-chemical signals (e.g., hormones~\cite{IMCB} and airinemes~\cite{EomDS}) and bio-mechanical forces (e.g., via cellular traction stresses exerted through the substrate~\cite{reinhart2008cell}).  In the context of these experimental research directions, the last two-three decades have witnessed the development of numerous mathematical models aimed at describing such non-local interactions across various cell and animal aggregations: from discrete (individual-based) models, that focus on changes in the position and velocity of individual organisms~\cite{judson1994rise,Mogilner2003_IBMnonlocal,Couzin2002_CollectiveMemory,Reynolds1987_FlocksHerdsSchools}, to continuum (integro-partial differential equation) models that focus on changes on the densities of populations of organisms~\cite{armstrong2006continuum,mogilner1999non}. While the discrete models have been investigated mainly computationally ~\cite{Couzin2002_CollectiveMemory,PinedaWeijerEftimie2015,Vicsek_2012} (due to a lack of methodological approaches for their analytical investigation), the continuum models have been investigated both analytically and computationally~\cite{topaz2006nonlocal,mogilner1999non,armstrong2006continuum}. However, the details of the numerical methods used to discretize these continuum models were not always available, since many researchers have been using classical approaches (i.e., finite difference and finite volume schemes combined with standard quadrature methods for the computation of integrals), which can be found in classical numerical analysis books~\cite{quarteroni2009numerical}. Nevertheless, this lack of details might impede some (early-career) researchers to further develop and investigate numerically various mathematical models that incorporate long-distance biological/ecological interactions. For this reason, in this study we aim to summarise different numerical approaches that could be used to simulate the spatial and spatio-temporal dynamics exhibited by this class of non-local PDE models in 1D and 2D. 
\\ 

\noindent\textbf{Terminology.} Before presenting briefly the various types of non-local mathematical models that have been developed over the years, we need to clarify the terminology used throughout this study. Random cell/bacterial/animal movement is described by \emph{diffusion terms}. Directed cell/bacterial/animal movement in response to a stimulus (e.g., food, stiffer substrates) is described by \emph{taxis terms}. In the literature there are different types of taxis: \emph{chemotaxis} (movement aligned with a chemical gradient), \emph{haptotaxis} (cell movement up a gradient of substrate-bound molecules), \emph{durotaxis} (cell movement up a gradient of extracellular matrix stiffness)~\cite{FortunatoSunyer2022}. All these taxis behaviors are modelled by similar terms: $\nabla \cdot (\chi u \nabla c)$, where $c$ describes the concentration/density of chemicals/molecules/etc. that determine the gradient $\nabla c$, $u$ is the density of individuals/cells/animals that move towards/away the chemical gradient $\nabla c$, and $\chi$ is the chemotactic sensitivity coefficient (that might be a constant or might depend on the concentration of the chemical $c$ or on other factors in the environment). Due to this directional transport of individuals/cells/animals with density $u$, %(transport described by the gradient $\nabla u$\AL{??? isn't transport with velocity $\vec{v}$ described by $\nabla \cdot (u \vec{v})$?})
these \emph{taxis models} are also referred to as \emph{advection models} or \emph{transport models}. While this terminology is common for local models~\cite{gerisch2006robust, berezovskaya1999bifurcations}, it has been used also to describe non-local fluxes that are the result of social interactions (e.g., attractive movement towards higher densities of conspecifics located further away), as those modeled in equation~\eqref{advection_part}~\cite{painter2023biological}. Nevertheless, some studies, mainly applied to cell biology, refer to these non-local models also as \emph{adhesion models}, since the non-local operator describes cell-cell or cell-matrix adhesion~\cite{painter2023biological, hillen2021non}.
Other studies, mainly applied to ecology, refer to these non-local models as \emph{aggregation models}~\cite{fabreges2019numerical, fellner2020solutions, james2015numerical, carrillo2021second,bernoff2013nonlocal}.
\\\\
The non-local interactions between different members of the same social community (e.g., a cell population or an animal population), which are located at various spatial positions (some closer, some further away) are described by integrals over different spatial ranges. These integrals can be included into the mathematical equations describing the changes in the cell/animal population densities through various terms that have different biological meanings: 
\begin{itemize}
\item \textbf{Non-local interactions in the reaction part}~\cite{lorz2013populational}, to describe the information received by an organism about its possibility of proliferation/elimination within a spatial range around its current position. The general structure of these models is:
\begin{equation}
\frac{\partial u}{\partial t} = \nabla \cdot \left( D \nabla u(t,\mathbf{x}) \right) + a[u],
\label{reaction_part}
\end{equation}
where $u(\textbf{x},t)$ is the density of cells/organisms at position $\textbf{x}\in \mathbb{R}^{n}$, $n\in\{1,2,3\}$, and time $t$, $D$ is the diffusion coefficient (which can be constant or density-dependent), $a[u]$ is the non-local source term often in the form $a[u] = F(u) + I[g(u)]$. Here, $F(u)$ describes local factors that impact proliferation/elimination and $I[g(u)] = \int_{\Omega} g(x,u(y,t)) \ dy$ describes non-local factors that impact proliferation/elimination (with $g$ a linear or non-linear function for the elimination/proliferation of population $u$, and  $\Omega$ the domain of interest for these non-local interactions, which might not be the whole $\mathbb{R}^{n}$). 
%%%
\item \textbf{Non-local interactions in the advection part}~\cite{bitsouni2017mathematical,painter2023biological,bitsouni2018non}, to describe the directional flux generated by interactions with conspecifics positioned further away from a reference individual. The terminology used for this class of non-local models, that have more-or-less the same structure, is quite varied: from \emph{non-local advection models} for non-local taxis (e.g., haptotaxis or chemotaxis) ~\cite{krasnianski2019nonlocal}, to \emph{aggregation models}~\cite{fellner2020solutions,james2015numerical,delarue2020convergence,bailo2020convergence}. The general form of these models is:
\begin{equation}
    \frac{\partial u}{\partial t} = \nabla \cdot \left( D \nabla u(t,\mathbf{x}) - \beta u(t,\mathbf{x})K[u] \right) + F(u),
    \label{advection_part}
\end{equation}
where $u(\mathbf{x},t)$ is the density of cells/organisms at position $\mathbf{x} \in \mathbb{R}^{n}$, $n\in\{1,2,3\}$, and time $t$, $D$ is the diffusion coefficient (which can be constant or density-dependent~\cite{krasnianski2019nonlocal}), $\beta$ is the advection coefficient (which can also be constant or density-dependent), $K[u]$ is the non-local operator and $F(u)$ is the reaction term.
In these models, the non-local operator $K[u]$ describes the relative force of interactions exerted by a cell/organism at position $\mathbf{y}$ on a cell/organism at position $\mathbf{x}$.
\\
There is also a second class of models that include non-local terms in the advection part: not directly through the convolution operator (as seen in equation (\ref{advection_part})), but rather through the gradient of the convolution~\cite{fabreges2019numerical}:%(see also the discussion at the beginning of Section~\ref{Sect:DiscussBioandEco})
\begin{equation}
\frac{\partial u}{\partial t} = \nabla \cdot \left( D \nabla u(t,\mathbf{x}) - \beta u(t,\mathbf{x})\nabla \mathcal{M}[u] \right) + F(u) .
    \label{advection_part_aggregation}
\end{equation}
Unlike the previous model, in equation (\ref{advection_part_aggregation}) the interaction potential $\mathcal{M}[u]$ represents a non-local measure of the cell population,  the gradient of which influences the movement of the population~\cite{painter2023biological}. The other terms and parameters in equation (\ref{advection_part_aggregation}) have the same meaning as in equation~\eqref{advection_part}.\\
\begin{remark}
    Note that for some specific cases, models~\eqref{advection_part} and ~\eqref{advection_part_aggregation} can turn out to be equivalent. For example, if in~\eqref{advection_part} we choose the following 1D kernel
    \begin{equation*}
        K[u] = \hat{\omega} * u \text{ with } \hat{\omega}(x)=\left\{
        \begin{array}{cl}
        1, & \quad \text{ if } 0 < x < r,\\
        -1, & \quad \text{ if } -r < x < 0,\\
        0, &\quad \text{else,}\\ 
        \end{array} \right.
    \end{equation*}
     and if in~\eqref{advection_part_aggregation} we choose the 1D kernel
     \begin{equation*}
         \mathcal{M}[u] = \tilde{\omega} * u \text{ with } \tilde{\omega}(x)=\left\{
        \begin{array}{cl}
        |x|, & \quad \text{ if } -r < x < r,\\
        r, &\quad \text{else,}\\ 
        \end{array} \right..
     \end{equation*}
     then we observe that we have $\partial_x \mathcal{M}[u] = K[u].$ Thus, models~\eqref{advection_part} and~\eqref{advection_part_aggregation} are equivalent.
\end{remark}
\end{itemize}

\begin{remark}
We need to mention that there are also other classes of models that consider non-local diffusion processes that result from particles/individuals jumping in/out of neighbouring positions~\cite{othmer1988models}. An example of such a model is given in~\cite{Chang2023_NonlocalDiffusionDengue}:
\begin{equation}
\frac{\partial u(t,\mathbf{x})}{\partial t} = \left(\int_{\mathbb{R}^{n}}\mathcal{J}(\mathbf{x}-\mathbf{y})u(t,\mathbf{y})d\mathbf{y} -u(t,\mathbf{x})\right)+ F(u),
\end{equation}
where $\mathcal{J}(\mathbf{x}-\mathbf{y})$ describes the probability distribution of jumping from position $\mathbf{y}\in\mathbb{R}^{n}$ to position $\mathbf{x}\in\mathbb{R}^{n}$~\cite{Chang2023_NonlocalDiffusionDengue}. Hence, the convolution $\mathcal{J}*u=\int_{\mathbb{R}^{n}}\mathcal{J}(\mathbf{x}-\mathbf{y})u(t,\mathbf{y})d\mathbf{y}$ gives the rate at which individuals arrive at position $\mathbf{x}$ from all other positions $\mathbf{y}$, while $-u(t,\mathbf{x})$ gives the rate at which individuals are leaving position $\mathbf{x}$. \\
Note that these non-local terms that describe jumping between different spatial locations could be interpreted as models with non-local interactions incorporated into the source term, as we can write $a[u]=\mathcal{J}*u+F(u)$.
\end{remark}

In this review we focus only on models with non-local terms inside the advection component since, from a numerical perspective, they pose particular challenges due to the advection operator combined with the particularities of biologically-inspired kernels (which can have singularities if one assumes very strong repulsion when $\textbf{y}=\textbf{x}$, to avoid cells/organisms occupying the same spatial position).
As mentioned before, since the published studies are usually lacking detailed discussions about the various numerical approaches used to discretise and simulate these non-local models, the goal of this review is to fill in this gap and to summarise these different numerical approaches while discussing their advantages and drawbacks, and emphasising the open numerical questions in the field. This way, we hope to provide a more structured approach to choosing appropriate numerical schemes for those non-local models, where the non-local interactions lead to an advective flux.

The rest of this review is structured as follows: in Section \ref{Sect:AnalyticalApproaches} we summarise briefly some of the analytical approaches used to investigate these non-local models (we do not focus more on this aspect, since a very detailed review of this topic was recently published ~\cite{chen2020mathematical}). In Section~\ref{Sect:NumericalApproaches} we summarise the numerical approaches (mainly based on finite differences and finite volume schemes) used to discretise these models. Since the finite element (FE) approaches are almost nonexistent, even if they are well-suited for 2D and 3D models that are important to describe the dynamics of tissues, in Section~\ref{Sect:FiniteElem} we will summarise the finite element method, and present a series of simulations for the dynamics of some very simple 1D and 2D models. In Section~\ref{Sect:OpenProblems} we discuss some open problems related to numerical approaches for this class of non-local models. We conclude this review with Section~\ref{Sect:Discussion}, where we present a final overview of these numerical methods, and mention some future research perspectives triggered by the recent advances on non-local models with applications to biology and medicine.

%%%%%%%%%%%%%%%%%%%%%%%%%%%%%%%%%%%%%%%%%%%%%%%%%%%%%%%%%%%%%%%%%%%%%%%%%%%%%%%%%%%%%%%%%%%%%%%%%%%%%
\section{Brief discussion of non-local transport-dominated models in biology and ecology}
\label{Sect:DiscussBioandEco}

In this section we discuss in more detail the biological background behind some simple non-local models developed in the context of cell biology (see Section~\ref{Subsect:CellBiology}) and ecology (see Section~\ref{Subsect:Ecology}). We also show the general form of the equations derived to describe these biological and ecological interactions. By presenting them in parallel, we aim to emphasise the similarities between the models used in these two different domains. While these aspects were discussed in more detail in other reviews~\cite{painter2023biological,chen2020mathematical,wang2023open}, it is important to present them here as well (but very briefly), for the completeness of the discussion. Especially since these details are not always known to researchers in the numerical analysis field.

%%%%%%%%%%%%%%%%%%%%%%%%%%%%%%%%%%%%%%%%%%%%%
\subsection{Non-local models for cell biology}
\label{Subsect:CellBiology}

In cell biology, cell-cell and cell-matrix adhesion are crucial processes in the organization and maintenance of tissues during development (e.g., morphogenesis) and during pathological situations (e.g., cancer, wound healing)~\cite{ScarpaMayor2016, friedl2012classifying}.
While these processes occur locally via adhesion molecules such as integrins and cadherins \cite{ren2011adhesion}, their effects are perceived not only locally but also at a distance (see also Figure~\ref{Fig:SchemeInteraction}). In general, cell-cell interactions can occur over different spatial ranges:
\begin{itemize}
    \item short-range repulsion: cells repel each other to avoid overcrowding~\cite{burger2007aggregation};
    \item medium-range alignment: cells tend to align due to traction forces ~\cite{trepat2009physical}, or due to substrate topography~\cite{LeclechBarakat2021}, which can align the actin fibres in the cell cytoskeleton;
    \item long-range attraction: cells can detect and interact with other cells positioned at large distances due to adhesion forces exerted by those cells, or through the generation of cell protrusions~\cite{yamashita2018specialized}.
\end{itemize}
The need to consider long-distance cell-cell interactions to understand collective cell dynamics is supported by various experimental studies.
For example, a 2009 study \cite{trepat2009physical} showed that traction forces are acting many cell rows behind the leading edge of the invading cell sheet (as they were measured up to $200\mu m$, equivalent to roughly 10 rows of cells, each with an average diameter of $20\mu m$).
In addition, cells can extend membrane protrusions over distances of up to ten times their diameter or even more, as is the case with cytonemes (that can reach up to $700\mu m$ in length \cite{yamashita2018specialized}), tunnelling nanotubes and tumour microtubes (that can reach up to $100 \mu m$ in length \cite{cancers12040857}). These membrane protrusions are used in long-distance cell-cell communication and movement coordination, via physical contact between cells positioned far away from each other \cite{cancers12040857} (see also Figure~\ref{Fig:SchemeInteraction}).

All these (and many other) biological experiments are emphasising the importance of incorporating non-local cell-cell and cell-ECM interactions that generate different force distributions, into the mathematical models for collective cell movement. To exemplify such a simple model, let us now define $u(\mathbf{x},t)$ to be the density of a cell population at time $t$ and position $\mathbf{x}\in \mathbb{R}^n$ with $n \in \{1,2,3\}$. The evolution of this density can be given by the following generic model:
\begin{equation*}
    \partial_t u = - \nabla \cdot J(\mathbf{x},t),
    \label{cellular density}
\end{equation*}
where $J$ is the population flux resulting from the movement of cells. This flux can be decomposed in two parts $J = J_d + J_a $, where $J_{d}$ is the diffusive part and $J_{a}$ is the advective part: 
\begin{itemize}
\item The diffusive part describes the local random movement of cells, which causes them to spread out in space following Fick's law~\cite{fick}. This spreading is described by the gradient of the density $u$, multiplied by a diffusion coefficient $D$ (here assumed to be constant):
\begin{equation}
    J_{d} = -D \nabla u.
    \label{diffusive flux}
\end{equation}
\item For the advective part, if we assume that a cell located at position $\mathbf{x}$ interacts with another cell located at position $\mathbf{y}\neq \mathbf{x}$, the movement generated by this interaction is the result of the traction forces caused by cell-cell adhesive bonds. Moreover, these adhesive bonds lead to a sequential  transmission of traction forces over larger spatial distances (see also Figure~\ref{Fig:SchemeInteraction}), and thus cells can sense the biomechanical effect of other cells at a distance. Therefore, we can assume that the movement of a generic cell located at $\mathbf{x}$ results from the sum of all interactions between the cell positioned at $\mathbf{x}$ with all other cells located within its perception radius $r$. Mathematically, this translates into the following advective flux:
\begin{equation}
    J_{a} = \ u(\mathbf{x},t) \int_{B_r} \alpha(\mathbf{y}) \mathbf{K_r}(\mathbf{x},\mathbf{y}) f(u(\mathbf{y},t)) \ d\mathbf{y},
    \label{advective flux}
\end{equation}
where $\alpha$ is the adhesion coefficient which could vary depending on the distance separating two cells. The product of $\mathbf{K_r}$ (a vector valued kernel) with $f$ (a scalar valued function) describes the induced direction of the cell movement and its magnitude. %\AL{ Shouldn't the  vector $K_r$ be in bold, as other vectors? }
\end{itemize}
Hence, the generic model can take the following form~\cite{armstrong2006continuum}:
\begin{equation}
    \partial_t u = D \Delta u -  \nabla \cdot \biggl( \ u(\mathbf{x},t) \int_{B_r} \alpha(\mathbf{y}) \mathbf{K_r}(\mathbf{x},\mathbf{y}) f(u(\mathbf{y},t)) \ d\mathbf{y} \biggr) + F(u),
    \label{nonlocal biology}
\end{equation}
where function $F$ was added here to describe the proliferation and decay of cells, to ensure a more general although still simple model. As it was shown in ~\cite{painter2023biological,armstrong2006continuum}, the coefficient $\alpha$ has a crucial importance for cellular aggregation: if this $\alpha$ coefficient is below a certain threshold, then a dispersed cell population remains dispersed over time; whereas if $\alpha$ is above this threshold, the population could aggregate.

Furthermore, when modelling the interactions between different populations of cells inside the same tissue, one needs to consider systems of coupled equations (where coupling can be done through local and/or non-local terms). Thus, model~\eqref{nonlocal biology} can be generalised to describe the dynamics of $p\in\mathbb{N}^{+}$ cell populations:
\begin{align}
    \partial_t u_i =& D_i \Delta u_i - \sum_{j=1}^{p} \nabla \cdot \biggl( u_i(\mathbf{x},t) \int_{B_{r_i}} \alpha_{i,j}(\mathbf{y}) \mathbf{K_{i,j}}(\mathbf{x},\mathbf{y}) f_{i,j}(u(\mathbf{y},t)) \ d\mathbf{y} \biggr)\nonumber\\ 
    & +F(u_{1}, ..., u_{i}, ..., u_{p}),
    \label{nonlocal system biology}
\end{align}
where $D_i$ is the diffusion coefficient for the $i^{th}$ population, $\alpha_{i,j}$ is the adhesion coefficient between the $i^{th}$ and $j^{th}$ population, $r_i$ the sensing radius of $i^{th}$ cell population, $\mathbf{K_{i,j}}$ and $f_{i,j}$ describe the induced direction of movement and its magnitude. Moreover, function $F(u_{1}, ..., u_{i}, ...u_{p})$ describes the proliferation and decay of these cells (sometimes as a result of interactions between cells $u_{i}$ and $u_{j}$). For different examples of such models~\eqref{nonlocal system biology}, that incorporate different functions (and interaction kernels), we refer the reader to~\cite{painter2023biological}. 

Before concluding this discussion on models for cell movement we need to emphasise that in many cases, these non-local models for cell-cell interactions are coupled with ODEs describing, for example, the remodelling of the extracellular matrix (ECM). The ECM does not move in space; it can only be produced and degraded by different cells. Therefore, these models for cell-ECM interactions take the form:
\begin{subequations}
\begin{align}
    \frac{\partial u_i}{\partial t} =& D_i \Delta u_i - \sum_{j=1}^{p}\nabla \cdot \biggl( u_i(\mathbf{x},t) \int_{B_{r_i}} \alpha_{i,j}(\mathbf{y}) \mathbf{K_{i,j}}(\mathbf{x},\mathbf{y}) f_{i,j}(u(\mathbf{y},t);e(t)) \ d\mathbf{y} \biggr)\nonumber \\
    & +F(u_{1}, ...,u_{i}, ...,u_{p}), \\
    \frac{de}{dt}=& p_{e}(u_{i}(\mathbf{x},t))-d_{e}(u_{i}(\mathbf{x},t))e(t).
    \end{align}
    \label{nonlocal system biology wth EDO}
\end{subequations}
In the equations above, the function $f_{i,j}$ can depend also on the density of the ECM $e$; and $p_{e}$ and $d_{e}$ functions describe the secretion and degradation of the ECM (by different cells in the system). All other terms are as described for model~\eqref{nonlocal system biology}.

 %%%%%%%%%%%%%%%%%%%%%%%%%%%%%%%%%%%%%%%%%%%%%%
\subsection{Non-local models for ecology}
\label{Subsect:Ecology}
Animals, like cells, can move collectively and form spatial and spatio-temporal social aggregations. The variety of these aggregations (e.g., swarms of insects, flocks of birds, schools of fish or herds of ungulates) has attracted the attention of both biologists and mathematicians who tried to identify common mechanisms behind them. Different empirical and mathematical studies have suggested that the three types of social interactions (i.e., repulsion, attraction and alignment) that are involved in the formation of these aggregations, might act on different spatial ranges~\cite{Vicsek_2012, article3}: 
\begin{itemize}
    \item at very short distances, animals move away from each other, to avoid physical contact~\cite{cavagna&al};
    \item at a slightly greater distance, animals tend to align their movements (e.g., walking in single file)~\cite{Baltzersen_2022};
    \item over large distances, animals tend to approach each other (if they are too far apart) in order to maintain group cohesion~\cite{larkin2008,cavagna2010scale,cavagna&al}.
\end{itemize}
The spatial distances over which these interactions occur depend on the communication ranges of animals (similar to cells' sensing radius discussed in the previous section). For example, migrating birds can coordinate their movement direction and speed across distances of more than 250 meters~\cite{larkin2008}, while swarming insects coordinate their movements within 200 mm~\cite{Gorbonos2016_LongRangeInsects}.
Moreover, such aggregations can persist for short periods of time (as for swarms of mosquitoes \cite{bengtsson2008aggregation}) or for longer periods of time, such as a season or an entire lifetime (as for wolves~\cite{mech2019wolves}). For example, over a wider time window, we notice that animals do not generally use all the available space but confine themselves to an area called home range, which they keep for a season or for their entire life~\cite{burt1943territoriality}. The presence of these home ranges for different animal species is associated with species-specific processes such as memory (as some animals tend to return to the same area after their migration or hibernation)~\cite{riotte2015memory}, or stigmergy (as some animals leave physical traces, such as scent marks, when they move across the landscape)~\cite{theraulaz1999brief}. 
Therefore, the simple mathematical models that we will review next to describe the spatio-temporal dynamics of various ecological aggregations, include not only local/non-local animal-animal interactions but also local/non-local animal-environment interactions.

To present these models for animal movement, let us define $u(\mathbf{x},t)$ to be the density of an animal population at time $t$ and position $\mathbf{x} \in \mathbb{R}^n$ with $n \in \{1,2,3\}$. The spatio-temporal changes in this density can be described by the following generic non-local reaction-diffusion-advection model~\cite{painter2023biological}:
\begin{equation}
    \partial_t u = D \Delta u - \nabla \cdot [ u \nabla \nu (\omega_r * f(u)) ] + F(u),\label{model:ecology}
\end{equation}
with $D$ the constant diffusion coefficient, $\nu$ the advection coefficient, the product between $\omega_r$ (scalar function) and $f$ having the same meaning as for the taxis-like model~\eqref{nonlocal biology}, and the function $F$ representing the births and deaths of the animal population.\\ 
Furthermore, it is quite common to couple the non-local equation with an ODE representing the process of creating memories and/or traces of passages, and their destruction~\cite{wang2023open}:
\begin{subequations}
    \begin{align}
        &\partial_t u = D \Delta u - \nabla \cdot [ u \nabla \nu (\omega_r * f(u;m)) ] + F(u),\\
        &\partial_t m = \gamma u - \delta m,
    \end{align}
    \label{ecology}
\end{subequations}
with $m(\mathbf{x},t)$ the cognitive map, and the constants $\gamma$ and $\delta$ describing the production/destruction of these memories. Model~\eqref{ecology} can also be easily adapted to describe interactions between multiple animal species (as it was done for model~\eqref{nonlocal system biology}):
\begin{subequations}
    \begin{align}
        \partial_t u_i = &D_i \Delta u_i - \nabla \cdot [ u_i \nabla \sum_{j=0}^p \nu_{i,j} (\omega_r * f(u_i;m_i) )] \nonumber\\
        &+ F(u_1,...,u_i,...,u_p),\\
        \partial_t m_i = &\gamma u_i - \delta m_i.
    \end{align}
    \label{ecology diff species}
\end{subequations}
\\
Throughout this Section~\ref{Subsect:CellBiology} we summarised the most general form of models used to describe nonlocal interactions in cell biology and in ecology, to emphasise the similarities of these models. However, we need to emphasise here also an important difference, which plays a particular role in numerical simulations: the spatial domain. While the ecological domain is relatively simple and can be described by a rectangular domain, the biological domain can take very complex geometrical shapes as cells form tissues and organs. The numerical methods derived to simulate the behaviour of these models should be able to accommodate these complex domains (for more on this see Section~\ref{Sect:ComplexGeom}).
Next, we summarise briefly some of the analytical approaches used to investigate these models (Section~\ref{Sect:AnalyticalApproaches}) followed by a more in-depth discussion of numerical approaches that can be used to simulate the dynamics of these models (Section~\ref{Sect:NumericalApproaches})

%%%%%%%%%%%%%%%%%%%%%%%%%%%%%%%%%%%%%%%%%%%%%%%%%%%%%%%%%%%%%%%%%%%%%%%%%%%%%%%%%%%%%%%%%%%%%%%%%%%%%
\section{Analytical approaches for non-local models}
\label{Sect:AnalyticalApproaches}

Since the overall goal of this article is to summarize the numerical approaches, in this Section we only mention (without going into details) certain existing analytical results for non-local models, where the non-locality appears in the advective part (models~\eqref{advection_part} and ~\eqref{advection_part_aggregation}).
For a recent detailed review of different analytical approaches applied to study such non-local models (with non-localities either in the reaction, diffusion or advection terms), we refer the reader to ~\cite{chen2020mathematical,painter2023biological}.

Many analytical studies focus on the existence and uniqueness of bounded solutions displayed by these non-local models~\eqref{advection_part} and ~\eqref{advection_part_aggregation}: from the existence (under certain hypotheses on the initial condition $u_0$ or even on the source function $f$) of positive weak and/or global solutions, ~\cite{fellner2020solutions,krasnianski2019nonlocal,bitsouni2017mathematical,dyson2013non}, to the convergence of these solutions towards solutions of the associated local models (i.e. where $K[u] = constant$ in Equation~\eqref{advection_part}) when the sensing radius $r$ of the cell approaches zero~\cite{armstrong2006continuum,krasnianski2019nonlocal}. In addition, various studies~\cite{armstrong2006continuum,painter2023biological} focused on pattern formation, as these models can produce very rich spatial and spatio-temporal patterns: from different types of stationary states to traveling fronts and pulses. Analytical approaches to investigate pattern formation start with a stability analysis of spatially homogeneous states (where all individuals/particles are uniformly spread over the domain) and/or heterogeneous states (where individuals/particles are aggregated at particular spatial positions). For example, equation~\eqref{advection_part} could exhibit spatial aggregations only for parameter $\beta$ above a threshold value~\cite{armstrong2006continuum,painter2023biological}. 

Pattern formation in more complex systems of non-local PDEs (with more than 2-3 population variables) in two and even three spatial dimensions is much more difficult to be performed analytically. In this case, numerical bifurcations approaches could help. Nevertheless, continuation algorithms for these non-local PDEs are very few and not standard (due to complexities of these equations), as we will discuss in Section~\ref{Subsect:NumContinuation}.

%%%%%%%%%%%%%%%%%%%%%%%%%%%%%%%%%%%%%%%%%%%%%%%%%%%%%%%%%%%%%%%%%%%%%%%%%%%%%%%%%%%%%%%%%%%%%%%%%%%%%
\section{Numerical approaches for non-local models}
\label{Sect:NumericalApproaches}
Let us now get to the heart of the matter, and discuss the numerical approaches implemented to solve such aggregation-diffusion problems. Although finite differences and finite volumes methods are often used to numerically solve non-local problems, the details of these approaches are rarely made explicit in the literature, for such models. This is why we will begin this section with a brief on the application of finite difference (Section~\ref{Subsect:FiniteDiff}) and finite volume (Section~\ref{Subsect:FiniteVol}) methods to non-local models. Note that since there are numerous such schemes that one could use to discretize these equations, here we will focus only on some simple commonly-used schemes that we also show how to implement numerically.  Then, in the section~\ref{Sect:FiniteElem} we will look at resolution of these non-local models by the finite element method, for which there are very few studies in the literature.

\begin{remark}
    In this article, we have chosen to deal only with ``classical" numerical approaches, namely finite difference, finite volume and finite element methods. Nevertheless, there are other numerical approaches that work well for non-local advection-dominated models, such as the spectral methods. We refer the reader to the papers of A. Gerisch~\cite{gerisch2010mathematical, gerisch2006robust} and A. Buttenschoen~\cite{buttenschon2018space} in which this method is well described and applied to non-local models. We prefer to not present spectral methods in detail in this study, since we want to focus mainly on the very few applications of finite element approaches for these non-local models, which further allows us to discuss the current open problems in this field.
\end{remark}

To present these different numerical methods, we focus on the following simple and generic non-local model~\eqref{nonlocal_model}, where the non-locality is included in the advective part. This model, from the paper of Armstrong, Painter and Sherratt~\cite{armstrong2006continuum}, is a special case of the model~\eqref{nonlocal biology} described previously in Section~\ref{Subsect:CellBiology} and of model~\eqref{model:ecology} described previously in Section~\ref{Subsect:Ecology}: 
\begin{equation}
    \frac{\partial u}{\partial t} = D \Delta u - \nabla \cdot (uK[u]),
    \label{nonlocal_model}
\end{equation} 
with the non-local term $K[u]$ taking the following form 
\begin{equation}
    K[u](\mathbf{x},t) = \alpha \int_{B_r(0)} g(u(\mathbf{x}+\mathbf{y},t)) \omega(||\mathbf{y}||_2) \frac{\mathbf{y}}{||\mathbf{y}||_2}\ d\mathbf{y},
    \label{K}
\end{equation}
with $\mathbf{x} \in \mathbb{R}^n$ with $n \in \{1,2,3\}$, $r$ the cell/animal sensing radius, $\alpha$ the coefficient describing the strength of interactions between cells/animals, $\omega$ the kernel describing the magnitude of these forces according to the distance of the cell/animal located in $\mathbf{y} \in \mathbb{R}^n$ relative to $\mathbf{x}$. Finally, $\mathbf{y} \longmapsto \frac{\mathbf{y}}{||\mathbf{y}||_2}$ is the vector indicating the direction of these forces.\\
Since the goal of this review is to showcase some numerical approaches for these non-local models, throughout this section we consider the following simple functions (which do not incorporate any of the biological realism discussed in Section~\ref{Subsect:CellBiology}): $g(u)=u$ and $\omega(\cdot)=1$.  Note that usually $\int_{B_r} \omega = 2$, but in our case $\omega$ is a constant so we let $\omega=1$ and integrate the constant $\frac{1}{2|B_r|}$ in $\alpha$. The simplicity of these choices will allow us to present clearly the numerical schemes (without being sidetracked by complicated mathematical expressions that could arise from the discretisation of more realistic but highly nonlinear functions $g$ and $\omega$). Also to keep the model as simple as possible, we will consider only the non-local PDE model \eqref{nonlocal_model}-\eqref{K}, and ignore the coupled PDE-ODE models discussed above in equations \eqref{nonlocal system biology wth EDO} and \eqref{ecology diff species}. These simplifications on model equations (which, as mentioned above, might not be very biologically realistic) also allow us to compare, for the purpose of this review, our numerical simulations obtained with different numerical schemes (implemented using an in-house Python code as well as the open-source FEniCS platform) with the numerical simulations in~\cite{armstrong2006continuum}.

Finally, the non-local model~\eqref{nonlocal_model} and~\eqref{K} has to be completed with boundary conditions, which need to consider the impact of the non-local operator [28]. Examples of biologically-realistic boundary conditions are: (i) periodic boundary conditions (where cells/animals leave the domain through one end and re-enter it through the opposite end), which are usually employed to describe ``infinite-like" domains where the boundaries do not affect the dynamics inside the domain; or (ii) zero-flux boundary conditions (where cells/animals cannot cross the domain boundaries). For the simulations performed throughout this Section, we will focus only on periodic boundary conditions.  Mathematically, by working on a 1D interval of length 2L (or a 2D square or a 3D cube with sides of lengths 2L), this amounts to looking for a solution $u$ that is 2L-periodic. Hence, the non-local operator $K[u]$ is also 2L-periodic by
definition \eqref{K}.

\begin{remark}
    For the description of the different numerical schemes in sub-sections~\ref{Subsect:FiniteDiff}, \ref{Subsect:FiniteVol} and~\ref{Subsect:FiniteElemDiscretization} below, we focus only on the following 1D model (although we present also simulations for 2D models):
\begin{equation}
    \partial_t u = -D \partial^2_{xx}u + \partial_x(uK[u]),
    \label{non_local model 1d}
\end{equation}
with 
\begin{equation}
    K[u](x,t) = \alpha \int^r_{-r} u(x+y,t) sgn(y) \ dy,
    \label{K 1D}
\end{equation}
In this $K[u]$ we have considered $\omega(|y|) = 1$, $g(u) = u$, and  $\alpha$ constant. In addition, we work with periodic conditions at the domain boundaries.
\end{remark}

%%%%%%%%%%%%%%%%%%%%%%%%%%%%%%%%%%%%%%%%%%%%%%
\subsection{Finite differences}
\label{Subsect:FiniteDiff}
Before going into the details of the discretization by finite difference of non-local models, let us cite a few articles in which this method is used to carry out simulations:~\cite{mogilner1999non, james2015numerical}.
To describe the finite differences method, we place ourselves in the following framework : the domain is $\Omega = [-L,L]$ discretized in $N$ sub-intervals $[x_i, x_{i+1}]$ with $i=0,..., N-1$. Hence, we note $h = 2L/N$ the mesh step and $\tau$ the time step such as $t_n = n \tau$. \\
For $n \in  \mathbb{N}$, we assume to have computed $(u_i^n)_{i=0,...,N}$ an approximation of $(u(x_i,t_n))_{i=0,...,N}$, and $(K_{i}^n)_{i=0,...,N}$ an approximation of $(K[u](x_i,t_n))_{i=0,...,N}$. At the next time step we obtain an approximation of $u(x_i, t_{n+1})$ denoted $u_i^{n+1}$ by using the following discretization of equation~\eqref{non_local model 1d} :
we consider an Euler finite difference scheme for the approximation of the time derivative
%\begin{itemize}
%  \item {\color{red} Euler} difference schemes for the temporal derivative
            \begin{equation*}
             \frac{\partial}{\partial t}u(x_i,t_n) \approx \frac{u_i^{n+1}-u_i^{n}}{\tau},
        \end{equation*}
        combined with an implicit centered difference scheme for the second-order spatial derivative (for stability reasons~\cite{quarteroni2009numerical})
%    \item second-order centered difference schemes for the second differential operator
            \begin{equation*}
                \frac{\partial^2}{\partial x^2}u(x_i,t_{n+1}) \approx \frac{u_{i+1}^{n+1}-2u_i^{n+1}+u_{i-1}^{n+1}}{h^2},
            \end{equation*}
   % \item centered difference scheme for the first order differential operators, to avoid the issues with the direction of the movement (up-wind or down-wind):
   and a mix implicit-explicit centered difference scheme for the transport term:
            \begin{equation*}
                \frac{\partial}{\partial x}[u(x_i,t_{n+1})K[u](x_i,t_n)] \approx \frac{u^{n+1}_{i+1}K_{i+1}^n-u_{i-1}^{n+1} K_{i-1}^n}{2h},
            \end{equation*}
%\end{itemize}
We could also consider the explicit centered difference scheme for the transport term:
\begin{equation*}
    \frac{\partial}{\partial x}[u(x_i,t_{n})K[u](x_i,t_n)] \approx \frac{u^{n}_{i+1}K_{i+1}^n-u_{i-1}^{n} K_{i-1}^n}{2h},
\end{equation*}
but as we will see in Section~\ref{Paragraph:JustificationScheme}, explicit scheme can leads to numerical oscillations, depending on the value of the time step and the parameters of the equation. In contrast, the mix explicit-implicit scheme for the transport term makes these oscillations disappear.\\
In the previous equations, $K^n_i$ is an approximation of the non-local operator evaluate at point $x_i$ and time $t_n$. This term could be approximated using the following trapezoid quadrature method (or any other quadrature methods~\cite{kress2012numerical}):
\begin{align}
    K[u](x_i,t_n) &= \alpha \int^{x_i+r}_{x_i-r} u(y, t_n) sgn(y-x_i) \ dy\nonumber\\
              &= \alpha \sum_{k=0}^{N_r-1} \biggl( \int^{x_{i+k+1}}_{x_{i+k}} u(y,t_n)  \ dy - \int^{x_{i-k}}_{x_{i-k-1}} u(y,t_n) \ dy \biggr) \nonumber\\
              &\approx \frac{\alpha h}{2} \sum_{k=0}^{N_r-1} \biggl( u_{i+k+1}^n + u_{i+k}^n - u_{i-k}^n - u_{i-k-1}^n \biggr)\nonumber \\
              &= \alpha h \biggl(  \frac{u_{i+N_r}^n - u_{i-N_r}^n}{2} + \sum_{k=1}^{N_r-1} (u_{i+k}^n - u_{i-k}^n) \biggr) = K_i^n,
    \label{approximation trapezes}
\end{align}
with $N_r = \lfloor \frac{r}{h} \rfloor$ the number of sub-intervals $[x_i, x_{i+1}]$ included in ``the sensitivity interval'' $[0,r]$. Moreover, thanks to the periodicity of $u$, $\forall j<0, \ u_j = u_{j+N}$, and similarly $\forall j>N, \ u_j = u_{j-N}$.\\
Thus, $\forall i = 1,...,N-2$, we have the following scheme:
\begin{equation*}
    u_i^{n+1}( 1 + 2 \beta) + u^{n+1}_{i+1}(-\beta + \gamma K^n_{i+1}) + u^{n+1}_{i-1}(-\beta - \gamma K^n_{i-1}) = u^n_i.
\end{equation*}
By considering also periodic conditions at the domain boundary (i.e., $i=0$ and $i=N-1$) leads to the following scheme
\begin{subequations}
\begin{align}
        &u_0^{n+1}( 1 + 2 \beta) + u^{n+1}_{1}(-\beta + \gamma K^n_{1}) + u^{n+1}_{N-1}(-\beta - \gamma K^n_{N-1}) = u^n_0,\\
        &u_i^{n+1}( 1 + 2 \beta) + u^{n+1}_{i+1}(-\beta + \gamma K^n_{i+1}) + u^{n+1}_{i-1}(-\beta - \gamma K^n_{i-1})= u^n_i,\nonumber\\
        &\hspace{3in} \forall i = 1,...,N-2,\\
        &u_{N-1}^{n+1}( 1 + 2 \beta) + u^{n+1}_{0}(-\beta + \gamma K^n_{N}) + u^{n+1}_{N-2}(-\beta - \gamma K^n_{N-2}) = u^n_{N-1},\\
        &u^{n+1}_N = u^{n+1}_0,
    \end{align}
    \label{scheme FD}
\end{subequations}
with $\beta= \frac{D \tau}{h^2}$, $\gamma = \frac{\tau}{2 h}$.\\
In Section~\ref{Subsect:NumericalSimulationVFandDF} we show the numerical implementation of this scheme for model~\eqref{non_local model 1d}-~\eqref{K 1D}.

%%%%%%%%%%%%%%%%%%%%%%%%%%%%%%%%%%%%%%%%%%%%%%%%
\subsection{Finite volumes}
\label{Subsect:FiniteVol}
We now present a discretization of the non-local equation~\eqref{nonlocal_model} using a finite volume approach, as these methods are widely used for such non-local equations~\cite{Bitsouni2018travelWaves2Pop,armstrong2006continuum,bitsouni2017mathematical,delarue2020convergence}.
We start by rewriting equation~\eqref{non_local model 1d} as:
\begin{equation}
    \partial_t u + \partial_x J(x) = 0,
    \label{eq_flux}
\end{equation}
with the flux $J=J_d+J_a$ defined previously~\eqref{diffusive flux}-~\eqref{advective flux} in Section~\ref{Subsect:CellBiology}. To describe the finite volumes method for the non-local problem~\eqref{non_local model 1d}, we start by considering the same domain as in the previous section: $\Omega = [-L,L]$ that we discretize in $N$ sub-intervals $C_i = [x_{i-1/2}, x_{i+1/2}]$ with $x_{i \pm 1/2} = -L +(i \pm 1/2)h$, for $i=1,..., N$. Hence, we note $h = 2L/N$ the mesh step and $\tau$ the time step such that $t_n = n \tau$.\\
Hence, by integrating equation~\eqref{eq_flux} on the interval $C_{i}$ (or a similar volume for higher spatial dimensions), we have:
\begin{equation*}
    \int_{C_i} \partial_t u(x,t) \ dx + \int_{C_i} \partial_x J(x,t) \ dx = 0.
\end{equation*}
We can use the same Euler finite difference scheme  as in Section~\ref{Subsect:FiniteDiff}, for the approximation of the time derivative, to obtain:
\begin{equation*}
    h \frac{u^{n+1}_i - u^n_i}{\tau} + J(x_{i+1/2},t_{n+1}) - J(x_{i-1/2},t_{n+1}) = 0,
\end{equation*}
with $u_i^n$ be the average of $u(x,t_n)$ on $C_i$. As in the previous Section~\ref{Subsect:FiniteDiff}, do not use a completely implicit scheme. Indeed, the diffusive flux is treated implicitly but for the transport flux we use a mix explicit-implicit discretization approach, to avoid numerical oscillations (see Section~\ref{Paragraph:JustificationScheme}). Hence, $J_{i \pm 1/2}^n$ is an approximation of the following flux:
\begin{equation*}
    J_{i \pm 1/2}^{n+1} \approx -D \frac{\partial u(x_{i \pm 1/2}, t_{n+1})}{\partial x} + u(x_{i \pm 1/2}, t_{n+1})K[u](x_{i \pm 1/2}, t_{n}).
\end{equation*}
So we have, $\forall i=1,...,N-1$
\begin{equation}
    u^{n+1}_i + \frac{\tau}{h}(J_{i+1/2}^{n+1} - J_{i-1/2}^{n+1}) = u^n_i.
    \label{FV scheme1}
\end{equation}
Then, we approximate the first spatial derivative with a central difference scheme (centered at $x_{i+1/2}$)
\begin{equation*}
    \frac{\partial u(x_{i + 1/2}, t_{n+1})}{\partial x} \approx \frac{u^{n+1}_{i+1}-u^{n+1}_i}{h},
\end{equation*}
and the value of the non-local term $K[u]$ and of $u$ at point $x_{i \pm 1/2}$ by
\begin{align*}
    &K[u](x_{i\pm1/2}, t_n) \approx \frac{K^{n}_{i\pm1}+K^{n}_i}{2}, \;\;\;\;\; u(x_{i\pm1/2},t_{n+1}) \approx \frac{u^{n+1}_{i\pm1}+u^{n+1}_i}{2},
\end{align*}
with $K^{n}_i$ an approximation of $K[u](x_i,t_n)$ calculated as in the previous Section~\ref{Subsect:FiniteDiff}. Hence, the scheme~\eqref{FV scheme1} becomes
\begin{align*}
    &u^{n+1}_i + \frac{\tau}{h} \biggl( -D \frac{u^{n+1}_{i+1}-u^{n+1}_i}{h} + \frac{u^{n+1}_{i+1}+u^{n+1}_i}{2} \frac{K^{n}_{i+1}+K^{n}_i}{2} \\
    &+ D \frac{u^{n+1}_{i}-u^{n+1}_{i-1}}{h} - \frac{u^{n+1}_{i}+u^{n+1}_{i-1}}{2} \frac{K^{n}_{i}+K^{n}_{i-1}}{2} \biggr) = u^n_i,
\end{align*}
and so
\begin{align*}
    &u^{n+1}_i \biggl( 1+2\beta + \gamma \frac{K^n_{i+1}-K^n_{i-1}}{2} \biggr) + u^{n+1}_{i+1} \biggl( -\beta + \gamma \frac{K^n_{i+1}+K^n_{i}}{2} \biggr) \nonumber \\
    &+ u^{n+1}_{i-1} \biggl( -\beta - \gamma \frac{K^n_{i}+K^n_{i-1}}{2} \biggr) = u^n_i, \ \forall i=1,...,N-2,
\end{align*}
with $\beta = \frac{D \tau}{h^2}$ and $\gamma = \frac{\tau}{2h}$. \\
Moreover, by adding periodic boundaries conditions, we obtain the following system:
\begin{subequations}
    \begin{align}
        &u^{n+1}_0 \biggl( 1+2\beta + \gamma \frac{K^n_{1}-K^n_{N-1}}{2} \biggr) + u^{n+1}_{1} \biggl( -\beta + \gamma \frac{K^n_{1}+K^n_{0}}{2} \biggr) \nonumber \\
        &\quad+ u^{n+1}_{N-1} \biggl( -\beta - \gamma \frac{K^n_{0}+K^n_{N-1}}{2} \biggr) = u^n_0,\\
        &u^{n+1}_i \biggl( 1+2\beta + \gamma \frac{K^n_{i+1}-K^n_{i-1}}{2} \biggr) + u^{n+1}_{i+1} \biggl( -\beta + \gamma \frac{K^n_{i+1}+K^n_{i}}{2} \biggr) \nonumber \\
        &\quad + u^{n+1}_{i-1} \biggl( -\beta - \gamma \frac{K^n_{i}+K^n_{i-1}}{2} \biggr) = u^n_i, \ \forall i=1,...,N-2,\\
        &u^{n+1}_{N-1} \biggl( 1+2\beta + \gamma \frac{K^n_{0}-K^n_{N-2}}{2} \biggr) + u^{n+1}_{0} \biggl( -\beta + \gamma \frac{K^n_{0}+K^n_{N-1}}{2} \biggr) \nonumber \\
        &\quad + u^{n+1}_{N-2} \biggl( -\beta - \gamma \frac{K^n_{N-1}+K^n_{N-2}}{2} \biggr) = u^n_{N-1},\\
        &u^{n+1}_N = u^{n+1}_0.
    \end{align}
    \label{scheme VF}
\end{subequations}
In Section~\ref{Subsect:NumericalSimulationVFandDF} we will see an example of numerical simulations obtained with the implementation of this finite volume scheme for model`\eqref{non_local model 1d}-~\eqref{K 1D}.

%%%%%%%%%%%%%%%%%%%%%%%%%%%%%%%%%%%%%%%%%%%%%%
\subsection{Numerical simulation}
\label{Subsect:NumericalSimulationVFandDF}
\begin{figure}[!ht]
     \centering
     \begin{subfigure}
         \centering
         \includegraphics[width=0.49\textwidth]{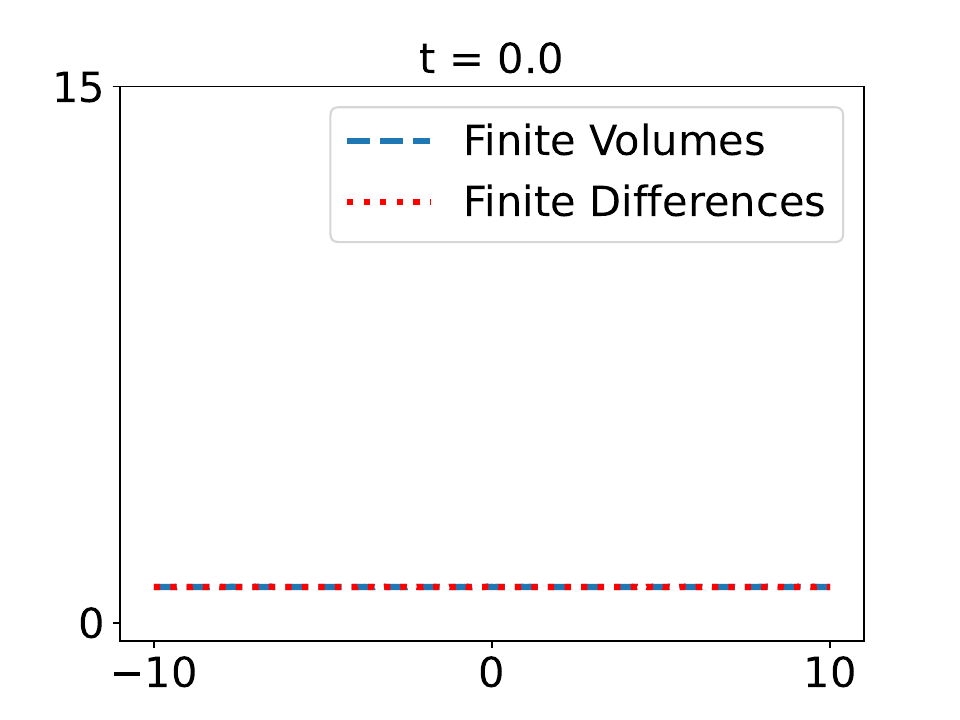}
     \end{subfigure}
     \begin{subfigure}
         \centering
         \includegraphics[width=0.49\textwidth]{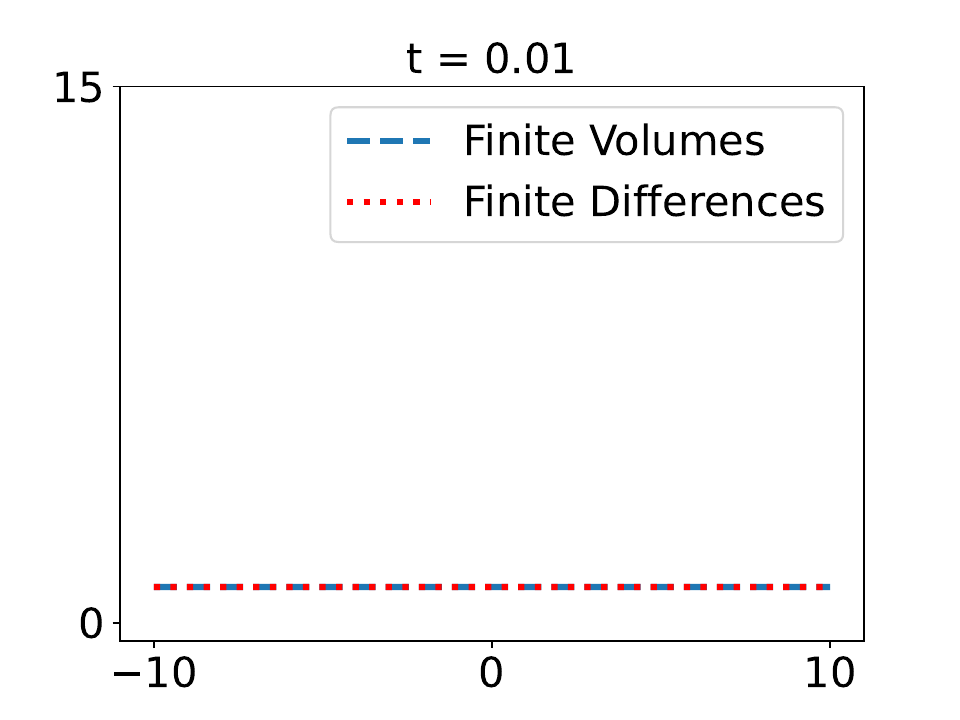}
     \end{subfigure}
     \begin{subfigure}
         \centering
         \includegraphics[width=0.49\textwidth]{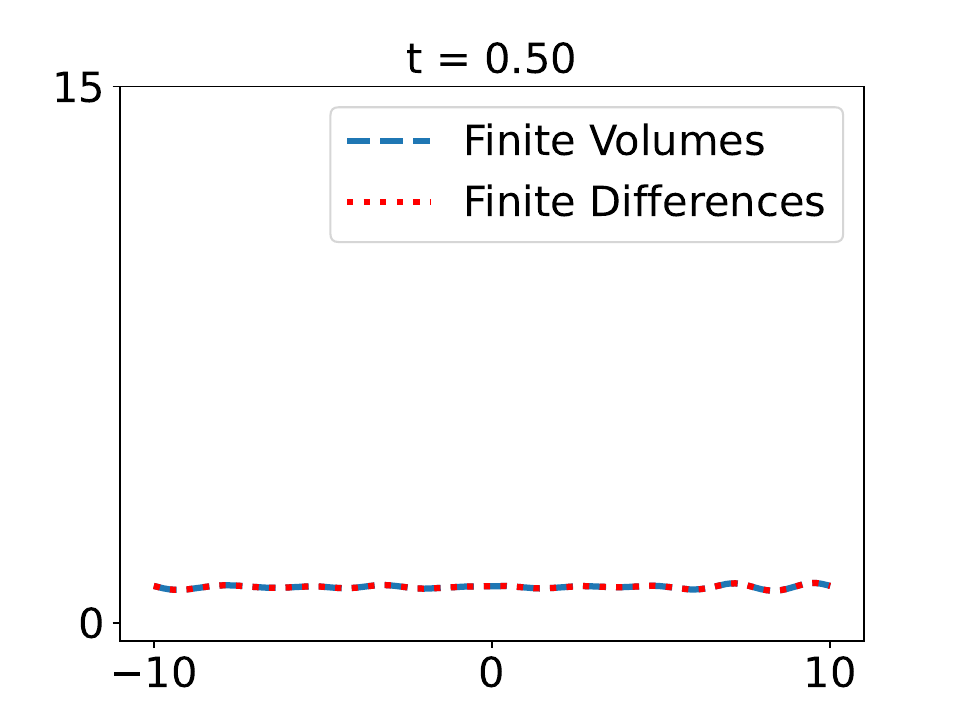}
     \end{subfigure}
     \begin{subfigure}
         \centering
         \includegraphics[width=0.49\textwidth]{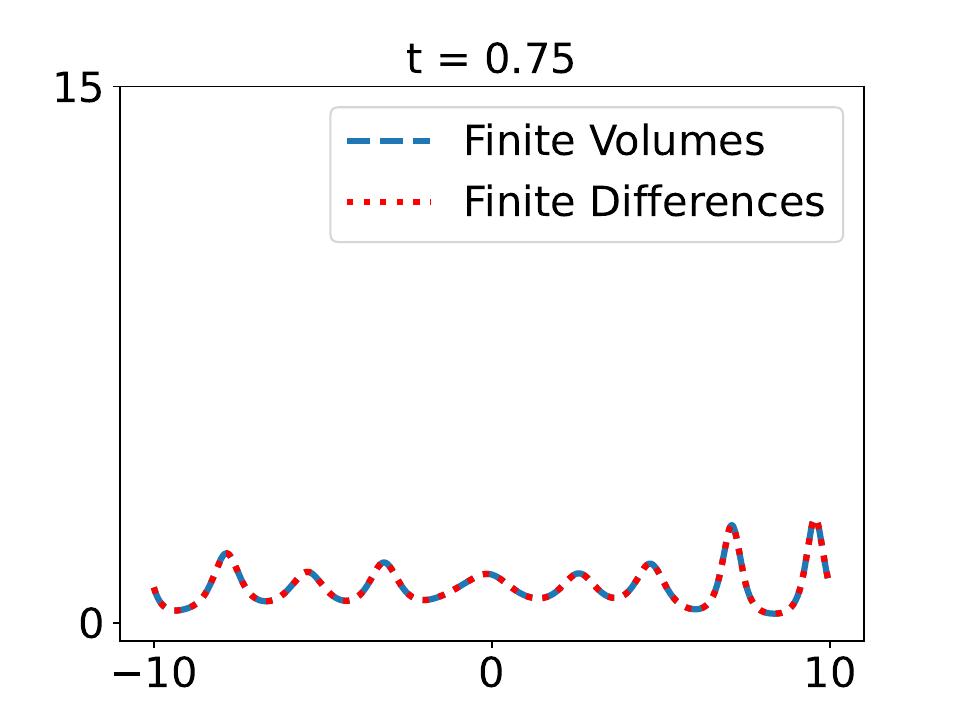}
     \end{subfigure}
     \begin{subfigure}
         \centering
         \includegraphics[width=0.49\textwidth]{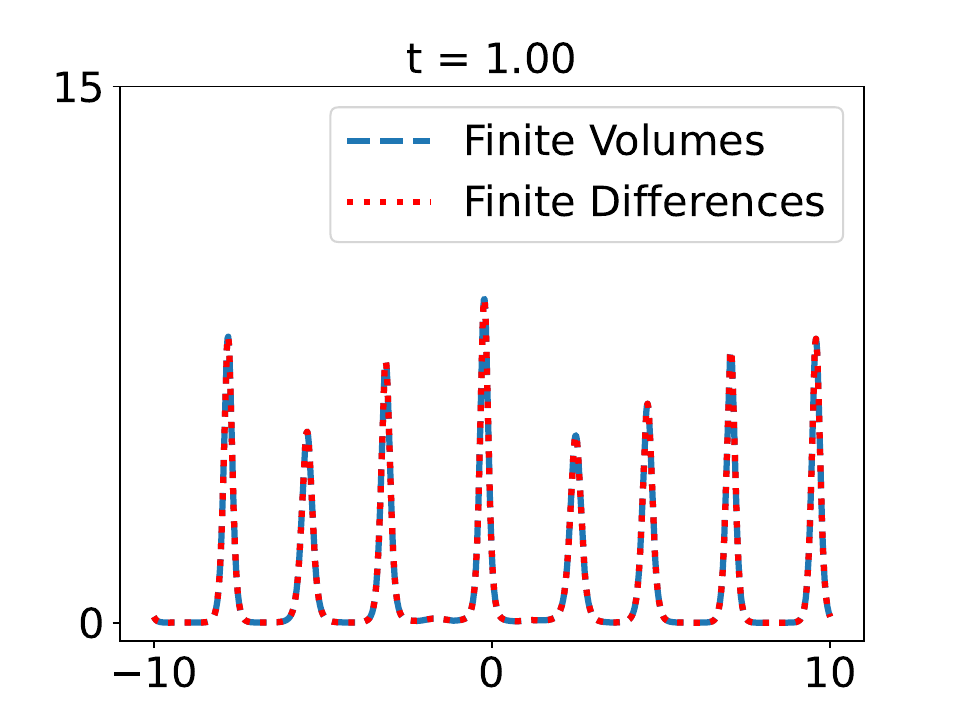}
     \end{subfigure}
     \begin{subfigure}
         \centering
         \includegraphics[width=0.49\textwidth]{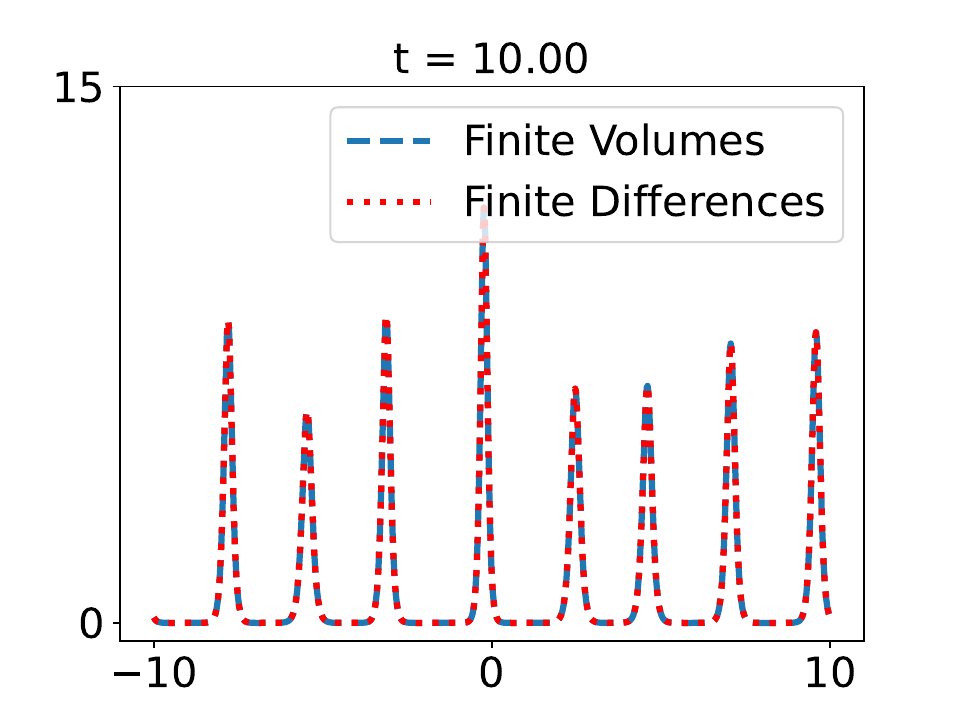}
     \end{subfigure}
     \caption{Numerical resolution of equation~\eqref{non_local model 1d} using numerical schemes~\eqref{scheme FD} (red dashed curve) and~\eqref{scheme VF} (blue continuous curve). We show here the solution curves at times $t=0$, $t=0.01$, $t=0.5$, $t=0.75$, $t=1$ and $t=10$. The parameter values are as in Table~\ref{Tab:Param1d}.}
    \label{Fig:FDandFV}
\end{figure}

To verify that the two schemes finite differences~\eqref{scheme FD} and finite volumes~\eqref{scheme VF} give the same result, we implement the resolution of these two schemes in Python (see Figure~\ref{Fig:FDandFV}). We take as initial condition
\begin{equation*}
    u_0(x) = 1 + \text{random}(0,10^{-2}),
\end{equation*}
and all parameters of the equation~\eqref{non_local model 1d} and of the discretization are the same in the two different schemes, see Table~\ref{Tab:Param1d}.
\begin{longtable}{|p{2.5cm}|p{3.5cm}|p{5.5cm}|}
\caption{Description of model~\eqref{non_local model 1d}-~\eqref{K 1D} parameters with their values used for the numerical simulations presented in Figures~\ref{Fig:FDandFV},~\ref{Fig:1D_FEniCS_results}, and~\ref{Fig:comparison_schemes}.} \label{Tab:Param1d} \\
\hline
Parameters & Value & Description \\ 
\hline \hline
$D$ & 1 & Diffusion coefficient\\
\hline
$\alpha$ & 10 & Adhesion coefficient\\
\hline
$r$ & 1 & Sensitivity radius\\
\hline
$L$ & 10 & $\Omega = [-L,L]$\\
\hline
$N$ & 1000 & Number of intervals\\
\hline
$\tau$ & 0.01 & Step time\\
\hline
\end{longtable}

%%%%%%%%%%%%%%%%%%%%%%%%%%%%%%%%%%%%%%%%%%%%%%%%%%%%%%%%%%%%%%%%%%%%%%%%%%%%%%%%%%%%%%%%%%%%%%%%%%%%%
\section{A new numerical approach for non-local models: finite elements method}
\label{Sect:FiniteElem}
Although various numerical approaches have already been considered for different non-local models~\cite{mogilner1999non, james2015numerical,Bitsouni2018travelWaves2Pop,armstrong2006continuum,bitsouni2017mathematical,delarue2020convergence}, very few deal with the finite element method due to the mathematical technicalities caused by the non-local terms~\cite{adebayo2023mathematical}. This is why we present here in detail the implementation of this method for the non-local aggregation-diffusion equation~\eqref{nonlocal_model} with periodic boundaries conditions seen previously. We start with the technical details of this numerical method; then we see the different numerical simulations carried out in both Python (in-house code) and FEniCS (open access computing platform) in the one-dimensional (1D) and two-dimensional (2D) cases, for one or two cell populations.

%%%%%%%%%%%%%%%%%%%%%%%%%%%%%%%%%%%%%%%%%%%%%
\subsection{Finite Elements discretization}
\label{Subsect:FiniteElemDiscretization}
To simplify the discretization by finite elements, we place ourselves in the 1D case~\eqref{non_local model 1d}. Let $\Omega = [-L,L]$ denote the 1D spatial domain and $I = [0,T]$ denote the temporal interval. Thus, the problem in which we are interested is the equation ~\eqref{non_local model 1d} where $K[u]$ is the non-local form~\eqref{K 1D}, coupled with periodic conditions at the edges:
\begin{equation}
    \begin{cases}
        \frac{\partial u(x,t)}{\partial t} = D  \frac{\partial^2 u(x,t)}{{\partial x}^2} - \frac{\partial}{\partial x} (u(x,t)K[u](x,t)),\\
        u(x,t=0) = u_0(x),\\
        u \text{ is $2L-$periodic}.
    \end{cases}
    \label{problem}
\end{equation}
For spatial discretization, we use the finite element method. Thus, we multiply the first line of the system~\eqref{problem} by a test function $v$ which belongs to the function space $H(\Omega)$. By defining:
\begin{equation*}
    H(\Omega) = \{ v \in H^1(\Omega), v \text{ is } 2L-\text{periodic}\}.
\end{equation*}
Next, we integrate over the spatial domain the first line of system~\eqref{problem} multiplied by $v$ and apply Green's formula. Thus, we obtain:
\begin{align*}
    \int_{-L}^L \frac{\partial u(x,t)}{\partial t} v(x) \ dx &= \int_{-L}^L \biggl[ D\frac{\partial^2 u(x,t)}{{\partial x}^2} - \frac{\partial}{\partial x} (u(x,t)K[u](x,t)) \biggr] v(x) \ dx,
\end{align*}
and so,
\begin{align*}
    \int_{-L}^L \frac{\partial u(x,t)}{\partial t} v(x) \ dx &= D \biggl( \biggl[ \frac{\partial u(x,t)}{\partial x} v(x,t)\biggr]_{-L}^L - \int_{-L}^L \frac{\partial u(x,t)}{\partial x} \frac{\partial v(x)}{\partial x} \ dx \biggr)\\
     - \biggl[ u(x,t)K[u](x,t)& v(x,t)\biggr]_{-L}^L  + \int_{-L}^L u(x,t)K[u](x,t) \frac{\partial v(x)}{\partial x} \ dx.
\end{align*}
But since $u$ and $v$ are $2L-$periodic, the terms in square brackets, in the second equality, are zero. Hence, the variational formulation of the problem~\eqref{problem} is :\\\\
Find $u \in H(\Omega)$ such that, for all $v \in H(\Omega)$,
\begin{subequations}
    \begin{align}
        &\int_{-L}^L \frac{\partial u(x,t)}{\partial t} v(x) \ dx + D \int_{-L}^L \frac{\partial u(x,t)}{\partial x} \frac{\partial v(x)}{\partial x} \ dx \nonumber \\
        &- \int_{-L}^L u(x,t)K[u](x,t) \frac{\partial v(x)}{\partial x} \ dx = 0.
    \end{align}
    \label{weak_form}
\end{subequations}
Now, for the time discretization, we use a standard Euler scheme (as in Section~\ref{Subsect:FiniteDiff} and~\ref{Subsect:FiniteVol}). Furthermore, we choose to treat the diffusion term implicitly, while we treat the advective part using a mixed explicit-implicit dicretisation approach. These choices were made in order to avoid stability problems of the numerical scheme. We will see later (Figure~\ref{Fig:comparison_schemes} in Section~\ref{Paragraph:JustificationScheme}) a numerical illustration of this instability which will justify our choice to not treat the advective part in a purely explicit manner. Thus, the weak form~\eqref{weak_form} becomes:
\\\\
Find $u^{n+1} \in H(\Omega)$ such that, for all $v \in H(\Omega)$,
\begin{equation}
       a_h(u^{n+1},v ; u^n)=  l_h(u^{n};v),
    \label{weak form}
\end{equation}
with
\begin{align*}
    &a_h(u^{n+1},v;u^n) = \int_\Omega u^{n+1}v \ dx + D\tau \int_\Omega  \frac{\partial u^{n+1}}{\partial x} \frac{\partial v}{\partial x} \ dx\ - \tau \int_\Omega u^{n+1}K(u^n) \frac{\partial v}{\partial x} \ dx,\\
    &l_h(u^n;v) = \int_\Omega u^n v \ dx.
\end{align*}
\\
Before presenting the problem in its discrete form, we discretize the domain $\Omega$ into $N$ intervals, $[x_i,x_{i+1}]=[x_0+ih,x_0+(i+1)h]$ for $i=0,...,N-1$, of length $h=\frac{2L}{N}$ and $\tau$ the time step such as $t_n = n\tau$.
Now, by defining the discrete space of finite elements $\mathbb{P}_1:=\{P(x) = ax+b |a,b \in \mathbb{R}\}$:
\begin{equation*}
    V_h = \{ v \in \mathcal{C}^0(\Omega), \ v_{|[x_{j}, x_{j+1}]} \in \mathbb{P}_1[X], \ \forall j=0,...,N-1; v \text{ is $2L-$periodic}\},
\end{equation*}
we can therefore write the discrete variational formulation of the problem~\eqref{problem}:
\begin{equation}
    \begin{cases}
       \text{Find }u_h^{n+1} \in V_h\text{ such that, }\forall v_h \in V_h,\\
       a(u_h^{n+1},v_h; u_h^n) = l(u_h^n;v_h).
    \end{cases}
    \label{variational form}
\end{equation}
But solving this problem~\eqref{variational form} amounts to solve the following problem (see Appendix~\ref{Appendix:CalculationDetails} for calculation details):
\begin{subequations}
    \begin{align}
        &u^{n+1}_0 \biggl(\frac{2}{3} + 2\beta + 2 \gamma(K^n_{1}-K^n_{N-1})\biggr) + u^{n+1}_{1} \biggl(\frac{1}{6} - \beta + 2 \gamma(2K^n_{0}+K^n_{1})\biggr) \nonumber\\
        &\quad + u^{n+1}_{N-1} \biggl(\frac{1}{6} - \beta - 2 \gamma(2K^n_{0}+K^n_{N-1})\biggr) = 4u^n_0 + u^n_{1} + u^n_{N-1},\\
        &u^{n+1}_i \biggl(\frac{2}{3} + 2\beta + 2 \gamma(K^n_{i+1}-K^n_{i-1})\biggr) + u^{n+1}_{i+1} \biggl(\frac{1}{6} - \beta + 2 \gamma(2K^n_{i}+K^n_{i+1})\biggr) \nonumber \\
        &\quad + u^{n+1}_{i-1} \biggl(\frac{1}{6} - \beta - 2 \gamma(2K^n_{i}+K^n_{i-1})\biggr) = 4u^n_i + u^n_{i+1} + u^n_{i-1}, \forall i=1,...,N-2,\\
        &u^{n+1}_{N-1} \biggl(\frac{2}{3} + 2\beta + 2 \gamma(K^n_{0}-K^n_{N-2})\biggr) + u^{n+1}_{0} \biggl(\frac{1}{6} - \beta + 2 \gamma(2K^n_{N-1}+K^n_{0})\biggr) \nonumber\\
        &\quad + u^{n+1}_{N-2} \biggl(\frac{1}{6} - \beta - 2 \gamma(2K^n_{N-1}+K^n_{N-2})\biggr) = 4u^n_{N-1} + u^n_{0} + u^n_{N-2},\\
        &u^{n+1}_N = u^{n+1}_0,
    \end{align}
    \label{FE system}
\end{subequations}
with $\beta = \frac{D \tau}{h^2}$ and $\gamma = \frac{\tau}{2h}$. In the previous system, we introduce $K_j^n$ which is the following approximation of the non-local operator $K[u](x_j,t_n)$ (see also Appendix~\ref{Appendix:CalculationDetails} for calculation details): 
\begin{equation*}
    K^n_j := K(u_h)(x_j,t_n) \approx - \alpha h \mathcal{F}^{-1}(\mathcal{F}(U^n_h)\mathcal{F}(\omega))(x_j,t_n),
\end{equation*}
where $\mathcal{F}$ (resp. $\mathcal{F}^{-1}$) is the Fourier transform (resp. inverse Fourier transform).

%%%%%%%%%%%%%%%%%%%%%%%%%%%%%%%%%%%%%%%%%%%%%%%%%
\subsection{Numerical simulations}
\label{Subsect:NumericalSimuFEM}
Now that we have described the dizcretisation with finite element method for the non-local problem~\eqref{nonlocal_model}, we can move on to the different numerical simulations carried out with an in-house Python code (which does not use finite element libraries) and with a FEniCS code (which uses already-implemented finite element libraries). 

%%%%%%%%%%%%%%%%%%%%%%%%
\subsubsection{One population model}
\begin{figure}[h!]
     \centering
     \begin{subfigure}
         \centering
         \includegraphics[width=0.49\textwidth]{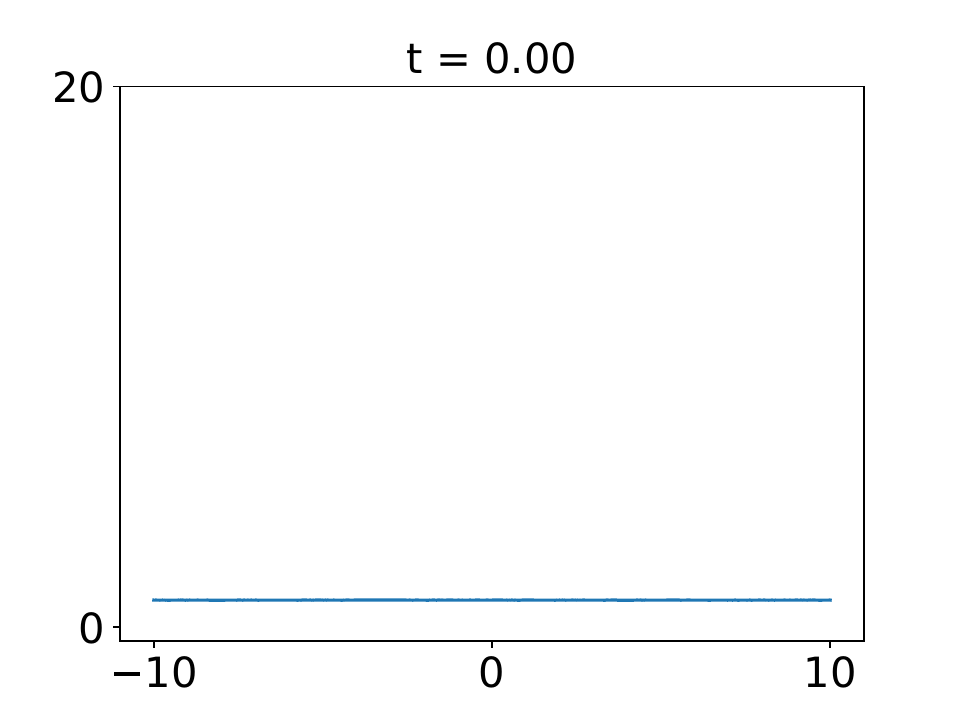}
     \end{subfigure}
     \begin{subfigure}
         \centering
         \includegraphics[width=0.49\textwidth]{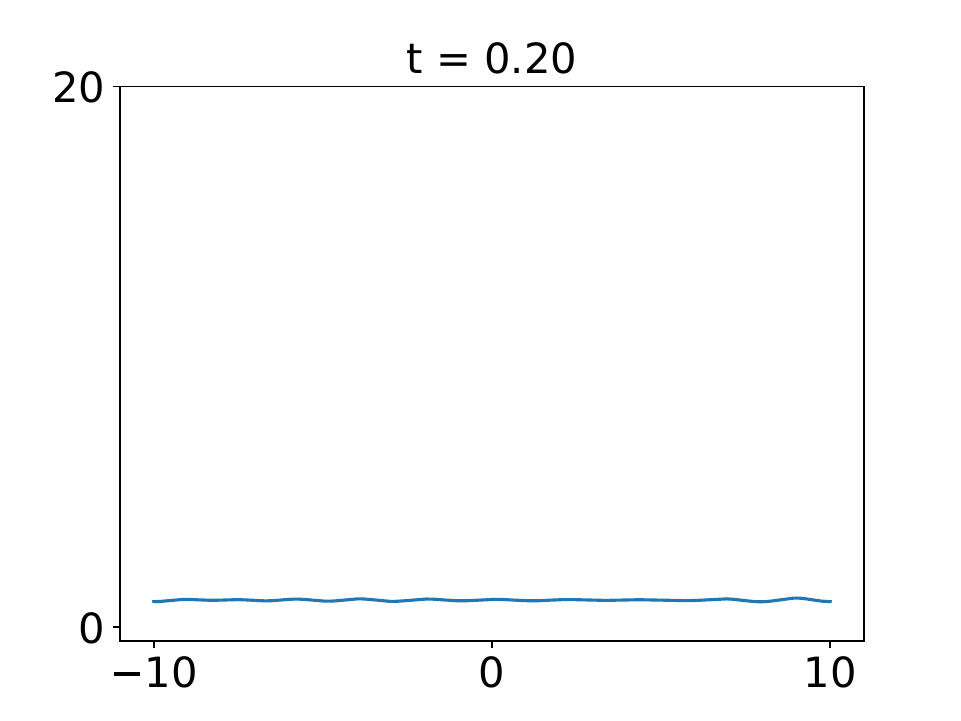}
     \end{subfigure}
     \begin{subfigure}
         \centering
         \includegraphics[width=0.49\textwidth]{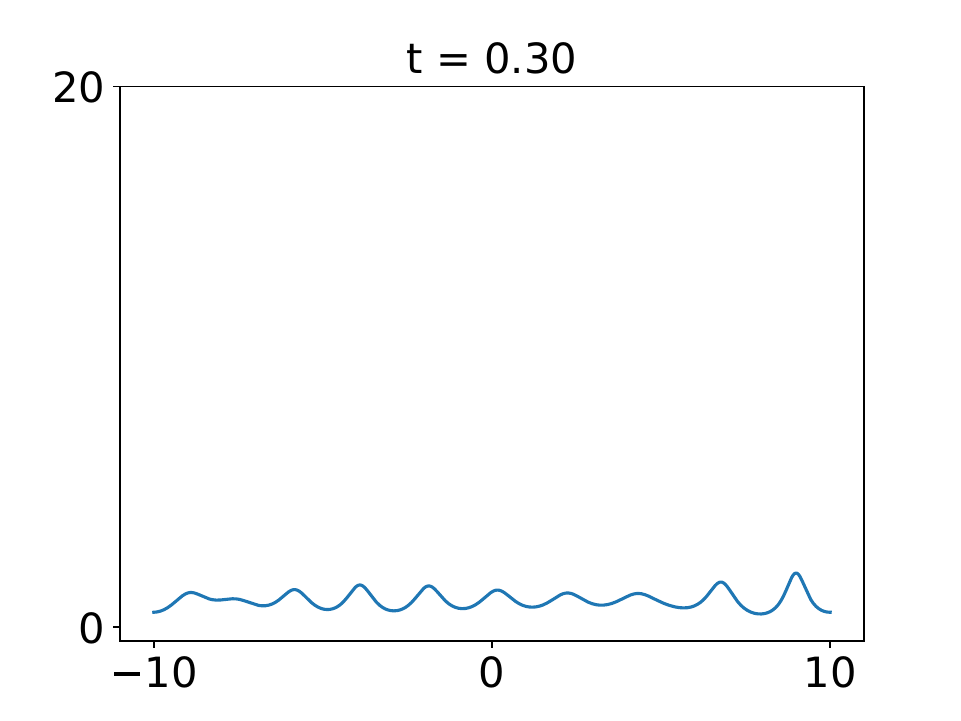}
     \end{subfigure}
     \begin{subfigure}
         \centering
         \includegraphics[width=0.49\textwidth]{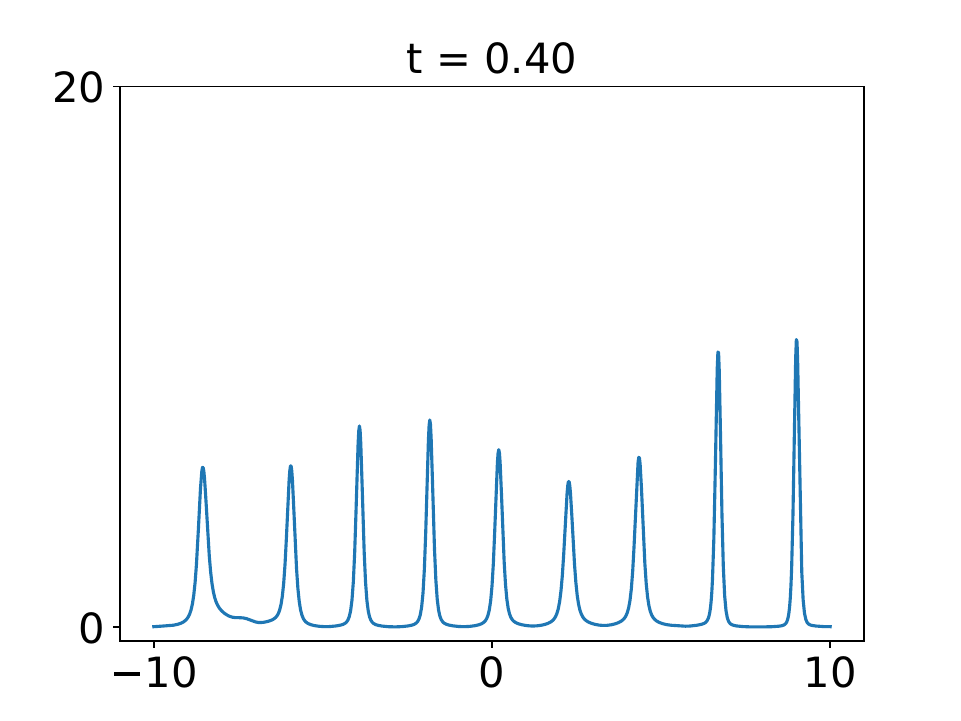}
     \end{subfigure}
     \begin{subfigure}
         \centering
         \includegraphics[width=0.49\textwidth]{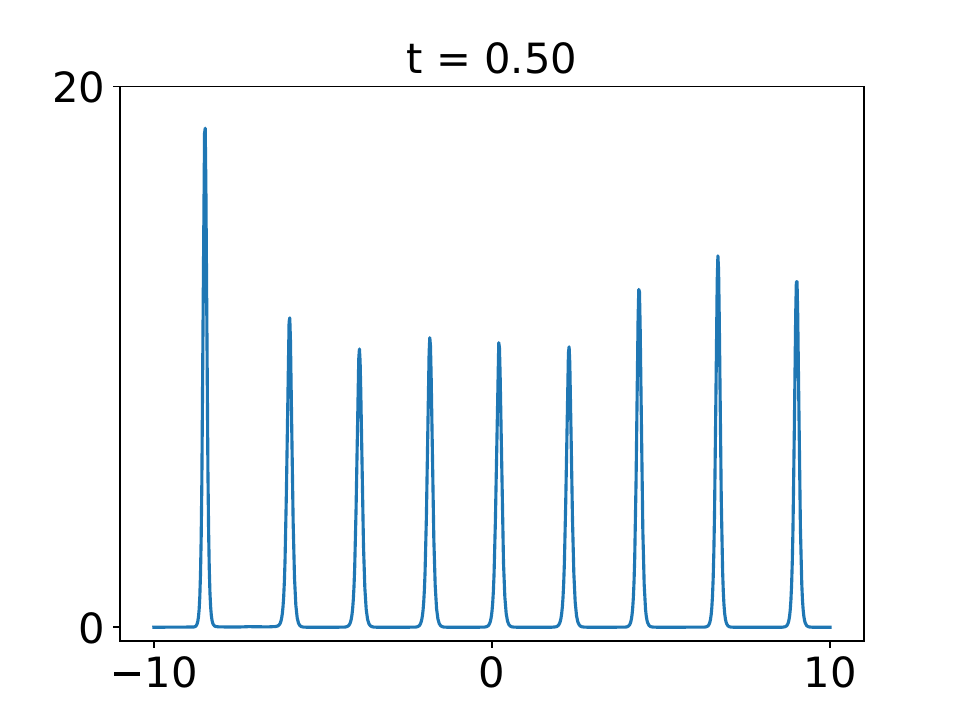}
     \end{subfigure}
     \begin{subfigure}
         \centering
         \includegraphics[width=0.49\textwidth]{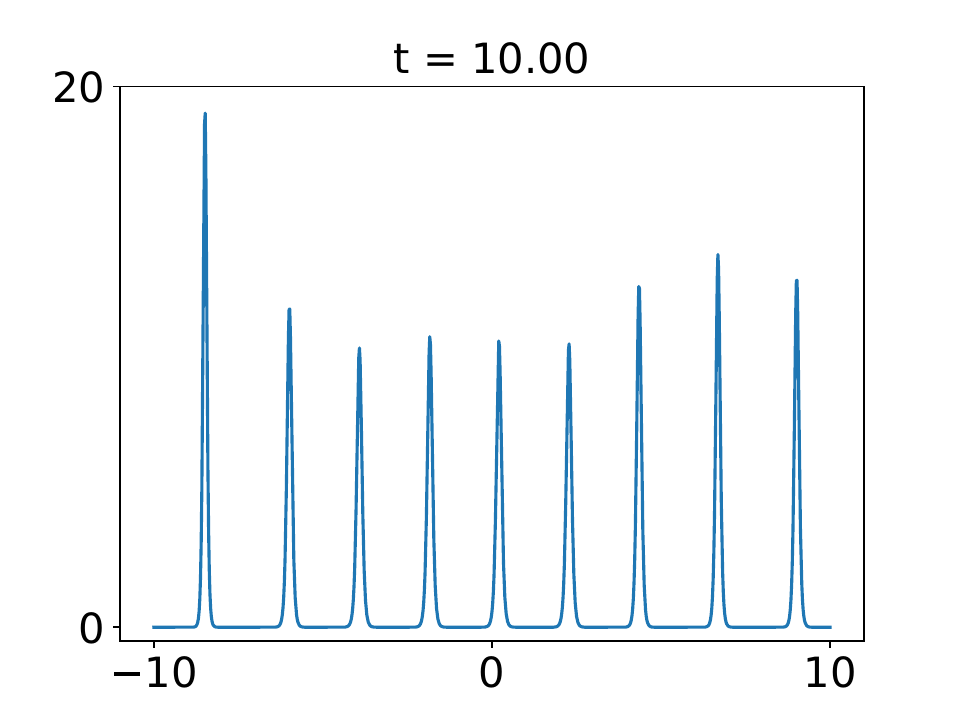}
     \end{subfigure}
     \caption{Numerical solution of the equation~\eqref{non_local model 1d} using FEM implemented in FEniCS, at different time. Parameter values are as in Table~\ref{Tab:Param1d}.}
    \label{Fig:1D_FEniCS_results}
\end{figure}

\paragraph{Uni-dimensional case.} For the 1D case~\eqref{non_local model 1d}, the adhesive flow $K[u]$ take the form~\eqref{K 1D}. In Figure~\ref{Fig:1D_FEniCS_results} we use the FEniCS software to numerically solve the discrete system\eqref{variational form}. We also solved~\eqref{variational form} using an in-house Python code ran with the same parameter values (see Table~\ref{Tab:Param1d})
and the results were the same (not shown here). For the
initial conditions, see Section~\ref{Subsect:NumericalSimulationVFandDF}.
 In the simulation~\ref{Fig:1D_FEniCS_results}, we can see that amplitude of the solution become stable for large time.

\begin{figure}[h!]
     \centering
     \begin{subfigure}
         \centering
         \includegraphics[width=0.49\textwidth]{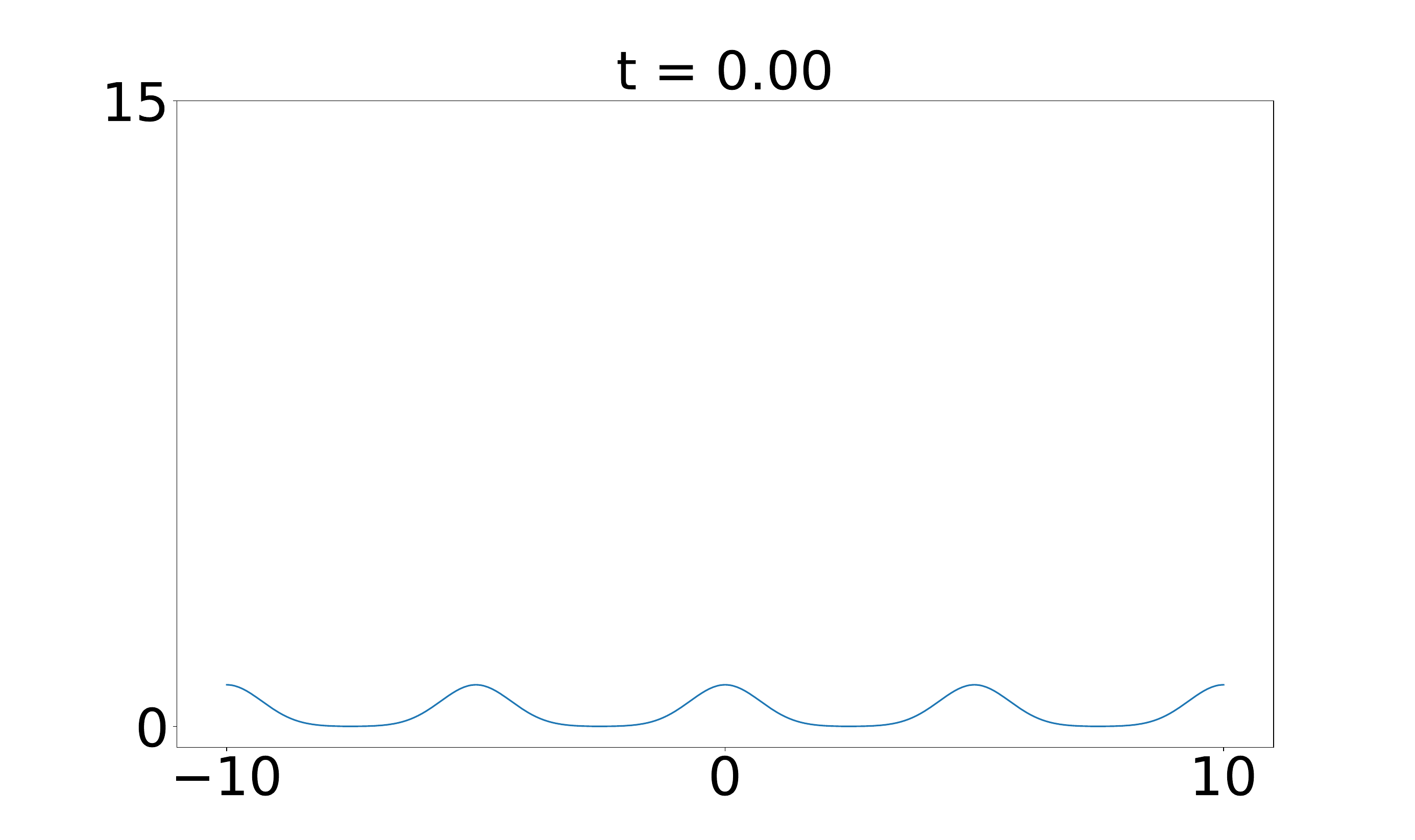}
     \end{subfigure}
     \begin{subfigure}
         \centering
         \includegraphics[width=0.49\textwidth]{Fig4a.pdf}
     \end{subfigure}
     \begin{subfigure}
         \centering
         \includegraphics[width=0.49\textwidth]{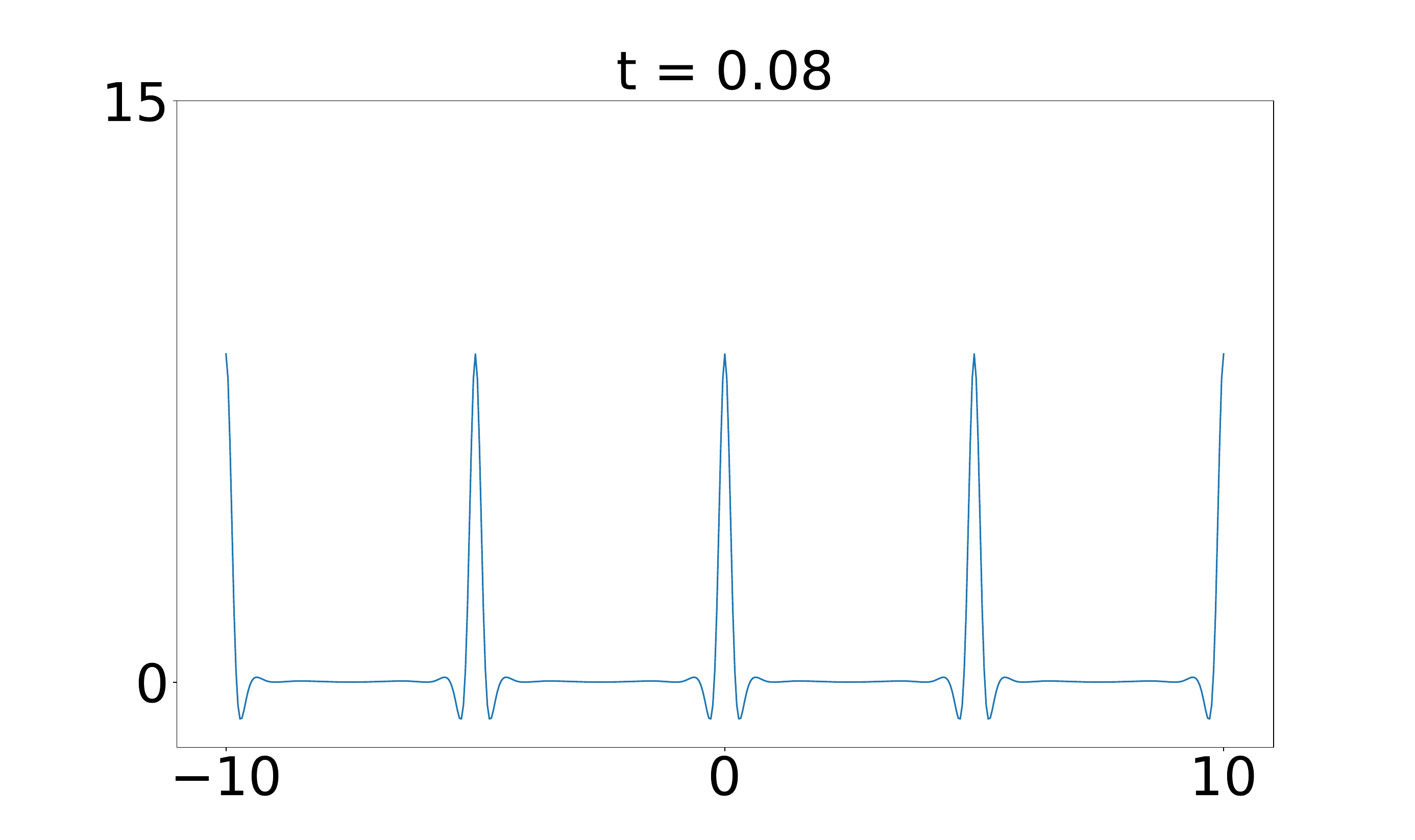}
     \end{subfigure}
     \begin{subfigure}
         \centering
         \includegraphics[width=0.49\textwidth]{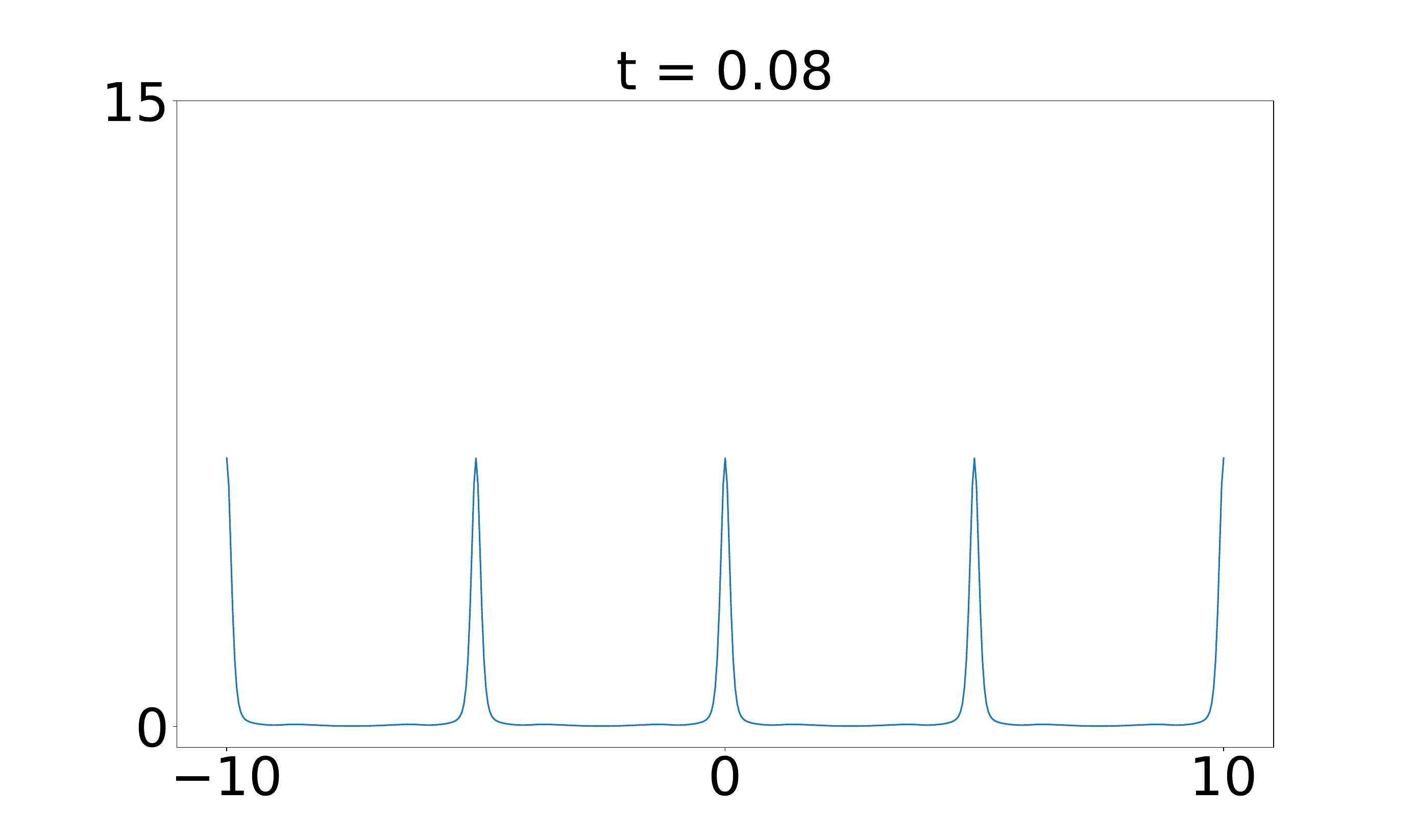}
     \end{subfigure}
     \begin{subfigure}
         \centering
         \includegraphics[width=0.49\textwidth]{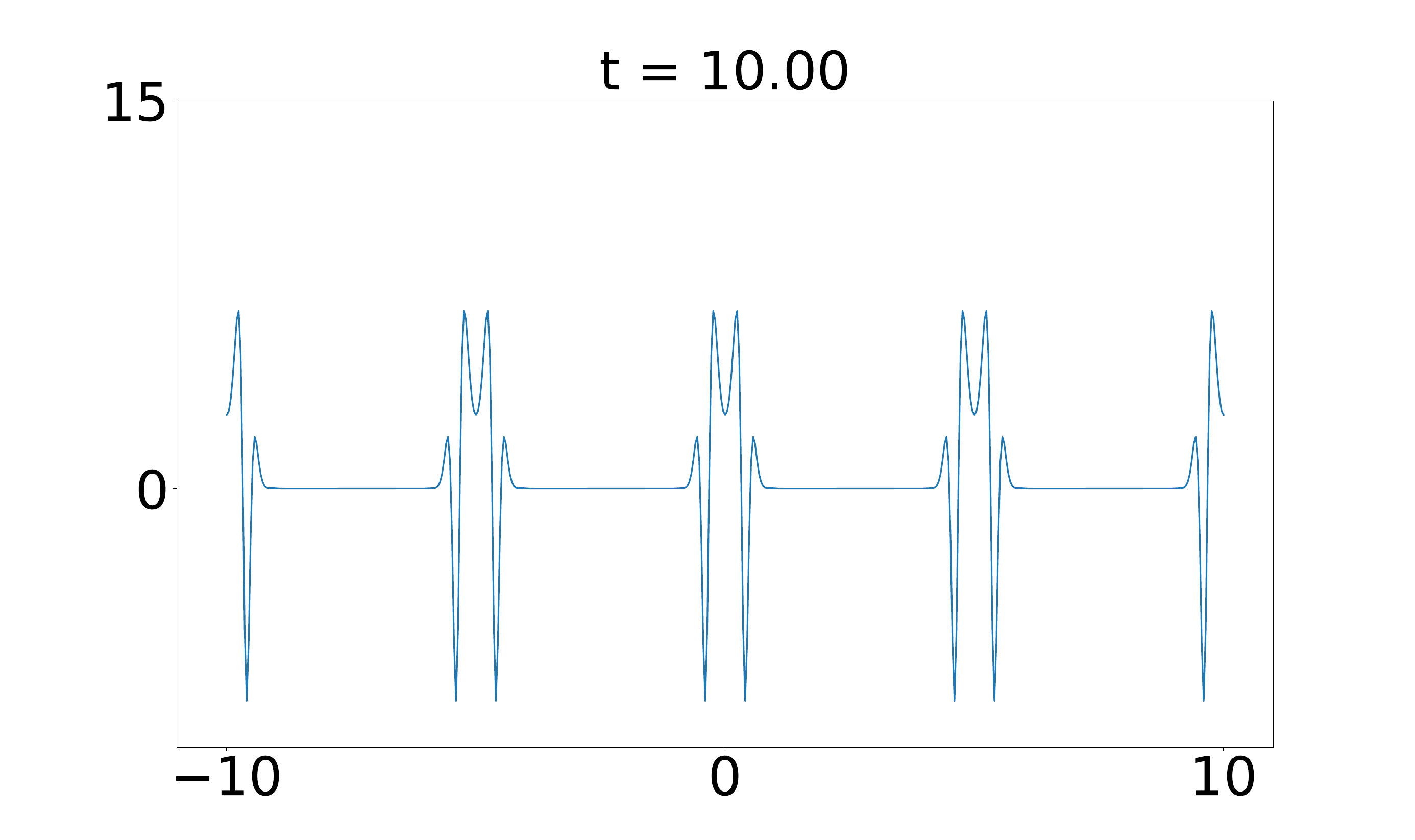}
     \end{subfigure}
     \begin{subfigure}
         \centering
         \includegraphics[width=0.49\textwidth]{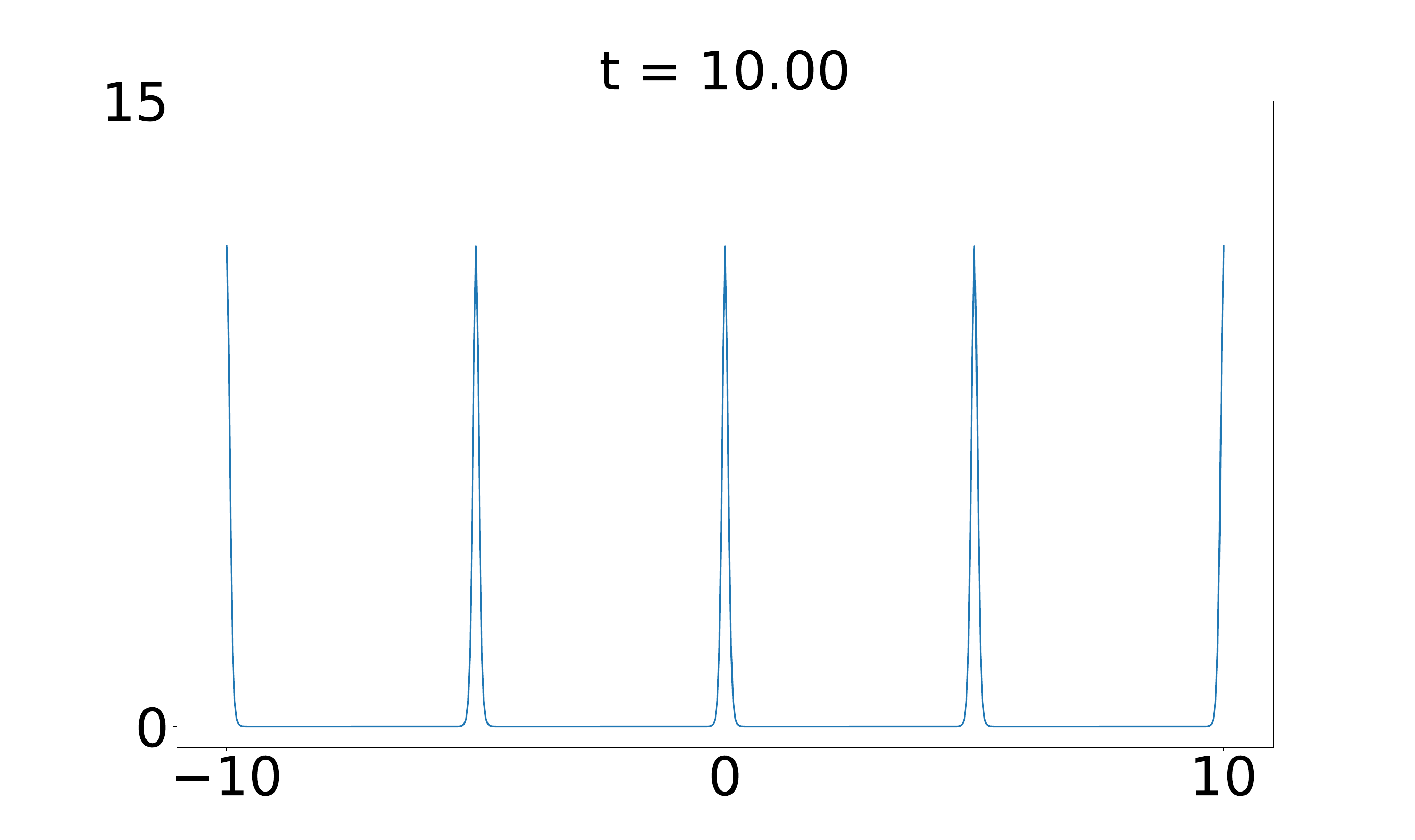}
     \end{subfigure}
     \caption{Numerical simulations of model~\eqref{non_local model 1d}-~\eqref{K 1D}, using the parameter values in Table~\ref{Tab:Param1d} (except $\alpha=15$) obtained with : (left column) the explicit scheme~\eqref{explicit}, (right column) the mix explicit-implicit scheme.}
    \label{Fig:comparison_schemes}
\end{figure}

\paragraph{Comparison between explicit scheme and  mix explicit-implicit scheme, for the advective part.}
\label{Paragraph:JustificationScheme}
This paragraph aims to justify our choice of choosing a mix explicit-implicit scheme to treat the advection part of the equation~\eqref{non_local model 1d}. Indeed, we have carried out two numerical simulations with these two schemes, for the same parameters values (see Table~\ref{Tab:Param1d}). We notice that for the explicit scheme for the advective part, i.e. the following scheme:
\begin{equation}
    \int_\Omega u^{n+1}v \ dx + D\tau \int_\Omega  \frac{\partial u^{n+1}}{\partial x} \frac{\partial v}{\partial x} \ dx\ = \tau \int_\Omega u^{n}K(u^n) \frac{\partial v}{\partial x} \ dx +\int_\Omega u^n v \ dx,
    \label{explicit}
\end{equation}
some numerical oscillations appear if $\tau$ is not very small, whereas those oscillations disappear when we use the mix explicit-implicit scheme~\eqref{weak form} for the advective part, see Figure~\ref{Fig:comparison_schemes}. For the initial conditions we have chosen:
\begin{equation*}
    u_0(x)= e^{-x^2} + e^{-(x-L)^2} + e^{-(x+L)^2} + e^{-(x-L/2)^2} +e^{-(x+L/2)^2}.
\end{equation*}
Thus, instead of using the explicit scheme and having to refine the time step $\tau$ to make these oscillations disappear, we choose to use the explicit-implicit mix scheme. These simulations illustrate our choice of scheme for the finite element part, but these results would have been similar for the finite difference and finite volume schemes.

\begin{figure}[h!]
     \centering
     \begin{subfigure}
         \centering
         \includegraphics[width=0.49\textwidth]{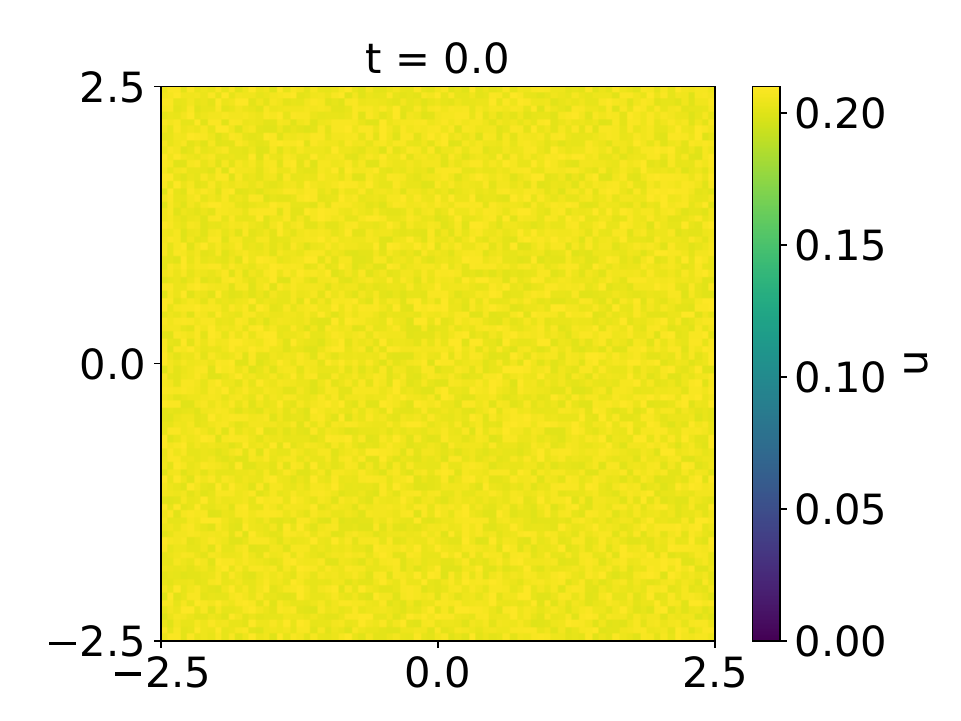}
     \end{subfigure}
     \begin{subfigure}
         \centering
         \includegraphics[width=0.49\textwidth]{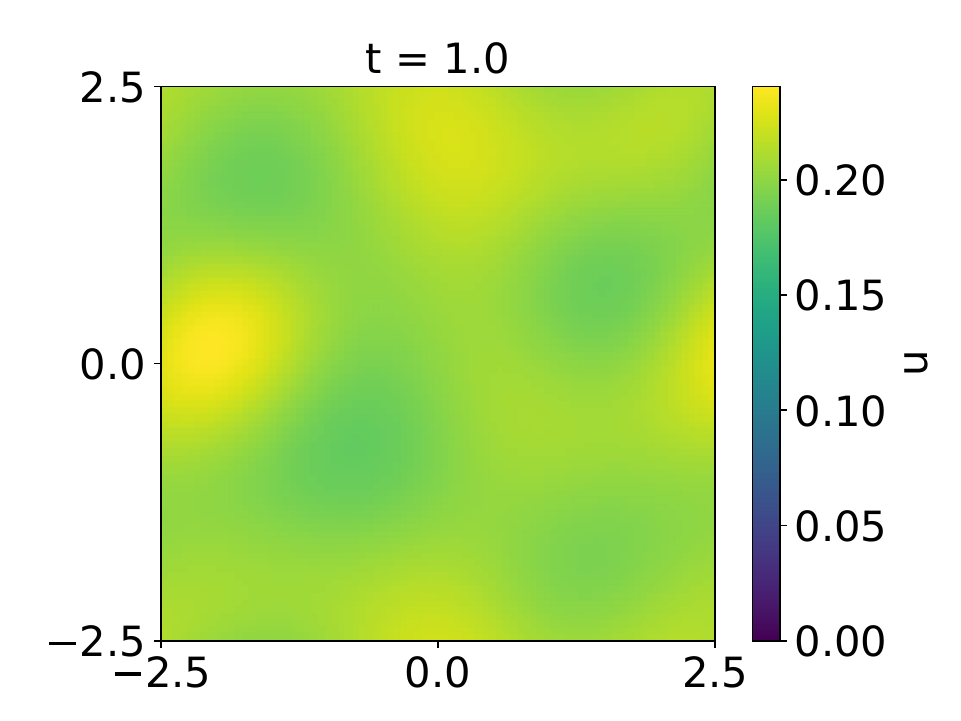}
     \end{subfigure}
     \begin{subfigure}
         \centering
         \includegraphics[width=0.49\textwidth]{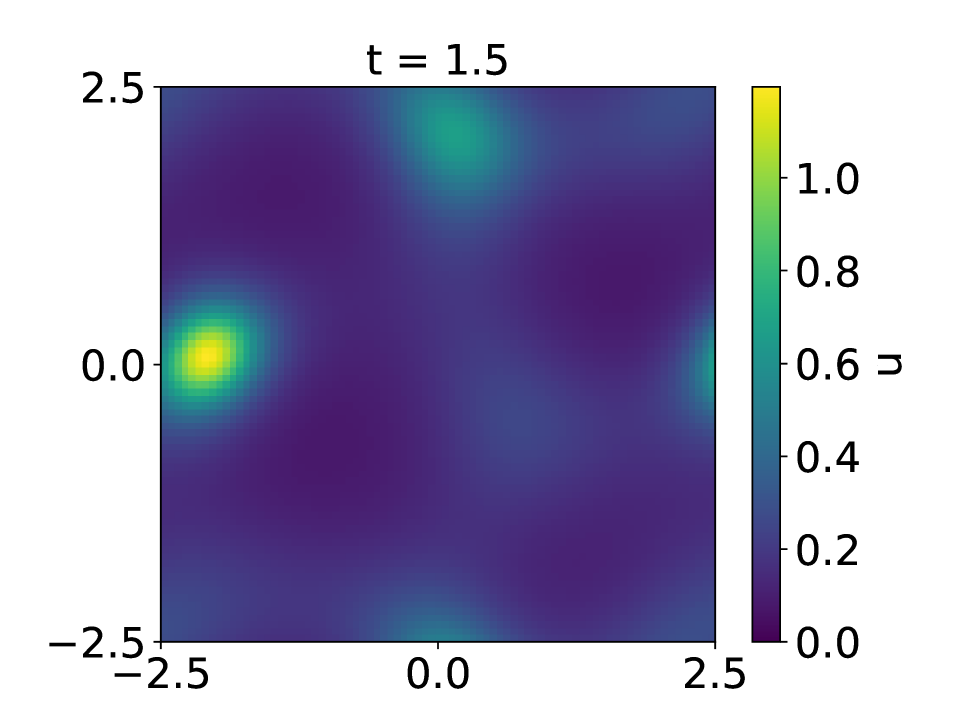}
     \end{subfigure}
     \begin{subfigure}
         \centering
         \includegraphics[width=0.49\textwidth]{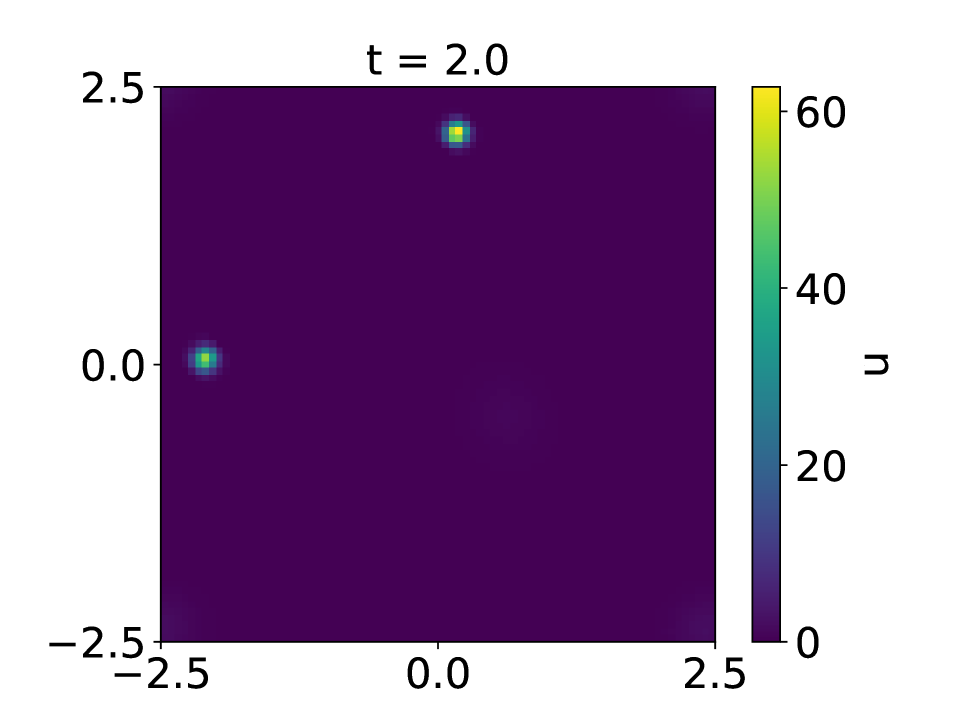}
     \end{subfigure}
     \begin{subfigure}
         \centering
         \includegraphics[width=0.49\textwidth]{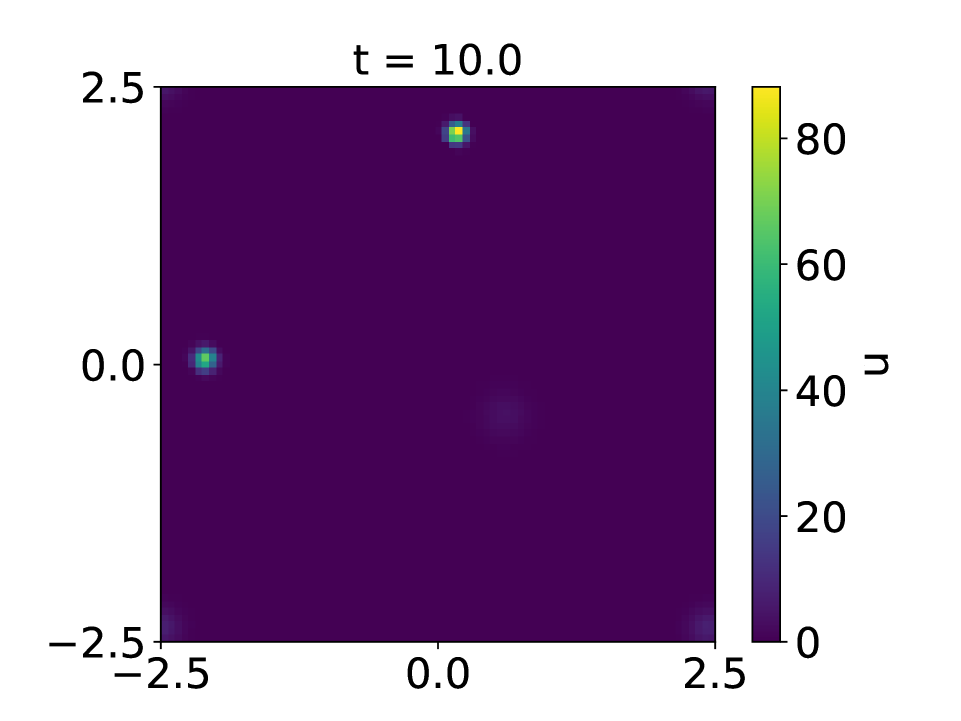}
     \end{subfigure}
     \begin{subfigure}
         \centering
         \includegraphics[width=0.49\textwidth]{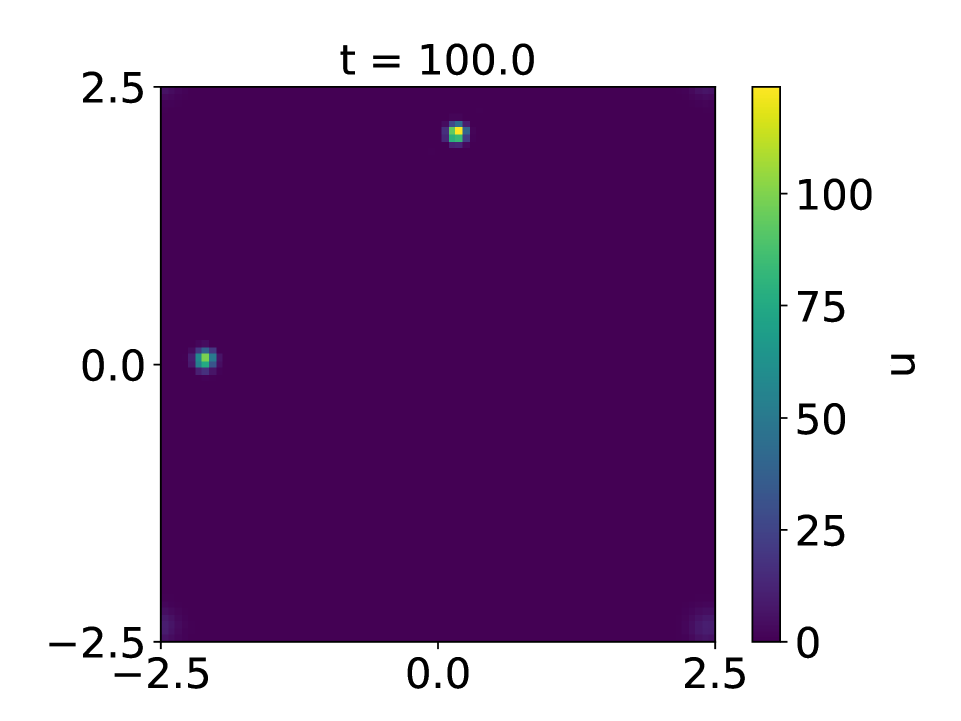}
     \end{subfigure}
     \caption{Numerical solution of the bi-dimensional equation~\eqref{nonlocal_model} using FEM implemented in FEniCS. We show the solution at six different times: $t=0$, $t=1$, $t=1.5$, $t=2$, $t=10$, $t=100$. The parameter values used for these simulations are listed in Table~\ref{Tab:Param2d}.}
    \label{Fig:2D_FEniCS_results}
\end{figure}

\paragraph{Bi-dimensional case.} Next we perform simulations with the finite element method for the 2D version of model~\eqref{nonlocal_model} using FEniCS software (see Table~\ref{Tab:Param2d} for parameters values). As for the uni-dimensional case, we also perform simulations using an in-house
Python code ran with the same parameter values (see Table~\ref{Tab:Param2d}) and the results were the same (not shown here). For the initial conditions we have chosen
\begin{equation*}
u_0(x,y) = 0.2 + \text{random}(0,10^{-2}).
\end{equation*}
Moreover, in the non-linear term $K[u]$, the linear form $g$ and the kernel $\omega$ are the same as for the 1D case (i.e. $g(u)=u$, $\omega(||y||_2) = \frac{1}{2|B_r|}$ and $f(u)=0$).
The results of these simulations are presented in Figure~\ref{Fig:2D_FEniCS_results}. We see that the small 2D perturbations give rise to spatial aggregations with finite densities, which are stable in time. Moreover, the solution is positive. Therefore, the numerical scheme behaves well.

\begin{longtable}{|p{2.5cm}|p{3.5cm}|p{5.5cm}|}
\caption{Description of model parameters with their values used for the 2D numerical simulations presented in Figure~\ref{Fig:2D_FEniCS_results}.} \label{Tab:Param2d} \\
\hline
Parameters & Value & Description \\ 
\hline \hline
$D$ & 1 & Diffusion coefficient\\
\hline
$\alpha$ & 10 & Adhesion coefficient\\
\hline
$r$ & 1 & Sensitivity radius\\
\hline
$L$ & 2.5 & $\Omega = [-L,L] \times [-L,L]$\\
\hline
$N$ & 80 & Number of intervals\\
\hline
$\tau$ & 0.1 & Step time\\
\hline
\end{longtable}

%%%%%%%%%%%%%%%%%
\subsubsection{Two populations model}
We now consider the modelling of two cell populations, $u$ and $v$, interacting through adhesion. In this new model, we must take into account the different attraction and adhesion relationships between cells. For each cell type, we assume that there are two adhesion forces: the first representing the adhesion force between cells of the same population, and the second representing the adhesion force between two cells of different populations. Thus, the non-local model for two populations of cells is as follows:
\begin{equation}
    \begin{cases}
        \frac{\partial u}{\partial t} = D \Delta u - \nabla \cdot (uK_u(u,v)) , \\
        \frac{\partial v}{\partial t} = D \Delta v - \nabla \cdot (vK_v(u,v)) ,
    \end{cases}
    \label{two population}
\end{equation}
where $D$ is the diffusion and, $K_u$ and $K_v$ are the following non-local forms:
\begin{align}
    K_{u}(\mathbf{x},t) &= S_u \int_{B_r(0)} g_{uu}(u(\mathbf{x}+\mathbf{y},t), v(\mathbf{x}+\mathbf{y},t)) \omega(||\mathbf{y}||_2) \frac{\mathbf{y}}{||\mathbf{y}||_2}\ dy \nonumber\\
    &+ C \int_{B_r(0)} g_{uv}(u(\mathbf{x}+\mathbf{y},t), v(\mathbf{x}+\mathbf{y},t)) \omega(||\mathbf{y}||_2) \frac{\mathbf{y}}{||\mathbf{y}||_2}\ dy,
    \label{Ku}
\end{align}
and,
\begin{align}
    K_{v}(\mathbf{x},t) &= S_v \int_{B_r(0)} g_{vv}(u(\mathbf{x}+\mathbf{y},t), v(\mathbf{x}+\mathbf{y},t)) \omega(||\mathbf{y}||_2) \frac{\mathbf{y}}{||\mathbf{y}||_2}\ dy \nonumber\\
    &+ C \int_{B_r(0)} g_{vu}(u(\mathbf{x}+\mathbf{y},t), v(\mathbf{x}+\mathbf{y},t)) \omega(||\mathbf{y}||_2) \frac{\mathbf{y}}{||\mathbf{y}||_2}\ dy,
    \label{Kv}
\end{align}
with $S_u$ (respectively $S_v$) the self-adhesive strength of population $u$ (respectively $v$), $C$ the cross-adhesive strength between the populations $u$ and $v$, $\omega$ describe the magnitude of these forces and $g_i$ for $i \in [uu, uv, vv, vu]$ functions defining the dependence of the adhesive force on cell density. This model, very common in the literature, has already been solved numerically in numerous papers~\cite{armstrong2006continuum, painter2023biological}.

\paragraph{Uni-dimensional case.}
We focused on the one dimensional case of the model~\eqref{two population}-~\eqref{Kv} coupled with the following initial condition
\begin{align*}
    &u_0(x) = 0.2 + \text{rand}(0,10^{-2}) ,\\
    &v_0(x) = 0.2 + \text{rand}(0,10^{-2}) .
\end{align*}
Moreover, unlike the simulations of the model for one population of cells, where we took a very simple form for the function $g$, we take here a more complicated form but more realistic from a biological point of view, for the functions $g_i$ (for $i \in [uu, uv, vv, vu]$)
\begin{equation}
    g_{uu} = g_{vu} =
    \begin{cases}
        u(1-u-v),  \text{ if } u+v < 1,\\
        0, \text{ else},
    \end{cases}
    \label{guu and gvu}
\end{equation}
and
\begin{equation}
    g_{vv} = g_{uv} =
    \begin{cases}
        v(1-u-v),  \text{ if } u+v < 1,\\
        0, \text{ else.}
    \end{cases}
    ,
    \label{gvv and guv}
\end{equation}
These two forms of $g_i$ functions, with, $i \in [uu, uv, vv, vu]$, make it possible to include a ``population pressure'' which reflects the fact that cells are attracted only towards regions where cell density is lower than 1~\cite{armstrong2006continuum}.
\begin{longtable}{|p{1.75cm}|p{2.5cm}|p{8cm}|}
\caption{Description of 1D case model~\eqref{two population}-~\eqref{Kv} parameters with their values used for the numerical simulations presented in Figure~\ref{Fig:TwoPopulations1DResults}.} \label{Tab:Param1dtwopop} \\
\hline
Parameters & Value & Description \\ 
\hline \hline
$D$ & 1 & Diffusion coefficient\\
\hline
$S_u$ & 25 or 200 & Self-adhesive strength of population $u$ \\
\hline
$S_v$ & 7.5 or 15 or 25 & Self-adhesive strength of population $v$ \\
\hline
$C$ & 0 or 5 or 30 or 50 & Cross-adhesive strength between population $u$ and $v$\\
\hline
$r$ & 1 & Sensitivity radius\\
\hline
$L$ & 10 & $\Omega = [-L,L]$\\
\hline
$N$ & 1000 & Number of intervals\\
\hline
$\tau$ & 0.01 & Step time\\
\hline
\end{longtable}
Now that the model for two populations is well defined, we can discretize it using the same approach as for the one-population model~\eqref{non_local model 1d}. Hence, we solved this model~\eqref{two population} using finite elements. Then we implemented its resolution with FeniCS. We tested the resolution of this finite element scheme with different parameter values $S_u$, $S_v$ and $C$.
\begin{figure}[h!]
     \centering
     \begin{subfigure}
         \centering
         \includegraphics[width=0.235\textwidth]{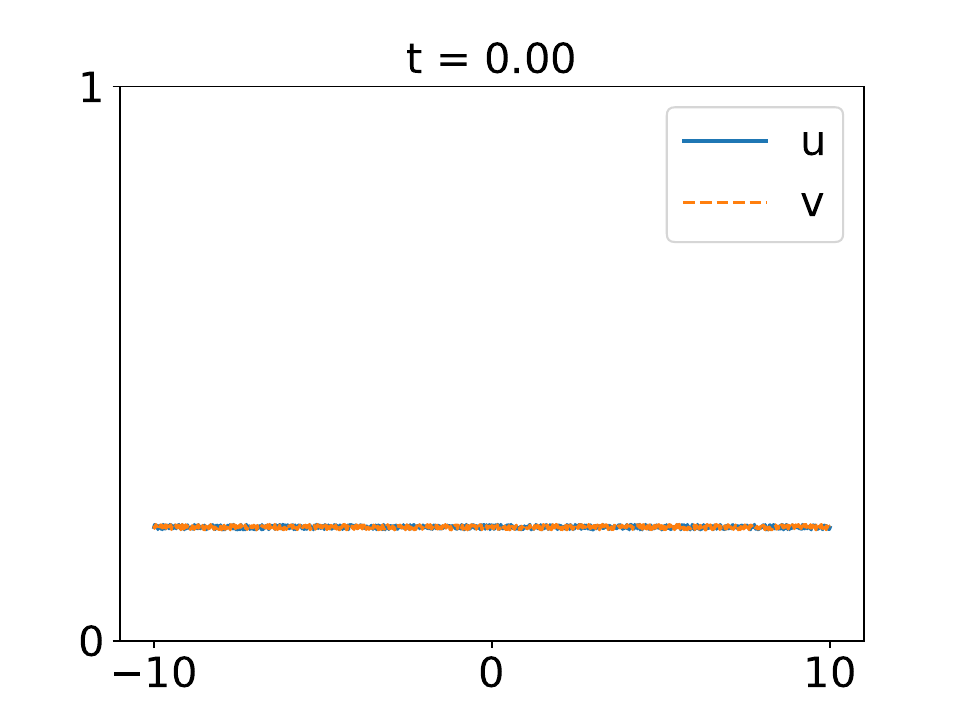}
     \end{subfigure}
     \begin{subfigure}
         \centering
         \includegraphics[width=0.235\textwidth]{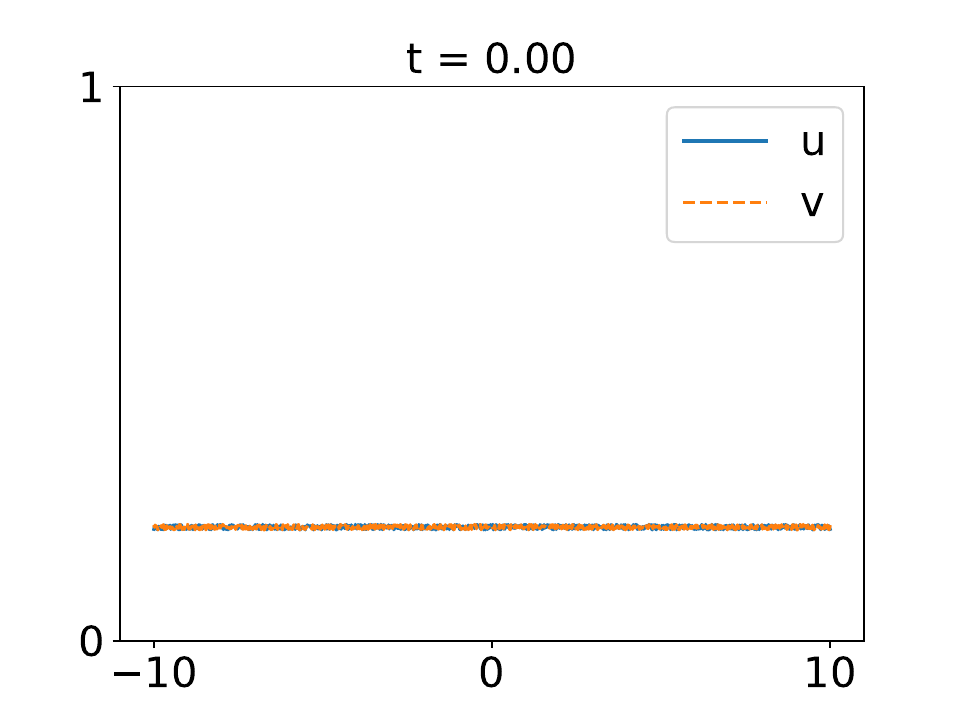}
     \end{subfigure}
     \begin{subfigure}
         \centering
         \includegraphics[width=0.235\textwidth]{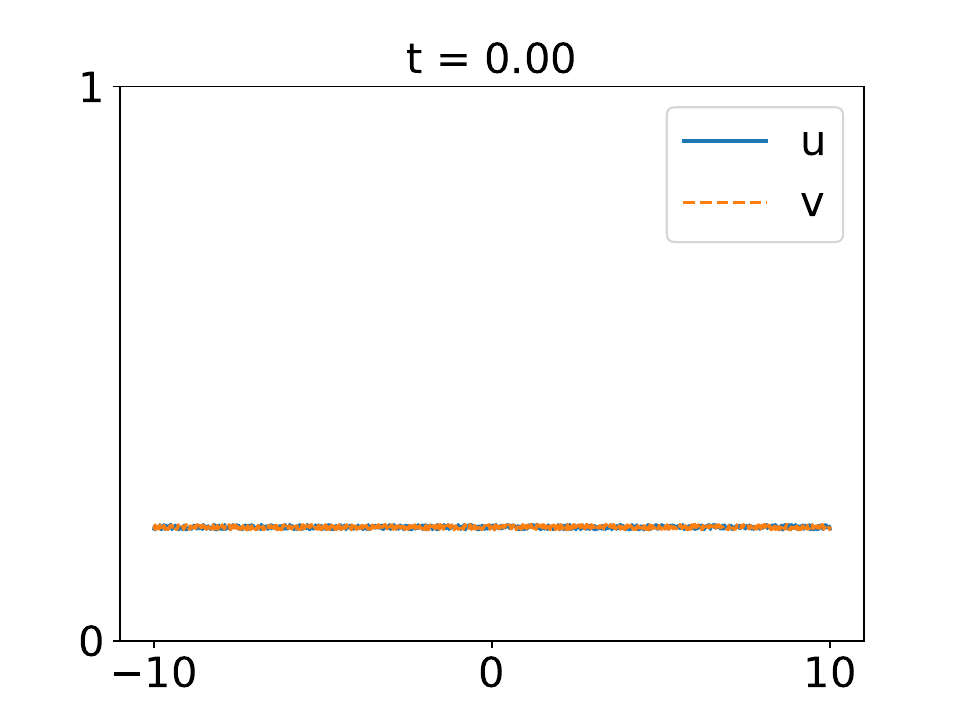}
     \end{subfigure}
     \begin{subfigure}
         \centering
         \includegraphics[width=0.235\textwidth]{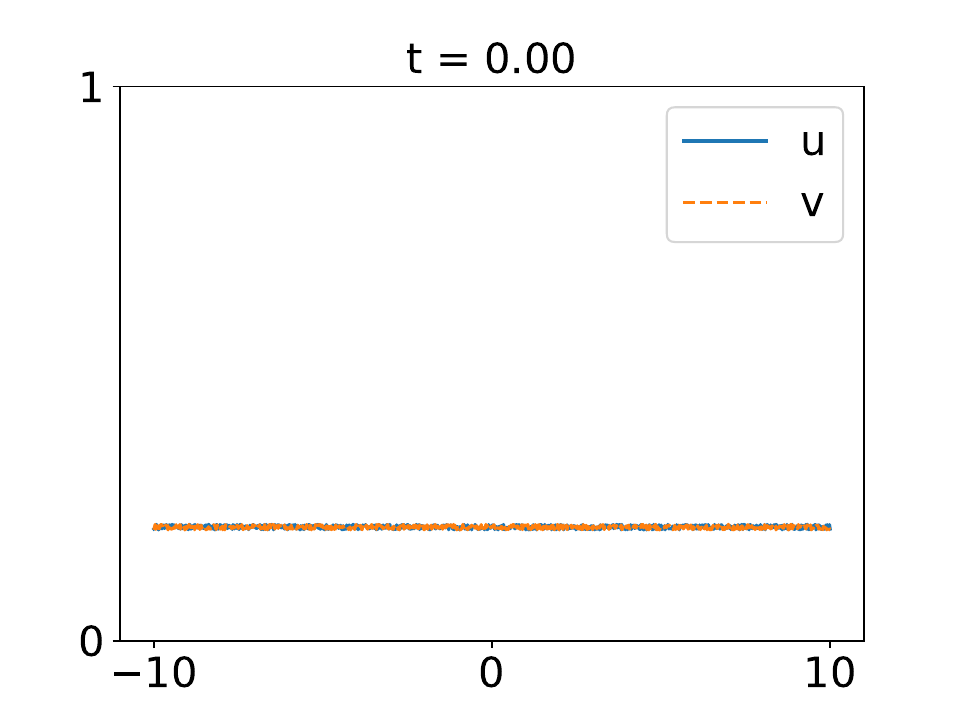}
     \end{subfigure}
     \begin{subfigure}
         \centering
         \includegraphics[width=0.235\textwidth]{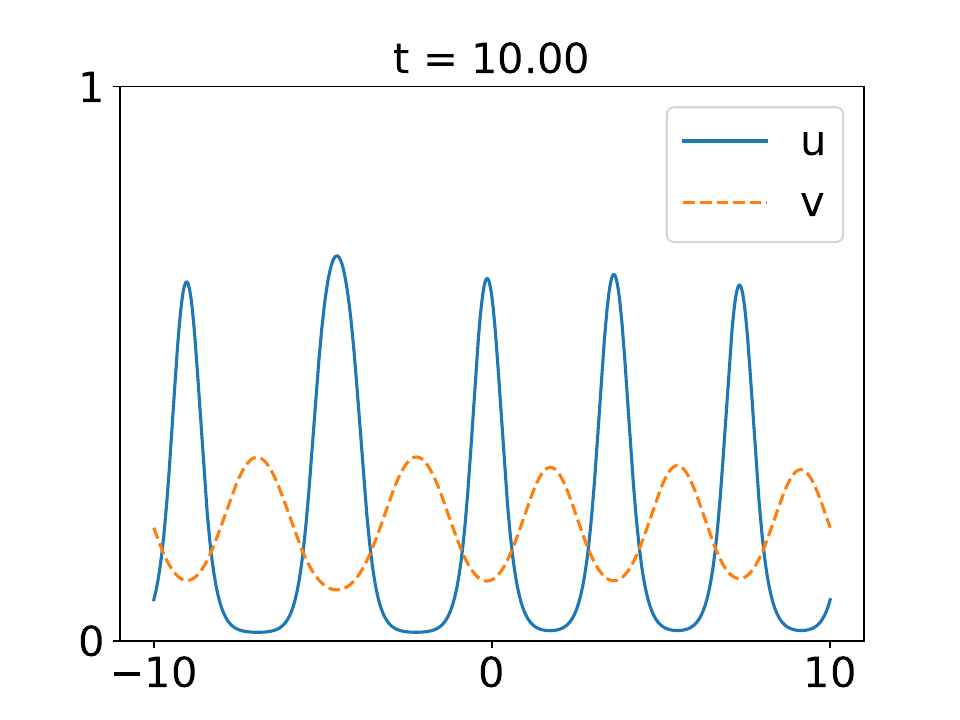}
     \end{subfigure}
     \begin{subfigure}
         \centering
         \includegraphics[width=0.235\textwidth]{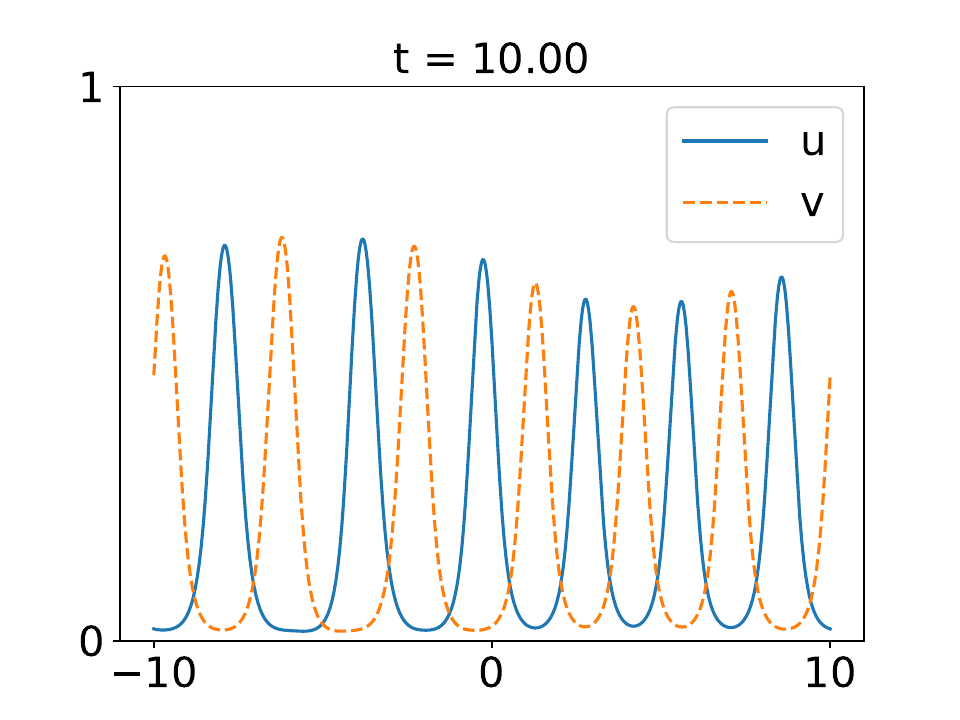}
     \end{subfigure}
     \begin{subfigure}
         \centering
         \includegraphics[width=0.235\textwidth]{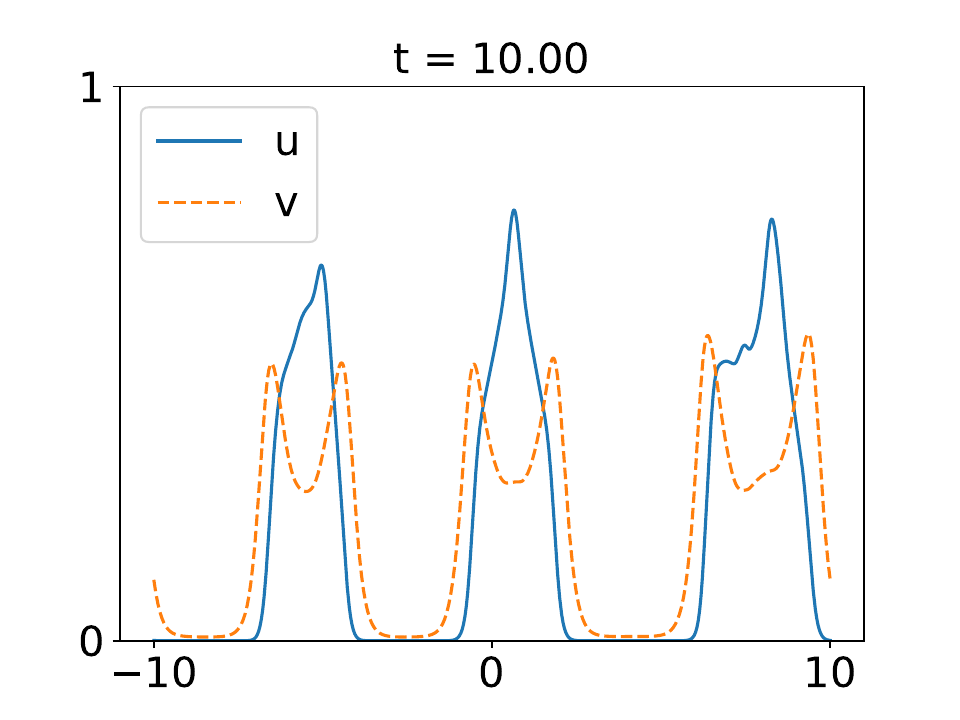}
     \end{subfigure}
     \begin{subfigure}
         \centering
         \includegraphics[width=0.235\textwidth]{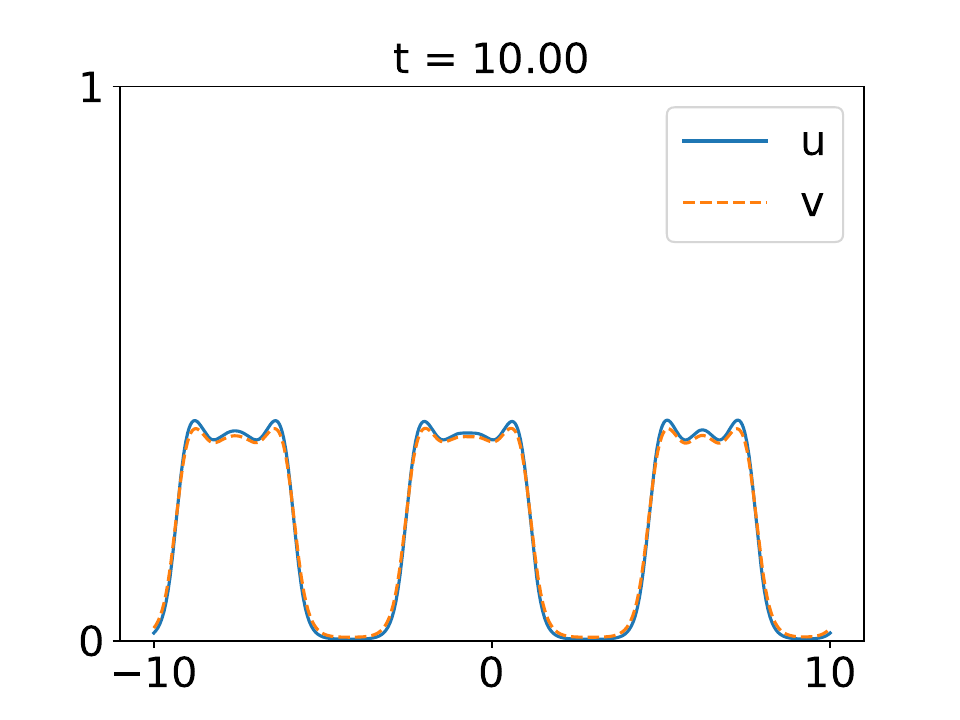}
     \end{subfigure}
     \begin{subfigure}
         \centering
         \includegraphics[width=0.235\textwidth]{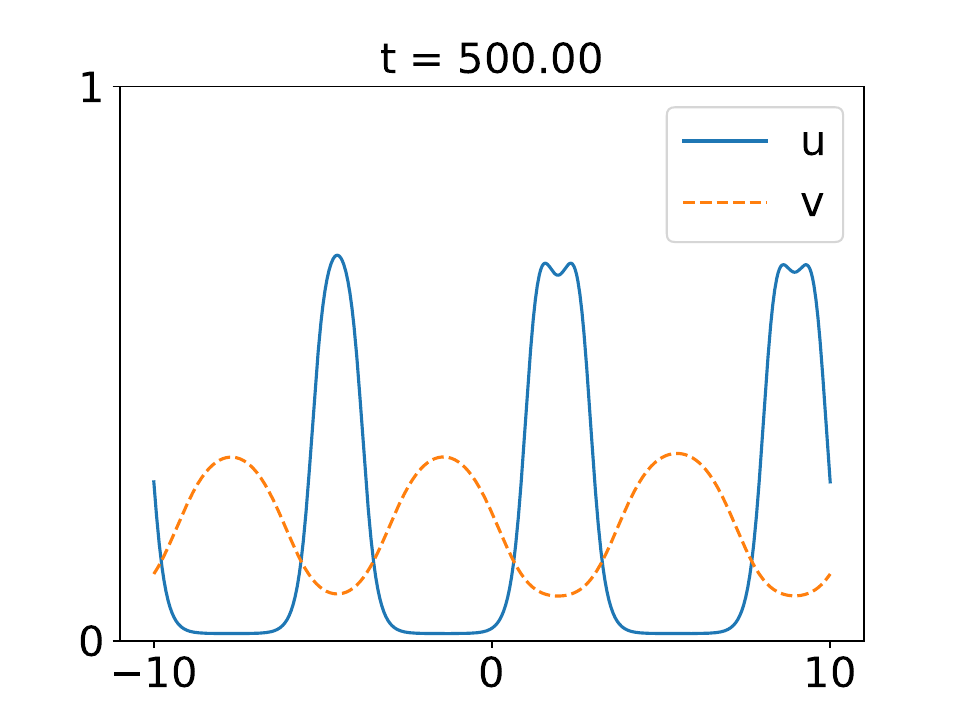}
     \end{subfigure}
     \begin{subfigure}
         \centering
         \includegraphics[width=0.235\textwidth]{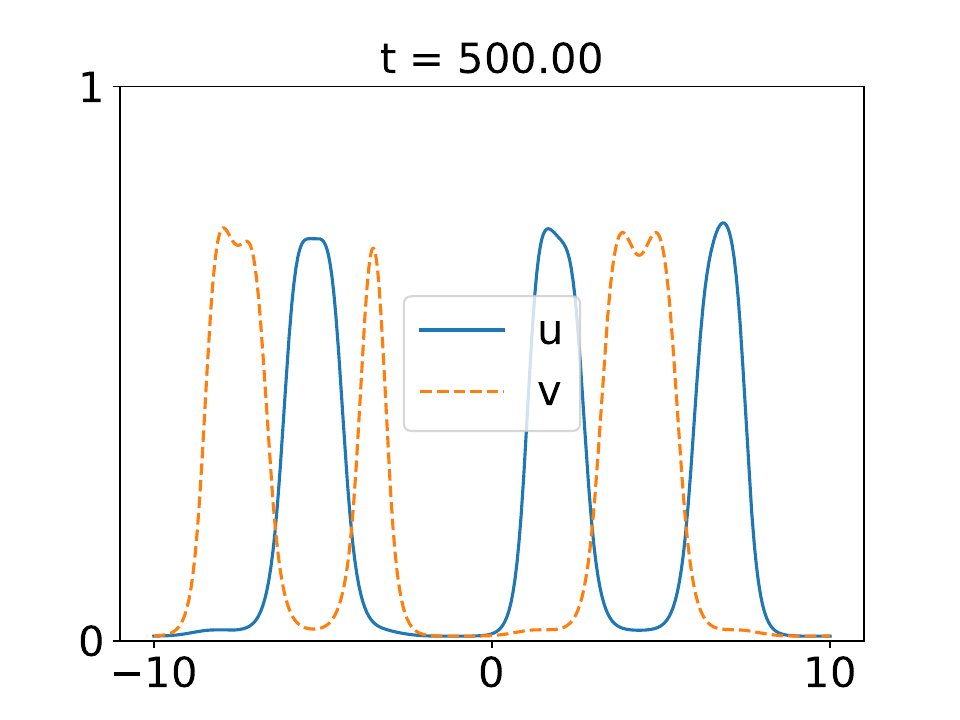}
     \end{subfigure}
     \begin{subfigure}
         \centering
         \includegraphics[width=0.235\textwidth]{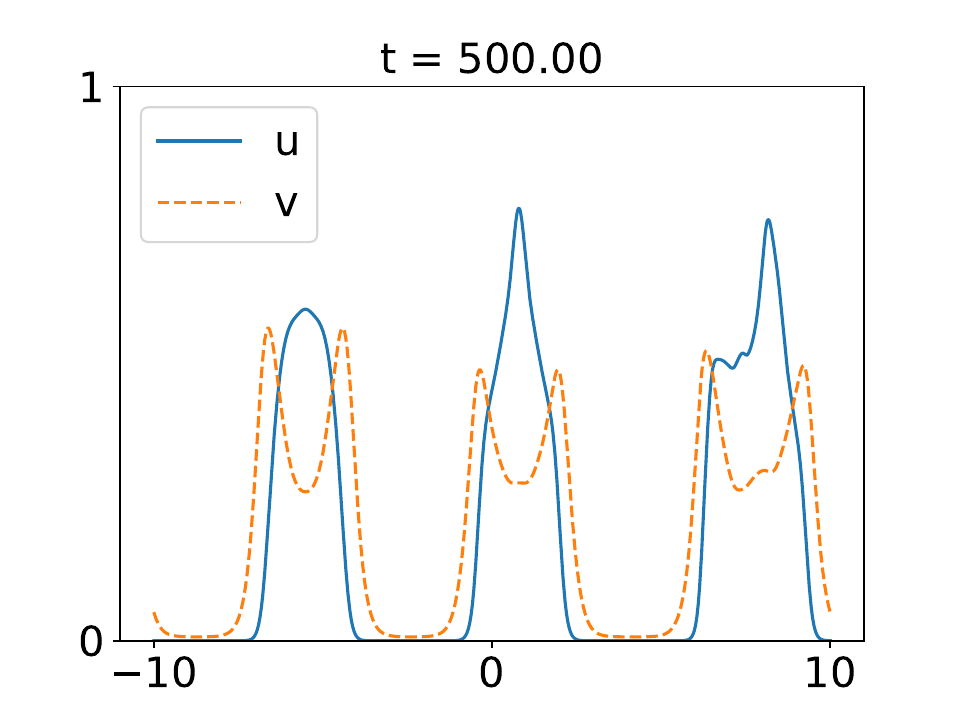}
     \end{subfigure}
     \begin{subfigure}
         \centering
         \includegraphics[width=0.235\textwidth]{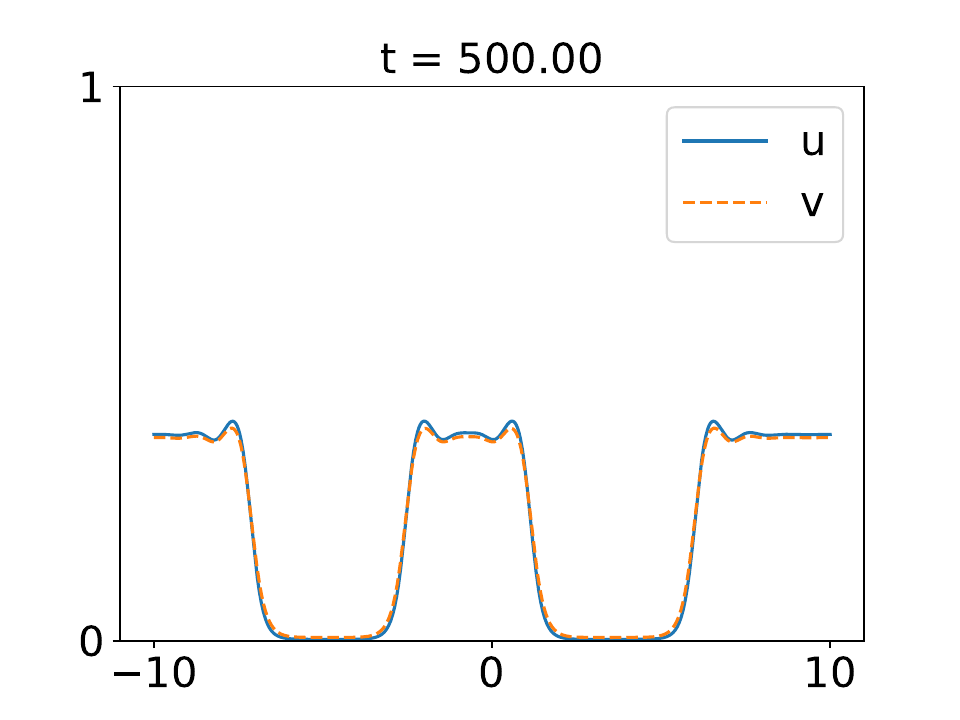}
     \end{subfigure}
     \caption{Numerical solutions of the uni-dimensional equation~\eqref{two population}, with FEM implemented in FEniCS, at different time. From left to right : First column $\Rightarrow$ $S_u=25$, $S_v=7.5$ and $C=0$, second column $\Rightarrow$ $S_u=25$, $S_v=25$ and $C=5$, third column $\Rightarrow$ $S_u=200$, $S_v=25$ and $C=50$, and last column $\Rightarrow$ $S_u=25$, $S_v=15$ and $C=30$. The rest of parameters are listed in Table~\ref{Tab:Param1dtwopop}.}
    \label{Fig:TwoPopulations1DResults}
\end{figure}
As expected, changing the values of the parameters $S_u$, $S_v$ and $C$ strongly impacts the behavior of the two populations over time. Furthermore, these results are in agreement with Steinberg's hypotheses on cell sorting, the Differential Adhesion Hypothesis~\cite{STEINBERG2007281, foty2004cadherin}, which states that:
\begin{align*}
    C=0 \rightarrow& \text{ Complet sorting};\\
    C < S_v  \text{ and } C<S_v \rightarrow& \text{ Partial engulfment of one population by the other};\\
    S_v<C<S_u \rightarrow& \text{ Engulfment of population $u$ by population $v$};\\
    \frac{S_u+S_v}{2} \leq C \rightarrow& \text{ Mixing}.
\end{align*}
Our simulations (Figure~\ref{Fig:TwoPopulations1DResults}) exactly reproduce these hypotheses: first column $\Rightarrow$ complete sorting, second column $\Rightarrow$ partial engulfment, third column $\Rightarrow$ engulfment of population $u$ by population $v$, and last column $\Rightarrow$ column mixing.

\paragraph{Bi-dimensional case.}
Analogous to what we have just done for the one-dimensional case, we implemented the finite element resolution, with FEniCS, of the two-dimensional problem with two populations of cells~\eqref{two population} coupled with the following initial condition:
\begin{align*}
    &u_0(x,y) = 0.1 + \text{rand}(0,10^{-2}),\\
    &v_0(x,y) = 0.1 + \text{rand}(0,10^{-2}).
\end{align*}
Moreover, functions $g_i$ for $i \in [uu, uv, vu, vv]$ have the same form as in the previous uni-dimensional case~\eqref{guu and gvu} and~\eqref{gvv and guv}.
\begin{longtable}{|p{1.75cm}|p{2.5cm}|p{8cm}|}
\caption{Description of 2D case model~\eqref{two population}-~\eqref{Kv} parameters with their values used for the numerical simulations presented in Figure~\ref{Fig:TwoPopulations2DResults}.} \label{Tab:Param2dtwopop} \\
\hline
Parameters & Value & Description \\ 
\hline \hline
$D$ & 1 & Diffusion coefficient\\
\hline
$S_u$ & 100 or 200 or 250 or 500 & Self-adhesive strength of population $u$ \\
\hline
$S_v$ & 50 or 100 or 200 & Self-adhesive strength of population $v$ \\
\hline
$C$ & 0 or 50 or 75 or 150 & Cross-adhesive strength between population $u$ and $v$\\
\hline
$L$ & 2.5 & $\Omega = [-L,L] \times [-L,L]$\\
\hline
$N$ & 80 & Number of intervals\\
\hline
$\tau$ & 0.01 & Step time\\
\hline
\end{longtable}
\begin{figure}[h!]
     \centering
     \begin{subfigure}
         \centering
         \includegraphics[width=0.235\textwidth]{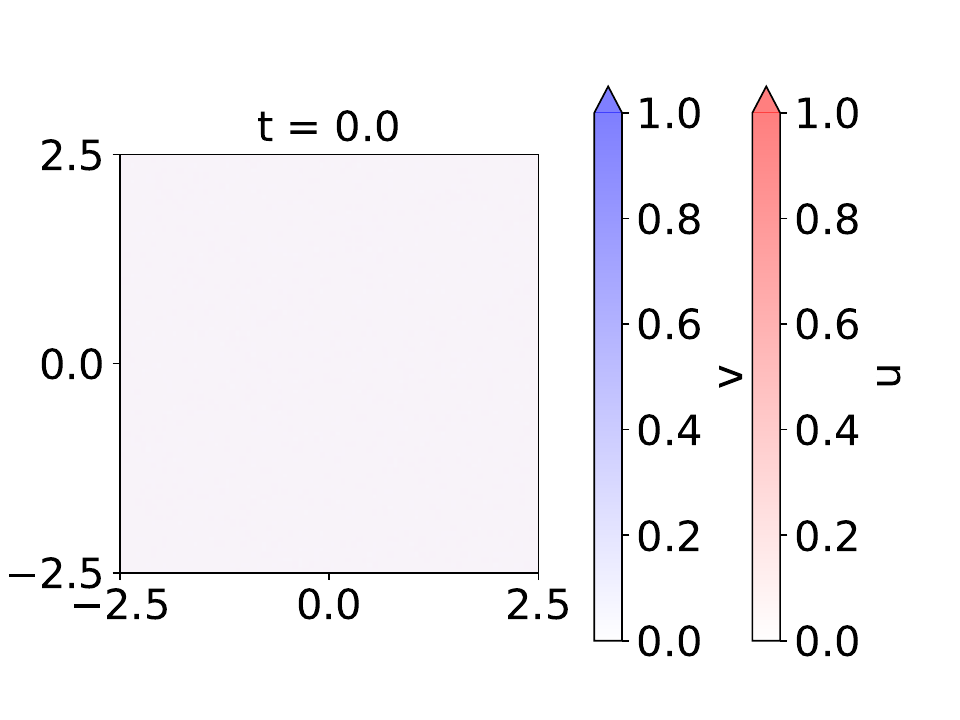}
     \end{subfigure}
     \begin{subfigure}
         \centering
         \includegraphics[width=0.235\textwidth]{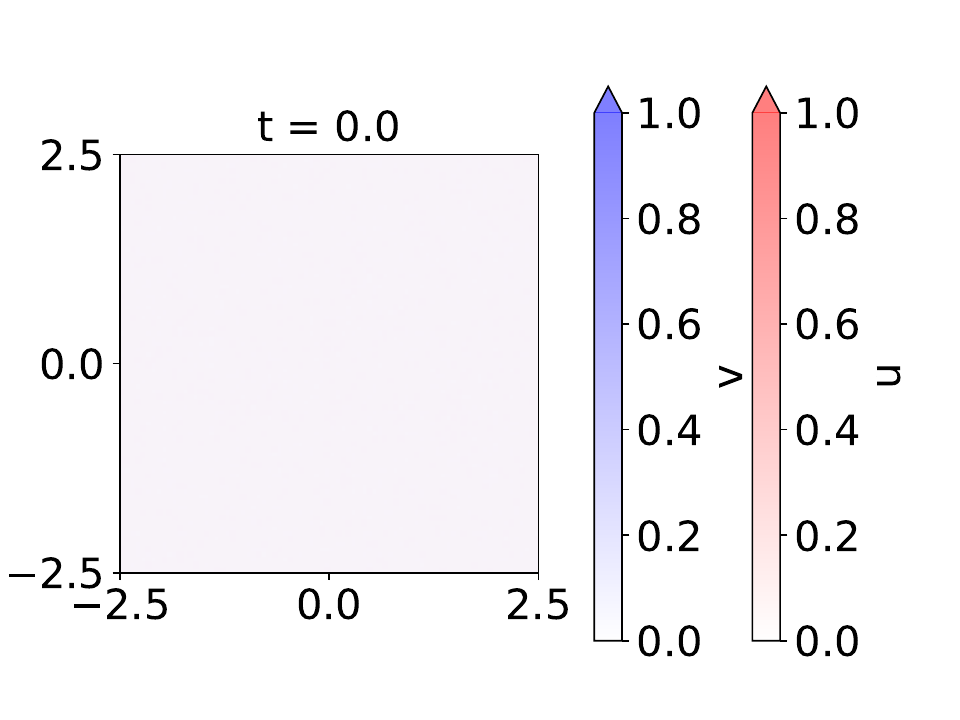}
     \end{subfigure}
     \begin{subfigure}
         \centering
         \includegraphics[width=0.235\textwidth]{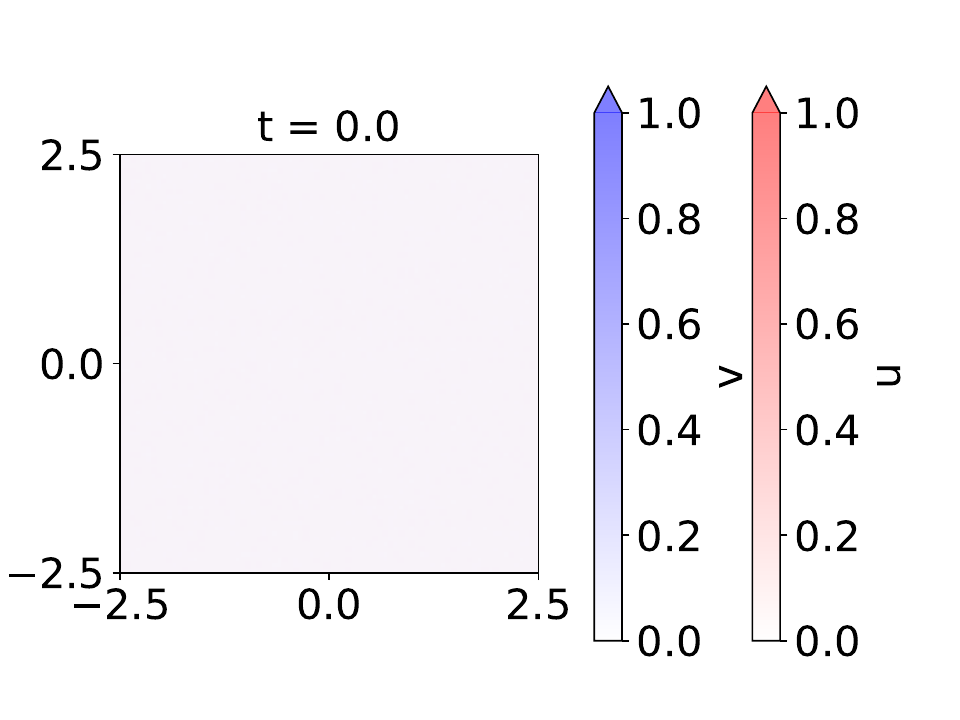}
     \end{subfigure}
     \begin{subfigure}
         \centering
         \includegraphics[width=0.235\textwidth]{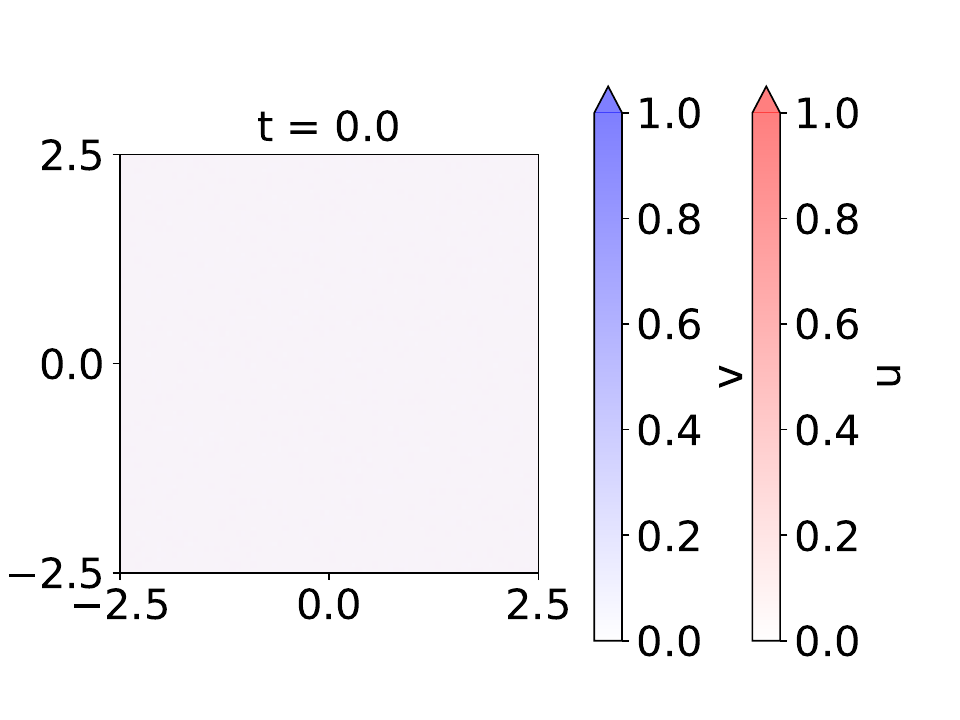}
     \end{subfigure}
     \begin{subfigure}
         \centering
         \includegraphics[width=0.235\textwidth]{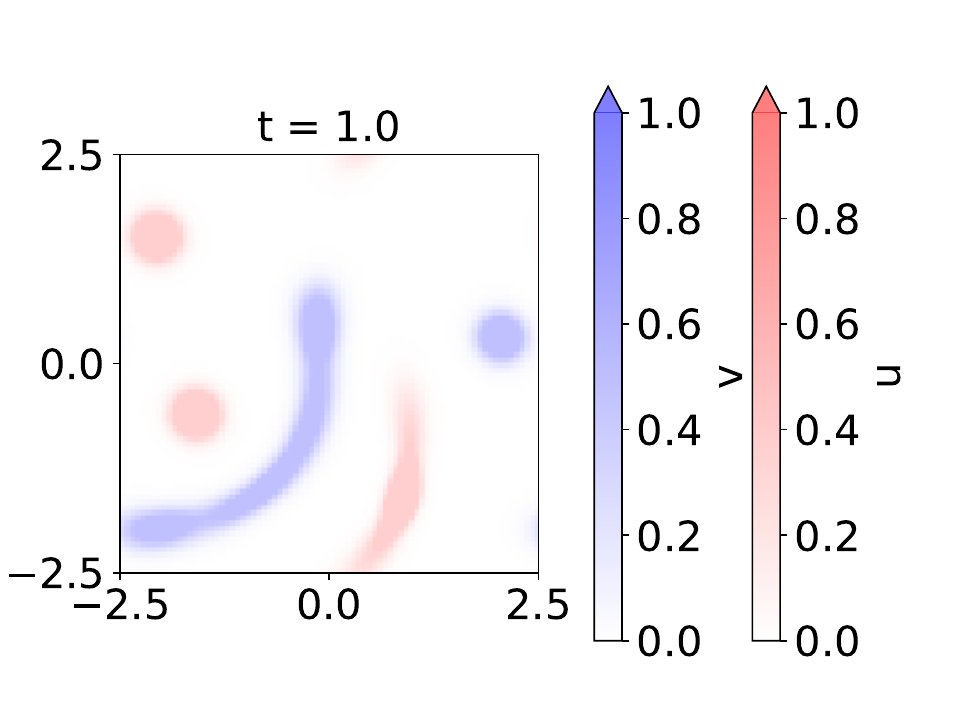}
     \end{subfigure}
     \begin{subfigure}
         \centering
         \includegraphics[width=0.235\textwidth]{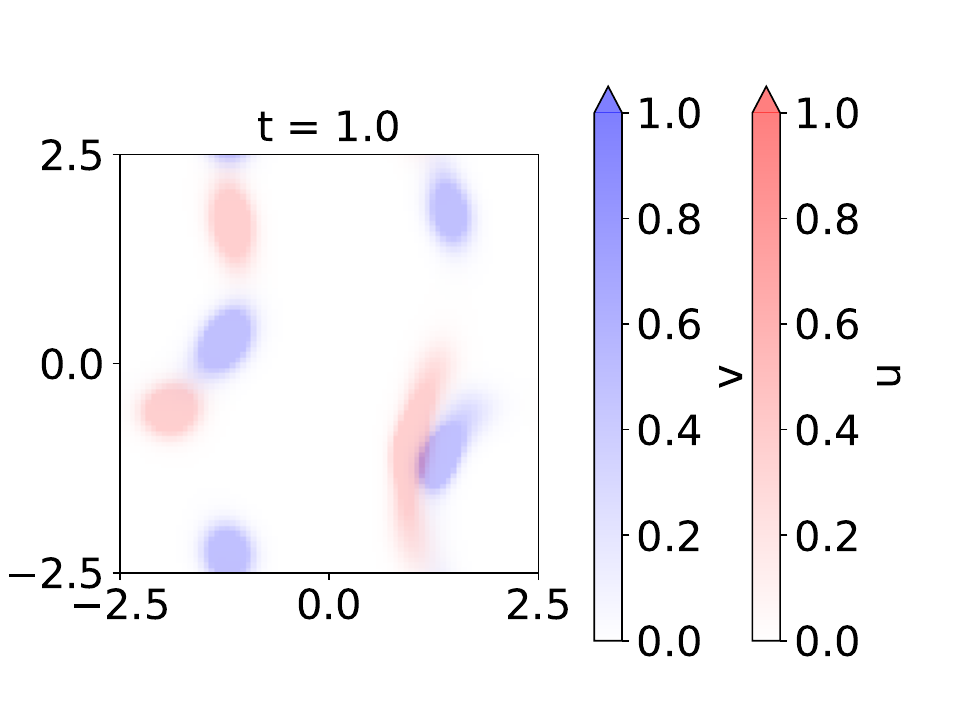}
     \end{subfigure}
     \begin{subfigure}
         \centering
         \includegraphics[width=0.235\textwidth]{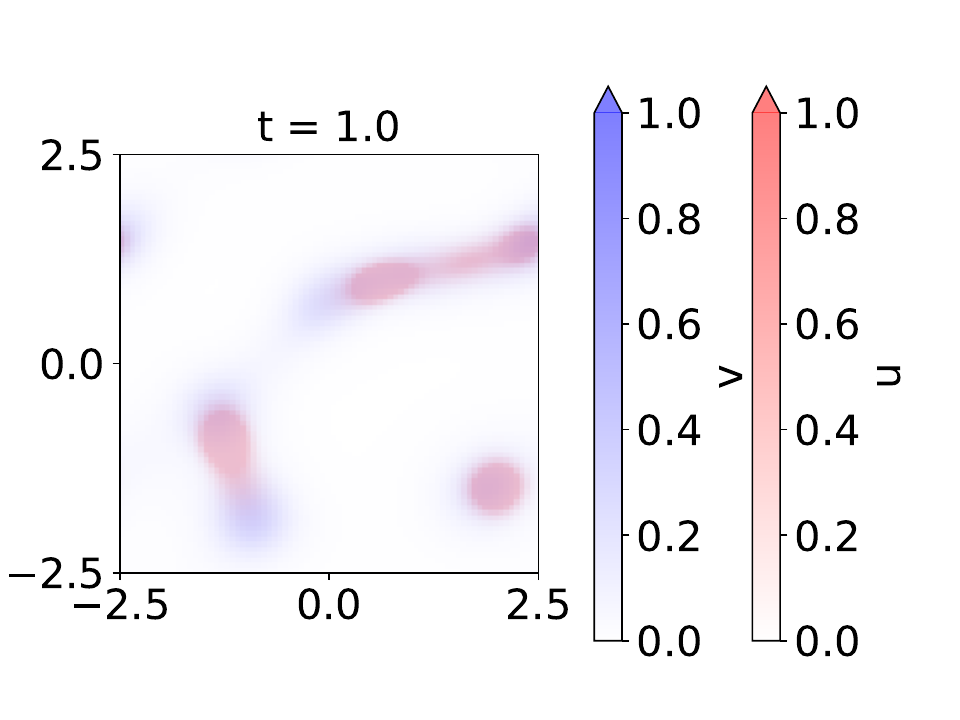}
     \end{subfigure}
     \begin{subfigure}
         \centering
         \includegraphics[width=0.235\textwidth]{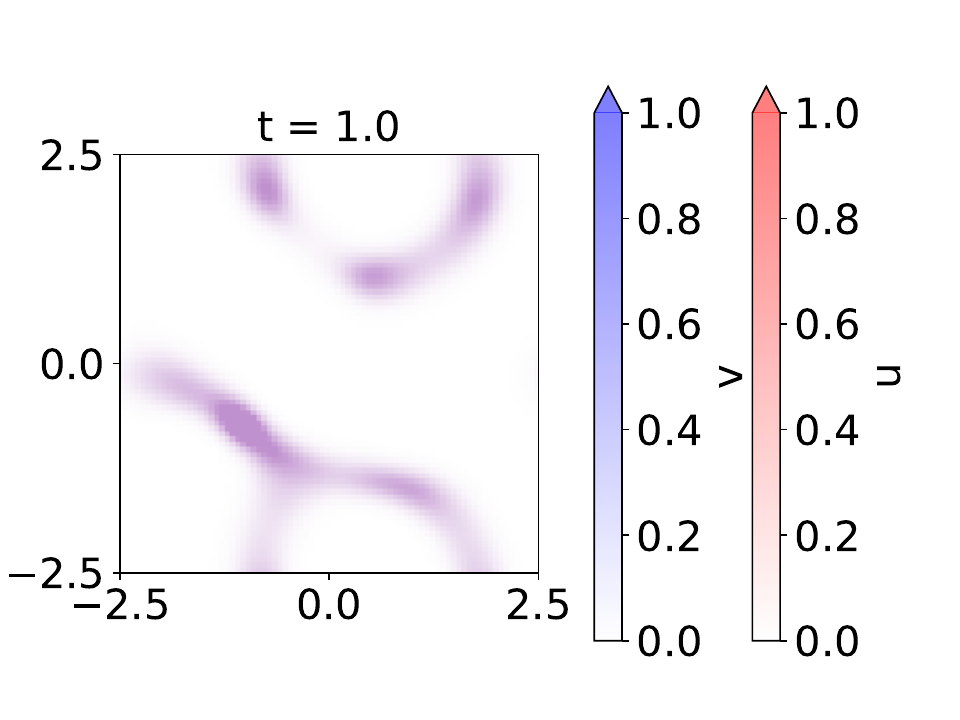}
     \end{subfigure}
     \begin{subfigure}
         \centering
         \includegraphics[width=0.235\textwidth]{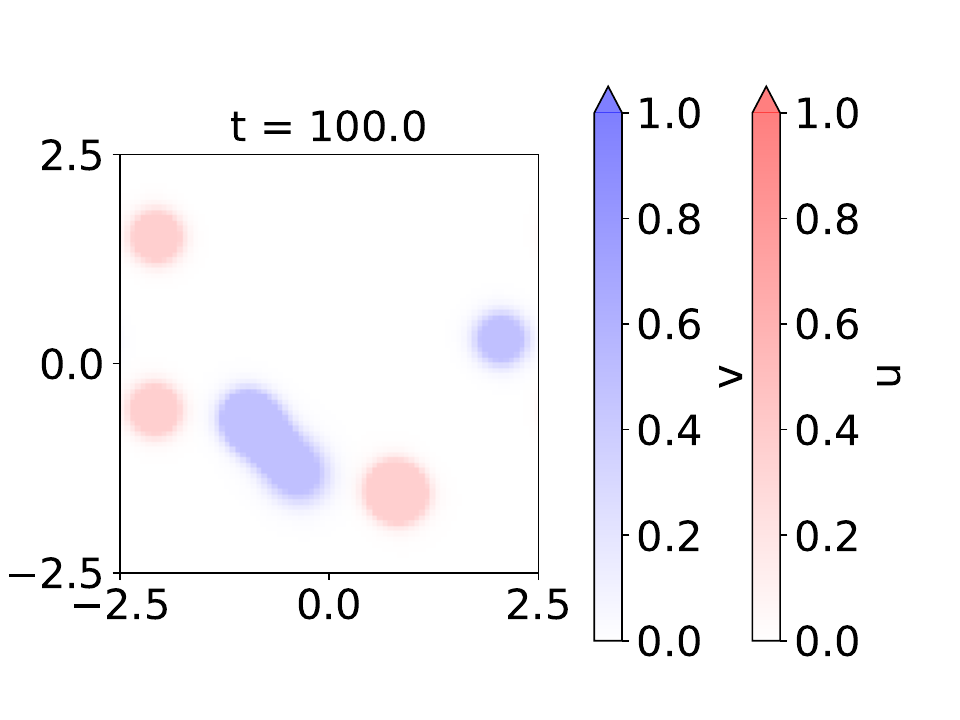}
     \end{subfigure}
     \begin{subfigure}
         \centering
         \includegraphics[width=0.235\textwidth]{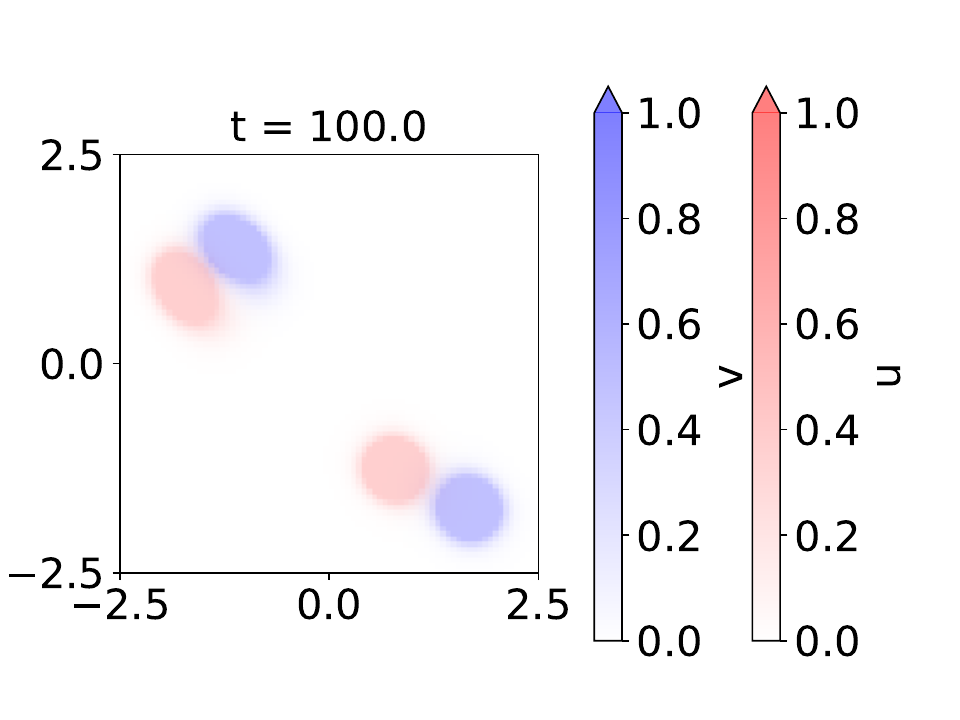}
     \end{subfigure}
     \begin{subfigure}
         \centering
         \includegraphics[width=0.235\textwidth]{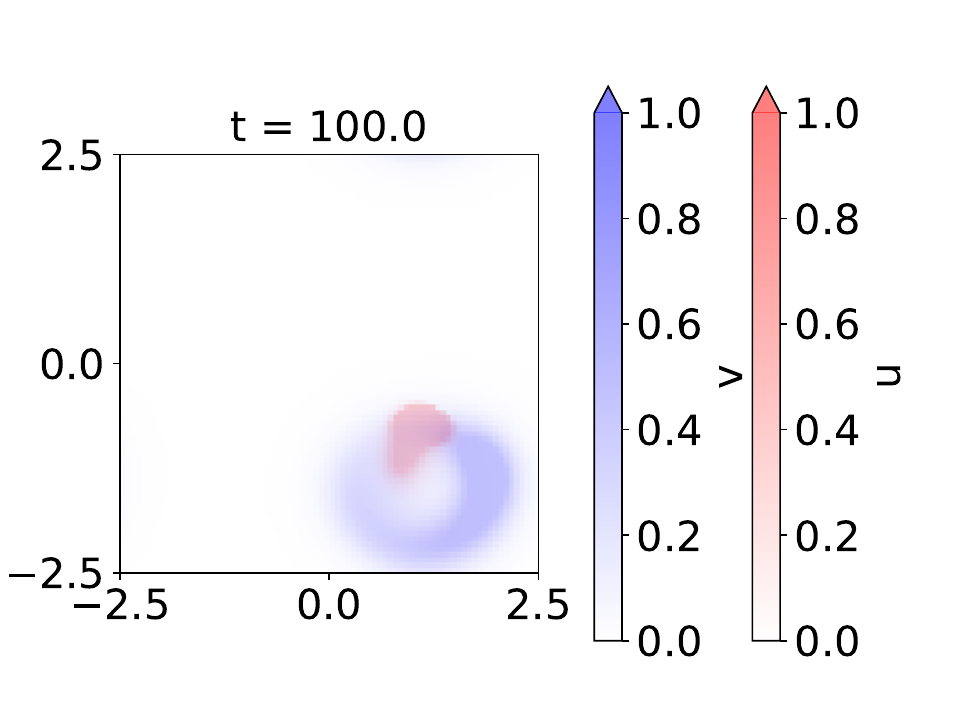}
     \end{subfigure}
     \begin{subfigure}
         \centering
         \includegraphics[width=0.235\textwidth]{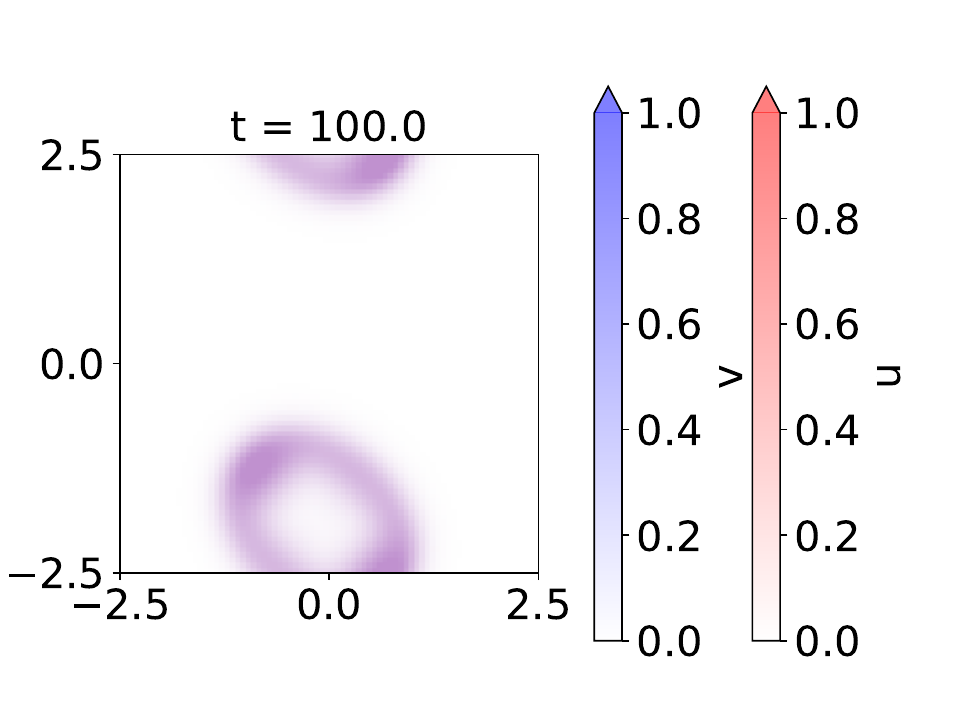}
     \end{subfigure}
     \caption{Numerical solutions of the bi-dimensional equation~\eqref{two population}, with FEM implemented in FEniCS, at different time. From left to right : First column $\Rightarrow$ $S_u=250$, $S_v=200$ and $C=0$, second column $\Rightarrow$ $S_u=200$, $S_v=200$ and $C=50$, third column $\Rightarrow$ $S_u=500$, $S_v=50$ and $C=75$, and last column $\Rightarrow$ $S_u=100$, $S_v=100$ and $C=150$. The rest of parameters are listed in Table~\ref{Tab:Param2dtwopop}.}
    \label{Fig:TwoPopulations2DResults}
\end{figure}
Here are the results obtained for different values of the parameters $S_u$, $S_v$ and $C$. These simulations (Figure~\ref{Fig:TwoPopulations2DResults}) once again confirm Steinberg's hypotheses~\cite{STEINBERG2007281, foty2004cadherin}:
\begin{itemize}
    \item First column: complet sorting because $C=0$
    \item Second column: partial engulfment because $C<S_u$ and $C<S_v$
    \item Third column: engulfment of population $u$ by population $v$ because $S_v < C < S_u$
    \item Last column: mixing because $\frac{S_u+S_v}{2} < C$
\end{itemize}

%%%%%%%%%%%%%%%%%%%%%%%%%%%%%%%%%%%%%%%%%%%%%%%%%%%%%%%%%%%%%%%%%%%%%%%%%%%%%%%%%%%%%%%%%%%%%%%%%%%%
\section{Open Numerical Problems}
\label{Sect:OpenProblems}

There are various open problems associated with the numerical implementation of 1D/2D/3D non-local models. In the following we discuss some of these issues.

%%%%%%%%%%%%%%%%%%%%%%%%%%%%%%%%%%%%%%%%%%%%%%%%%%%%%%%%%%
\subsection{In-house codes vs. open-source computing platforms}
\label{Subsect:In-houseVSOpen-source}

There are various open-source computational platforms that can be used to implement faster these non-local mathematical models: FEniCS~\cite{logg2012automated}, DUNE~\cite{bastian2008generic}, FreeFEM~\cite{hecht2012}, Chaste~\cite{cooper2020chaste}, etc. 
Each of these platforms have their own particularities. For this review study, we chose to focus on FEniCS (open source, that already has function implementing the finite element method) and on Python (in-house code, for which we implemented ourselves the finite element method -- see Appendix~\ref{Appendix:CalculationDetails}). Below we summarise some of the advantages and challenges we encountered when discretising and simulating numerically models~\eqref{nonlocal_model} and~\eqref{two population}.\\\\
Open source software (FEniCS in our case) generally works well for classic local problems. They are practical and powerful for this type of models. Indeed, these problems have already been, for the most part, implemented and it is therefore possible to find in the  documentation~\cite{LangtangenLogg2017} or the online forums (~\url{https://fenicsproject.discourse.group/}), the answers to the questions that one may ask. However, when we want to solve a more complicated problem, we may encounter implementation difficulties. This was the case for the numerical resolution of non-local problems~\eqref{nonlocal_model} and~\eqref{two population}. Indeed, on several points, the implementation of these non-local models posed some problems with the FEniCS software, because we needed to understand what is behind various functions of the main library \emph{dolfin}. Since FEniCS is
coded in Python with libraries from C++, understanding the implementation of the objects we are manipulating is not very intuitive. Let us present a concrete problem that we faced: 
when working with periodic conditions at the boundaries of the domain, the numbering of the degrees of freedom (DOFs) is no longer the same as when we use more classical conditions at the boundaries (e.g. Dirichlet, Neumann). In the one-dimensional case with classic conditions, the DOFs are numbered from left to right (node $x_0$ corresponds to DOF number 0, node $x_1$ to node 1, etc...). However, with periodic conditions, this order is no longer respected, the DOFs are reorganized in such a way as to make the solver more efficient. Thus, when we wanted to assign specific values to a FEniCS function (as was the case for the initial condition $u_0$ or the non-local term $K[u]$ at each time iteration), this posed a problem because the order of the vector coefficients no longer corresponded to the order of numbering of the DOFs. The wrong coefficient was then assigned to the DOF. We therefore had to find a function that reorganized the DOFs in the right order and then perform the assignment. But it was not easy to find as we had difficulty understanding in what order and why FEniCS did not organize the DOFs in logical order.\\\\
Unlike these open-source platforms, in-house codes (Python in our case) allow for more freedom and control of the discretization method chosen, and on the actual implementation of these methods. However, this requires having a fairly solid background in the different numerical methods because it is up to the person who implements to calculate by hand the coefficients of the matrices and vectors to be implemented to solve the problem. In the case of finite elements, the programmer must calculate the integrals himself and then implement them efficiently in matrices. These manual calculations can be a source of errors and the choice of implementation can be crucial if you do not want the program to take too long to run (e.g. use libraries that facilitate calculation with sparse matrices like \emph{scipy.sparse}). These problems specific to in-house codes do not concern open-source platforms. Indeed, when using these platforms, it is generally not necessary to have in-depth knowledge of the finite element method. We just need to know the variational formulation of our problem and then implement it. This therefore reduces the risk of errors in calculations and implementation. In addition, the speed of code compilation no longer depends on the programmer's skills, since the main calculations (e.g. calculations of integrals, matrix inversion) are already implemented in the software and the programmer does not have access to them.

Nevertheless, a problem that is common to open-source software and in-house codes is the efficient and inexpensive implementation of the non-local term at each mesh point and each iteration time. This is why throughout this article, we have seen two possible choices to approach the non-local term $K[u]$: see equations~\eqref{approximation trapezes} and~\eqref{approximation FFT}. The first method~\eqref{approximation trapezes}, which uses the trapezoid quadrature rule, is time-consuming. While with the second approximation~\eqref{approximation FFT}, is much faster since it uses the discrete Fourier transform, and more precisely an algorithm (Fast Fourier Transform) allowing the discrete Fourier transform of a vector to be calculated very efficiently. 
%But this method can only be used when the non-local term can be rewritten in the form of a convolution, which is not necessarily the case (see, for example, the models in~\cite{}).

Finally, we managed to exactly reproduce the results obtained with FEniCS (Figures~\ref{Fig:1D_FEniCS_results} and~\ref{Fig:2D_FEniCS_results}), by implementing the in-house Python code. This therefore shows that FEniCS software works correctly for non-local models.

%%%%%%%%%%%%%%%%%%%%%%%%%%%%%%%%%%%%%%%%%%%%%%
\subsection{Numerical oscillations}
\label{Subsect:NumericalOscillations}
\begin{figure}[h!]
     \centering
     \begin{subfigure}
         \centering
         \includegraphics[width=0.49\textwidth]{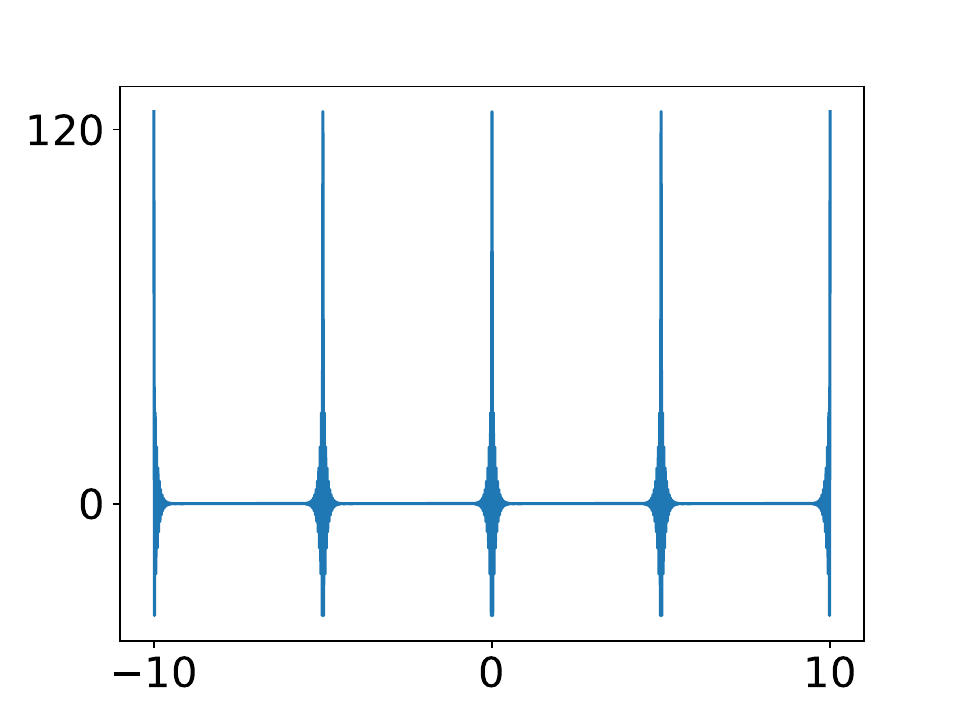}
     \end{subfigure}
     \begin{subfigure}
         \centering
         \includegraphics[width=0.49\textwidth]{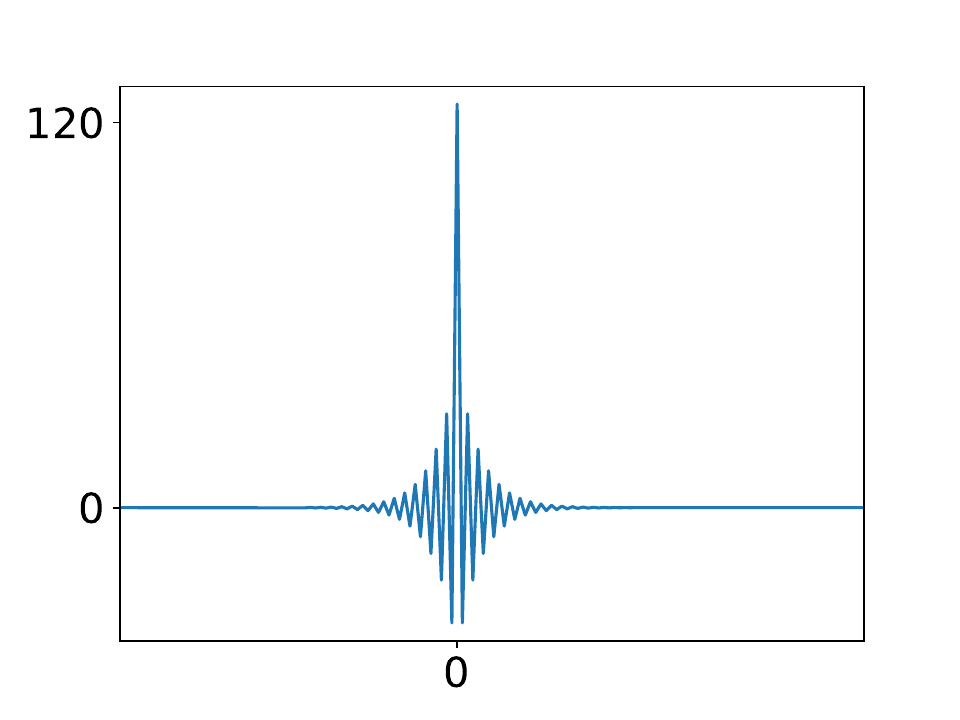}
     \end{subfigure}
     \caption{Finite element solution of the equation~\eqref{non_local model 1d} at time $t=1.0$ (left) and zoom on one of the interesting parts showing oscillations (right). Equation and mesh aparmeters are the same as in Table~\ref{Tab:Param1d}, except $\alpha=13$.}
    \label{Fig:Oscillations}
\end{figure}

It is well known that for the following local advection-diffusion equation:
\begin{equation*}
    \partial_t u = D \Delta u - \nabla \cdot (\beta u), 
\end{equation*}
numerical oscillations appear when we solve this equation numerically with a mesh step $h$ too large. Indeed, if $1<\frac{|\beta| h}{2D}$ oscillations appear and it is enough to refine the mesh to make them disappear~\cite{quarteroni2009numerical}. Thus, we wondered if such oscillations could appear with the non-local model~\eqref{nonlocal_model}. Indeed, by playing with the parameter $\alpha$ of the non-local form $K[u]$, oscillations appear (as we can see in Figure~\ref{Fig:Oscillations}) when the advective part dominates the diffusive part (i.e. when $D << \alpha$).
However, refining the mesh can become very time-consuming (even more for the non-local model than for the local model). Thus, stabilization methods were created for the local model to make these oscillations disappear without reducing the mesh step $h$: Discontinuous Galerkin~\cite{quarteroni2009numerical}, Streamline Upwind Petrov-Galerkin~\cite{quarteroni2009numerical}, Continuous Interior Penalty~\cite{burman2014robust}, etc. To our knowledge, there is no such methods for non-local models. Thus, we could try, in the near future, to adapt the existing stabilization methods for the local model to the non-local problem. Or even to completely create new stabilizing methods for such non-local models.

%%%%%%%%%%%%%%%%%%%%%%%%%%%%%%%%%%%%%%%%%%%%%%%%%
\subsection{Numerical continuation for non-local advection-dominated models}
\label{Subsect:NumContinuation}

To have a better understanding of how these aggregation patterns (shown in Figures~\ref{Fig:1D_FEniCS_results}-~\ref{Fig:TwoPopulations2DResults}) form, persist and sometimes disappear to give rise to other patterns, one can focus on numerical bifurcation theory, and in particular on numerical continuation approaches, which can trace solution branches as model parameters are varied~\cite{AllgowerGeorg1990,Rheinboldt2000_NumericalContinuation}. The general idea behind these algorithms is to use a predictor-corrector numerical scheme combined with singularity detection approaches to track the parametrised solution branches~\cite{Kuehn2014_GluingNumericalContinuation}. While there are various software packages for continuation problems arising in Ordinary Differential Equations (e.g., \texttt{AUTO}, \texttt{XPP}, \texttt{MatCont})
\cite{Kuehn2014_GluingNumericalContinuation}, this is not the case for Partial Differential Equations (PDEs), where most of the algorithms developed are applied to the specific problems investigated. One of the reasons for this lack of a general continuation software for PDEs is the fact that there are different spatial discretization methods for each class of PDEs (with specific spatial domain characteristics and boundary conditions)~\cite{Kuehn2014_GluingNumericalContinuation}. Also, the large number of degrees of freedom for spatial 2D and 3D PDE models also requires specific continuation algorithms~\cite{Uecker2021_ContinuationBifurcPDE}. 

In regard to PDEs, we note that for classical (local) PDE models one can use the \texttt{MATLAB} package \texttt{pde2path}, based on the finite element method (FEM) with $\mathbb{P}_{1}$ Lagrangian elements~\cite{Uecker2021_ContinuationBifurcPDE}. (We note that the FEM approaches are relatively common for these continuation algorithms due to the algebraic system that results from the discretization of the FEM model: $F(U,\lambda)=0$, where $\lambda$ is a model parameter and $F:\mathbb{R}^n\times \mathbb{R} \to \mathbb{R}$ a nonlinear function~\cite{Kuehn2014_GluingNumericalContinuation,Uecker2021_ContinuationBifurcPDE}.)
However, more complex PDEs, such as these non-local advection/transport-dominated models that we discussed here, pose various convergence problems~\cite{LeEftimie}. In general, for these non-local PDEs, specific in-house continuation algorithms have been developed for specific parameter ranges~\cite{buono2019kinetic}. 

For the future, we expect that the fast development of these non-local aggregation-diffusion models, and the various 2D and 3D patterns that they can generate, will lead to the development of new continuation algorithms for numerical bifurcation studies.

%%%%%%%%%%%%%%%%%%%%%%%%%%%%%%%%%%%%%%%%%%%
\subsection{Data availability and Parameter identification}
\label{Subsect:ParamData}

There are numerous studies that look at parameter and function identification for local models using various approaches: from Tikhonov regularisation~\cite{BraunNatalini2022} to Bayesian approaches~\cite{FalcoCohenCarrilloBaker2023,CampilloFunoletVenkataramanMadzvamuse2019}, and references therein. However, there are very few studies that try to identify parameters for non-local models of the type discussed here; among them we mention ~\cite{CarrilloEstradaRodriguez2024,app131910598,Maher2022}. For example, the study in \cite{app131910598}, focused on non-local models for ecological movement of animal populations inside home ranges. The authors first used synthetic data (the equilibrium solutions of the model) to explore parameter inference via maximum likelihood estimation and Stochastic Gradient Descent (SGD)~\cite{app131910598}. Then, they used location data for thirteen meerkat populations (from the Kalahari Meerkat Project), to perform a kernel density estimation to use initial population distributions based on this dataset.
The study in \cite{Maher2022} focused on local and non-local models for collective movement of cells during cancer invasions, and used Tikhonov regularisation approaches to identify the unknown mutation function (describing the transition from a ``normal" less invasive cancer cell sub-population to a ``mutated" more invasive sub-population). The authors in~\cite{Maher2022} used synthetic data (to which different levels of noise was added).
A very recent study~\cite{CarrilloEstradaRodriguez2024} proposed a Partial Inversion algorithm to estimate the non-local kernel (interaction potential) for a simple aggregation-diffusion equation, from synthetic noisy discrete trajectory data.

We emphasise that most of these studies use synthetic data for parameter/function identification. This is because spatial (and mainly spatio-temporal) data is not always available. However, recent advances in imaging techniques will improve this aspect through the generation of numerous datasets quantifying spatial interactions in both ecology \cite{Greer2021_ImagingZooplankton} and cell biology~\cite{Lewis2021_CancerBiologyImaging}.

\subsection{Complex domain geometries}\label{Sect:ComplexGeom}
The simulations presented in the previous sections were performed on simple rectangular domains. However, the power of finite element method relies in how it deals with complex geometries. And biological and medical applications of advection-dominated equations have considered for now only simple domains, due to the numerical difficulties associated with simulations on more complex geometries. But such complex geometries are very common in medical applications: the structure of tissues, the shape of organs, etc.. For this reason we focused this review on the finite element method: to emphasise the numerical possibilities of simulating such non-local models on more realistic 2D and 3D geometries.

%%%%%%%%%%%%%%%%%%%%%%%%%%%%%%%%%%%%%%%%%%%%%%%%%%%%%%%%%%%%%%%%%%%%%%%%%%%%%%%%%%%%%%%%%%%%%%%%%%%%
\section{Summary and discussion}
\label{Sect:Discussion}

The goal of this review paper was to summarise the current numerical approaches used to discretize and simulate numerically non-local aggregation PDEs, that are becoming more and more common in the literature as a result of more and more experimental studies that emphasise the importance of non-local interactions in cell biology and ecology.\\
Before discussing the different numerical approaches to solving non-local equations, we started this review by explaining the importance and origin of non-local models in biology and ecology. Then, we briefly reviewed the different existing analytical approaches for these non-local models: existence results, stability analysis, bifurcations, etc... We did not focus more on these analytical aspects as they were recently reviewed in~\cite{chen2020mathematical}.\\

Afterwards, we started the presentation of different numerical approaches that could be used to discretise such non-local models. For this part, we focused on a very simple and generic non-local model, that did not incorporate any of the biological details discussed in Section~\ref{Subsect:CellBiology}. The reason for using this simple model (for numerical illustration purposes) was two folded: (i) it allowed us to present clearly the various discretisation schemes (see equations \eqref{scheme FD}, \eqref{scheme VF}, \eqref{FE system}); (ii) it allowed us to compare our numerical results (obtained with these different schemes) with the simulations presented in~\cite{armstrong2006continuum}, for which we did not have any numerical details.
In regard to these numerical approaches, first we discussed two classical finite differences and finite volumes schemes, then a less common finite element scheme. We then presented a series of numerical simulations for different non-local models solved using the finite element method. We concluded this review with a discussion on open numerical problems related to these non-local models. In this part we first compared two ways of implementing the resolution of these non-local models via two
finite element schemes, one leading to numerical oscillations. Then, we discussed open problems related to numerical continuation algorithms for the numerical investigation of different bifurcation patterns exhibited by these non-local models. Finally we discussed briefly open problems related to numerical approaches that can be used to for parameter identification for such non-local models.

Through this review of numerical approaches for non-local advection-dominated problems, we aimed to emphasise the opportunities to further develop the finite element approaches to discretise and simulate these non-local models not only in 1D, but mainly in 2D and 3D, and on various geometries that are important in biological and medical applications (e.g., complex domains simulating various tissues and organs~\cite{magomedov2020application} on which non-local interactions among cells, molecules, etc., seem to gain interest~\cite{lee2001non}). By discussing (in Section~\ref{Sect:OpenProblems}) various open problems related to the numerical implementation of the finite element schemes, stabilisation approaches for these schemes, as well as the need to develop new schemes for numerical continuation methods (to investigate the bifurcation of the complex spatio-temporal patterns that can be exhibited by these models), we aimed to suggest some new numerical research directions that could advance the current research in the field.
This review complemented other recent reviews in the literature related to non-local models: from more modeling-focused and biologically-focused
reviews~\cite{lee2001non,painter2023biological,wang2023open} 
to more analytically-focused reviews~\cite{hillen2021non}, to reviews focused on numerical approaches for non-local models (spectral methods~\cite{hillen2021non}, finite volume methods~\cite{james2015numerical} and finite difference methods~\cite{delarue2020convergence}). Our overall goals here were (i) to show that finite element methods could also be used to discretise and simulate non-local advection-dominated models; (ii) to open this field of non-local models in biology to researchers focused on numerical analysis, to generate new advances in the field.

%Rather, it should relate the results to the existing literature; explain what was gained through this paper; talk about relevance for applications or math or numerics; and have a wider look at the field.

%%%%%%%%%%%%%%%%%%%%%%%%%%%%%%%%%%%%%%%%%%%%%%%%%%%%%%%%%%%%%%%%%%%%%%%%%%%%%%%%%%%%%%%%%%%%%%%%%%%%
\section*{Acknowledgements}
This work (JM, AL, RE) has been supported by EIPHI Graduate School (contract ANR-17-EURE-0002) and by the  Bourgogne-Franche-Comté Region. 
RE and AL work was also partially supported by the ANR grant ANR-21-CE45-0025-01. RE and JM work was also partially supported by the Région Bourgogne Franche-Comté “Accueil de Nouvelle Équipe de Recherche (ANER) 2022” grant number 2022-7-13828.

\section*{Conflict of interest}
All authors have no conflicts of interest.

%%%%%%%%%%%%%%%%%%%%%%%%%%%%%%%%%%%%%%%%%%%%%%%%%%%%%%%%%%%%%%%%%%%%%%%%%%%%%%%%%%%%%%%%%%%%%%%%%%%%
\begin{appendices}

\section{Calculation details of FEM}\label{Appendix:CalculationDetails} 

\begin{figure}[h!]
    \centering
    \includegraphics[width=\linewidth]{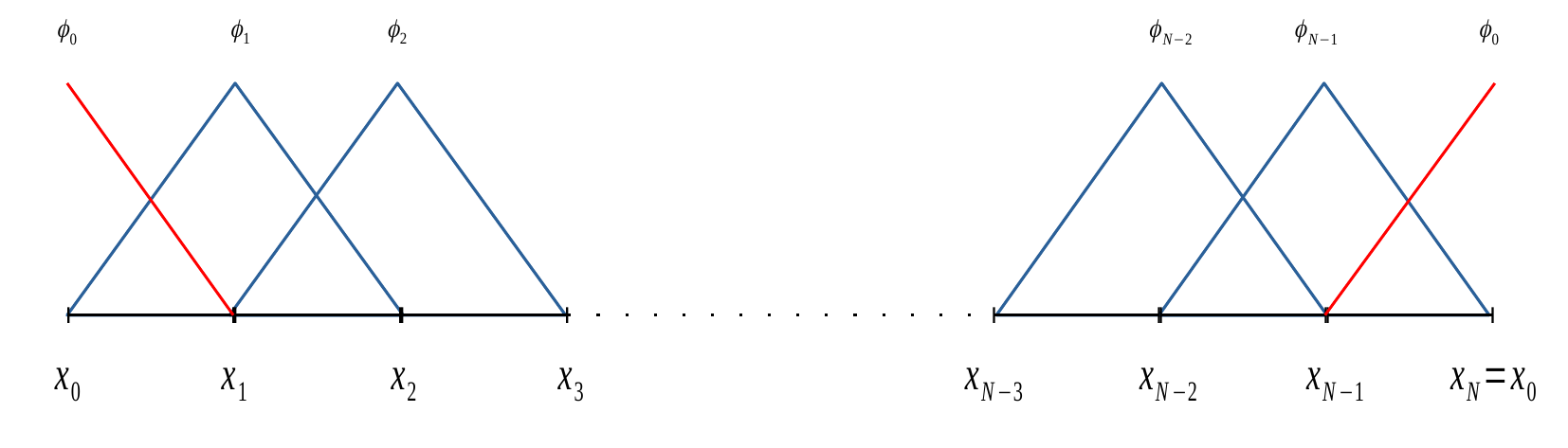}
    \caption{The basis functions $\phi_i$ associated to node $x_i$ with $i \in \{0,...,N-1\}$. The basic function $\phi_0$ is drawn in red since it is this which translates the periodic conditions at domain boundaries.}
    \label{Fig:basic functions}
\end{figure}

Given the discrete variational formulation~\eqref{variational form}, we can compute the matrix $A^n_h$ and the right-hand side $B^n_h$, such that:
\begin{align*}
        &A^n_{i,j} = a(\phi_i,\phi_j;u_h^n),\\
        &B^n_i = l(u_h^n;\phi_i),
\end{align*}
with $\phi_i$ for $i=1,...,N-1$  are the classical basis functions (see Figure~\ref{Fig:basic functions}):
\begin{equation*}
    \phi_i(x)=
    \begin{cases}
        \frac{x-x_{i-1}}{h}, \text{ if } x \in [x_{i-1},x_i ],\\
        \frac{x_{i+1}-x}{h}, \text{ if } x \in [x_i,x_{i+1}],\\
        0, \text{ otherwise},  
    \end{cases}
    \text{ and } \phi_0(x) = 
    \begin{cases}
        \frac{x-x_{N-1}}{h}, \text{ if } x \in [x_{N-1},x_N ],\\
        \frac{x_{1}-x}{h}, \text{ if } x \in [x_0,x_{1}],\\
        0, \text{ otherwise}.  
    \end{cases}
\end{equation*}
Thus, to find the solution to the problem~\eqref{variational form} at each iteration in time, it suffices to solve the matrix system $A_h U_h^{n+1} =B_h^n$. \\
To solve the previous system, we have to make several approximations: first we apply some quadrature rules to calculate the integral
\begin{equation*}
    \int_\Omega \phi_i(x)K[u^n](x,t) \frac{\partial \phi_j(x)}{\partial x} \ dx,
\end{equation*}
because we cannot calculate its exact value. Then, we need to calculate the value of $K(u^n)$ at each point of the mesh to be able to apply the quadrature rule. To approximate the value of the non-local term, we could have used the same approximation as in Section~\ref{Subsect:FiniteDiff}, but we chose to try to find another method, since this non-local term seems to be similar to a convolution. Hence, we first rewrite it in the form of a convolution:
\begin{equation*}
    K[u](x,t)= -\alpha \int_{-L}^L u(x-y,t)\tilde{\omega}(y)=- \alpha (u \ast \tilde{\omega})(x),
\end{equation*}
with
\begin{equation}
    \tilde{\omega}(x) = 
    \begin{cases}
        sgn(x), \quad \text{if } x \in [-r,r]\\
        0, \quad \text{else}.
    \end{cases}
\end{equation}
Under this form, we can therefore use the \textit{Discret Fourier Transform} to calculate the non-local term $K[u]$. First, we can discretize the convolution with the trapezoidal rule:
\begin{align*}
&K[u](x,t) \approx - \alpha h \sum_{k=0}^{N-1} \frac{u(x-x_{k+1},t) \tilde{\omega}(x_{k+1}) + u(x-x_k,t) \tilde{\omega}(x_k)}{2} \\
&= - \alpha h \biggl( \frac{u(x-x_0,t) \tilde{\omega}(x_0) + u(x-x_{N},t) \tilde{\omega}(x_{N})}{2} +\sum_{k=1}^{N-1} u(x-x_k,t) \tilde{\omega}(x_k) \biggr).
\end{align*}
Next, thanks to the periodicity of $u$ and $\tilde{\omega}$, we have:
\begin{align*}
    K[u](x,t) &= - \alpha h \biggl( u(x-x_0,t) \tilde{\omega}(x_0) + \sum_{k=1}^{N-1} u(x-x_k,t) \tilde{\omega}(x_k) \biggr) \\
    &= - \alpha h \sum_{k=0}^{N-1} u(x-x_k,t) \tilde{\omega}(x_k).
\end{align*}
So, at any point of the mesh,
\begin{align*}
K[u](x_j,t)&= - \alpha h \sum_{k=0}^{N-1} u(x_j-x_k,t) \tilde{\omega}(x_k)\\
&= - \alpha h \sum_{k=0}^{N-1} u(-L+jh -(-L+kh),t) \tilde{\omega}(x_k)\\
&= - \alpha h \sum_{k=0}^{N-1} u(L+x_{j-k},t) \tilde{\omega}(x_k).
\end{align*}
Hence, to obtain the value of $K[u]$ at any point of the mesh, we have to calculate $u$ at the points $x_i + L$.
In practice, thanks to the periodicity of $u$ and $K[u]$, we calculate $K[u](x_{n}+L,t)= \alpha h \sum_{k=0}^{N-1} u(x_{n-k},t) \tilde{\omega}(x_k )$, then we translate all the values from L. \\
So, $\forall  j=0,...,N-1$ we calculate:
\begin{equation*}
    (K[u])^n_j =  - \alpha h  \sum_{k=0}^{N-1} u^n_{j-k} \tilde{\omega}_k,
\end{equation*}
where $\tilde{\omega}_k = \tilde{\omega}(x_k)$ for $k=0,...,N-1$.\\
Knowing that $u$ and $\tilde{\omega}$ are N-periodic, we use the following DFT:
\begin{equation*}
(\mathcal{F}(u^n))_j = \sum_{k=0}^{N-1} u_k^n e^{- \frac{2 \pi i}{N} kj}.
\end{equation*}
Hence, 
\begin{align*}
(\mathcal{F}(K[u]^n))_m &= \sum_{s=0}^{N-1} (K[u])^n_s e^{- \frac{2 \pi i}{N} sm}\\
&= - \alpha h \sum_{s=0}^{N-1}\sum_{k=0}^{N-1}u_{s-k}^n \tilde{\omega}_k e^{- \frac{2 \pi i}{N} sm}\\
&= - \alpha h \sum_{k=0}^{N-1} e^{-\frac{2 \pi i}{N}km} \tilde{\omega}_k \sum_{s=0}^{N-1} e^{-\frac{2 \pi i}{N}(s-k)m} u_{s-k}^n \\
&= - \alpha h \mathcal{F}(u^n)_m \mathcal{F}(\tilde{\omega})_m.
\end{align*}
And so,
\begin{equation}
    K^n_i := K(u_h)(x_i,t_n) \approx - \alpha h \mathcal{F}^{-1}(\mathcal{F}(U^n_h)\mathcal{F}(\tilde{\omega})).
    \label{approximation FFT}
\end{equation}
We can now calculate the coefficients of the matrix $A_h^n$, by using Simpson's rule when is necessary, at each iteration in time. Thanks to the form of the basis functions $\phi_i$, we already know which coefficients of the matrix $A^n_h$ will be non-zeros. In fact, we have (see Figure~\ref{Fig:basic functions}):
\begin{equation*}
    \begin{cases}
        \text{supp}(\phi_0)=[x_0,x_1] \cup [x_{N-1},x_0],\\
        \text{supp}(\phi_i)=[x_{i-1},x_{i+1}], \ \forall i \in \{1,...,N-2 \},\\
        \text{supp}(\phi_{N-1}) = [x_{N-2},x_{N-1}] \cup [x_{N-1},x_0],
    \end{cases}
\end{equation*}
thus the non-zero coefficients are: 
\begin{equation*}
    \begin{cases}
        A^n_{0,1}, \ A^n_{0,0}, \ A^n_{0,N-1},\\
        A^n_{i,i+1}, \ A^n_{i,i}, \ A^n_{i,i-1}, \ \forall i \in \{1,...,N-2\},\\
        A^n_{N-1,0}, \ A^n_{N-1,N-1}, \ A^n_{N-1,0},
    \end{cases}
\end{equation*}
so the matrix $A^n_h$ is almost tridiagonal. Hence, let's start by calculating the coefficients of the diagonal. For $i=1,...,N-2$:
\begin{align*}
    A^n_{i,i} &= a(\phi_i,\phi_i,u^n_h)\\
    &= \int^{x_i}_{x_{i-1}} [ (\phi_i(x))^2 +D \tau (\frac{\partial \phi_i(x)}{\partial x})^2 - \tau \phi_i(x) K(u^n_h)(x) \frac{\partial \phi_i(x)}{\partial x} ] \ dx\\
    & \quad + \int_{x_i}^{x_{i+1}} [ (\phi_i(x))^2 +D \tau (\frac{\partial \phi_i(x)}{\partial x})^2 - \tau \phi_i(x) K(u^n_h)(x) \frac{\partial \phi_i(x)}{\partial x} ] \ dx\\
    &\approx \frac{1}{h^2} \biggl( \biggl[\frac{x^3}{3} - x^2x_{i-1} + xx_{i-1}^2 \biggr]^{x_i}_{x_{i-1}} + D \tau h  - \frac{\tau h^2}{6} (2K^n_i + K^n_{i-1})\\
    & \quad + \biggl[\frac{x^3}{3} - x^2x_{i+1} + xx_{i+1}^2 \biggr]_{x_i}^{x_{i+1}} + D \tau h  + \frac{\tau h^2}{6} (2K^n_i + K^n_{i+1}) \biggr)\\
    &= \frac{2h}{3} + \frac{2D \tau }{h} + \frac{\tau}{6}(K^n_{i+1} - K^n_{i-1}).
\end{align*}
And with a similar calculation,
\begin{align*}
    &A^n_{0,0} = \frac{2h}{3} + \frac{2D \tau }{h} + \frac{\tau}{6}(K^n_{1} - K^n_{N-1}),\\
    &A^n_{N-1,N-1} = \frac{2h}{3} + \frac{2D \tau }{h} + \frac{\tau}{6}(K^n_{0} - K^n_{N-2}).
\end{align*}
Now, we calculate the coefficients of the superior diagonal, for $i=0,...,N-2$:
\begin{align*}
    A^n_{i,i+1} &= a(\phi_i,\phi_{i+1},u^n_h)\\
    &= \int_{x_i}^{x_{i+1}} [ \phi_i(x)\phi_{i+1}(x) +D \tau \frac{\partial \phi_i(x)}{\partial x} \frac{\partial \phi_{i+1}(x)}{\partial x} - \tau \phi_i(x) K(u^n_h)(x) \frac{\partial \phi_{i+1}(x)}{\partial x} ] \ dx\\
    &\approx \frac{1}{h^2} \biggl( \biggl[ \frac{x^2}{2}(x_{i+1}+x_i) - xx_i x_{i+1} - \frac{x^3}{3} \biggr]_{x_{i-1}}^{x_i} - D \tau h + \frac{\tau h^2}{6} (2K^n_i + K^n_{i+1})  \biggl)\\
    &= \frac{h}{6} - \frac{D \tau}{h} + \frac{\tau}{6}(2K^n_i + K^n_{i+1}),
\end{align*}
and 
\begin{equation*}
    A^n_{N-1,0} = \frac{h}{6} - \frac{D \tau}{h} + \frac{\tau}{6}(2K^n_{N-1} + K^n_{0}).
\end{equation*}
The same for the lower diagonal, for $i=1,...,N-1$:
\begin{align*}
    A^n_{i,i-1} &= a(\phi_i,\phi_{i-1},u^n_h)\\
    &= \int^{x_i}_{x_{i-1}} [ \phi_i(x)\phi_{i-1}(x) +D \tau \frac{\partial \phi_i(x)}{\partial x} \frac{\partial \phi_{i-1}(x)}{\partial x} - \tau \phi_i(x) K(u^n_h)(x) \frac{\partial \phi_{i-1}(x)}{\partial x} ] \ dx\\
    &\approx \frac{1}{h^2} \biggl( \biggl[ \frac{x^2}{2}(x_{i-1}+x_i) - xx_i x_{i-1} - \frac{x^3}{3} \biggr]^{x_{i+1}}_{x_i} - D \tau h - \frac{\tau h^2}{6} (2K^n_i + K^n_{i-1})  \biggl)\\
    &= \frac{h}{6} - \frac{D \tau}{h} - \frac{\tau}{6}(2K^n_i + K^n_{i-1}),
\end{align*}
and 
\begin{equation*}
    A^n_{0,N-1} = \frac{h}{6} - \frac{D \tau}{h} - \frac{\tau}{6}(2K^n_{0} + K^n_{N-1}).
\end{equation*}
To complete our systeme, we have to calculate the right hand side vector $B_h^n$. We use, for any interior node of the mesh, Simpson's quadrature rule (but we could also have used the trapezoid method):
\begin{align*}
    B^n_i &\approx \int_\Omega u^n_h(x) \phi_i(x) \ dx \\
    &= \frac{1}{h} \biggl[ \int^{x_{i+1}}_{x_i} u^n_h(x) (x_{i+1}-x) \ dx + \int_{x_{i-1}}^{x_i} u^n_h(x) (x - x_{i-1}) \ dx \biggr]\\
    &\approx \frac{h}{6} \biggl[ u^n_{i-1} + 4 u^n_i + u^n_{i+1} \biggr],
\end{align*}
and 
\begin{align*}
    B^n_0 &\approx \frac{h}{6} \biggl[ u^n_{N-1} + 4 u^n_0 + u^n_{1}\biggr],\\
    B^n_{N-1} &\approx \frac{h}{6} \biggl[ u^n_{N-2} + 4 u^n_{N-1} + u^n_{0} \biggr].
\end{align*}
Now that we know the form of the matrix $A^n_h$ and the vector $B^n_h$, we rewrite the matrix system $A^n_h U_h = B^n_h$ as the system~\eqref{FE system}.

\end{appendices}

%%%%%%%%%%%%%%%%%%%%%%%%%%%%%%%%%%%%%%%%%%%%%%%%%%%%%%%%%%%%%%%%%%%%%%%%%%%%%%%%%%%%%%%%%%%%%%%%%%%%
\bibliography{Bibliography}

\end{document}